\documentclass[10pt,oneside]{amsart} 
\usepackage{amsmath}
\usepackage{amsthm}
\usepackage{amsfonts}
\usepackage{amssymb,amscd,epsf,verbatim,mathtools,framed}
\usepackage{mathrsfs}
\usepackage{graphicx}
\usepackage{latexsym}
\usepackage{lscape}
\usepackage[colorlinks=true]{hyperref}
\hypersetup{colorlinks, citecolor=blue, filecolor=black, linkcolor=red, urlcolor=green}
\usepackage{epstopdf}
\usepackage{tikz}
\usetikzlibrary{calc}
\usetikzlibrary{matrix,arrows,decorations.pathmorphing}
\usepackage{tikz-cd}
\usepackage{color}
\usepackage{multirow}  
\usepackage{scalefnt}
\usepackage[shortlabels]{enumitem}

\newcommand{\cl}{\overline}

\newcommand{\Z}{\mathbb{Z}}
\newcommand{\F}{\mathbb{F}}
\newcommand{\N}{\mathbb{N}}

\newcommand{\Q}{\mathbb{Q}}

\renewcommand{\L}{\mathbb{L}}

\newcommand{\mC}{\mathcal{C}}
\newcommand{\mD}{\mathcal{D}}

\newcommand{\mF}{\mathcal{F}}

\newcommand{\mH}{\mathcal{H}}

\newcommand{\mM}{\mathcal{M}}

\newcommand{\mL}{\mathcal{L}}
\newcommand{\mO}{\mathcal{O}}

\newcommand{\mW}{\mathcal{W}}

\newcommand{\fm}{\mathfrak{m}} 
\newcommand{\fU}{\mathfrak{U}}
\newcommand{\fX}{\mathfrak{X}}

\newcommand{\ul}{\underline}

\newcommand{\ra}{\rightarrow}

\newcommand{\sq}{\widetilde}
\newcommand{\minus}{\backslash}
\DeclareMathOperator{\coker}{coker}
\DeclareMathOperator{\im}{im}
\DeclareMathOperator{\Hom}{Hom}

\DeclareMathOperator{\Ainf}{\text{A}_{\textup{inf}}}
\DeclareMathOperator{\AAinf}{\mathbb{A}_{\textup{inf}}}
\DeclareMathOperator{\Ainfx}{\mathbb{A}_{\textup{inf}, \fX_C}}
\DeclareMathOperator{\Acris}{\text{A}_{\textup{cris}}}

\DeclareMathOperator{\Bdr}{\text{B}_{\textup{dR}}}

\DeclareMathOperator{\Bst}{\text{B}_{\textup{st}}}
\DeclareMathOperator{\logcris}{log-cris}

\newcommand{\gr}[1]{\langle {#1} \rangle} 

\DeclareMathOperator{\spec}{Spec}
\DeclareMathOperator{\spf}{Spf}
 
\DeclareMathOperator{\etale}{\textup{\'etale}}
\DeclareMathOperator{\Etale}{\textup{\'Etale}}
\DeclareMathOperator{\ett}{\textup{\'et}}
\DeclareMathOperator{\pet}{\textup{pro\'et}}
\DeclareMathOperator{\Et}{\textup{\'Et}}

\DeclareMathOperator{\isom}{\;\xrightarrow{\: {}_{\sim} \:} \;}

\DeclareMathOperator{\uX}{\underline{X}}
\DeclareMathOperator{\uY}{\underline{Y}}
\DeclareMathOperator{\AD}{\textup{DA}^{\log}}
\DeclareMathOperator{\ADs}{\textup{DA}_{\mathfrak{sat}}^{\log}}
\DeclareMathOperator{\ADc}{\textup{DA}_{\mathfrak{str}}^{\log}}
\DeclareMathOperator{\ADp}{\textup{DA}^{\log,p}}
\DeclareMathOperator{\ADps}{\textup{DA}_{\mathfrak{sat}}^{\log,p}}
\DeclareMathOperator{\ADpc}{\textup{DA}_{\mathfrak{str}}^{\log,p}}
\DeclareMathOperator{\str}{\mathfrak{str}}
\DeclareMathOperator{\sat}{\mathfrak{sat}}

\newcommand{\bi}{\begin{itemize}}
\newcommand{\ei}{\end{itemize}}
\newcommand{\bt}{\begin{theorem}}
\newcommand{\et}{\end{theorem}}
\newcommand{\bbt}{\begin{theorem*}}
\newcommand{\eet}{\end{theorem*}}
\newcommand{\bp}{\begin{proposition}}
\newcommand{\ep}{\end{proposition}}
\newcommand{\bl}{\begin{lemma}}
\newcommand{\el}{\end{lemma}}
\newcommand{\bbl}{\begin{lemma*}}
\newcommand{\eel}{\end{lemma*}}
\newcommand{\bc}{\begin{corollary}}
\newcommand{\ec}{\end{corollary}}
\newcommand{\beg}{\begin{example}}
\newcommand{\eeg}{\end{example}}
\newcommand{\br}{\begin{remark}}
\newcommand{\er}{\end{remark}}
\newcommand{\bbr}{\begin{remark*}}
\newcommand{\eer}{\end{remark*}}
\newcommand{\bd}{\begin{definition}}
\newcommand{\ed}{\end{definition}}
\newcommand{\be}{\begin{enumerate}}
\newcommand{\ee}{\end{enumerate}}
\newcommand{\bex}{\begin{exercise}}
\newcommand{\eex}{\end{exercise}}
\newcommand{\bproof}{\begin{proof}}
\newcommand{\eproof}{\end{proof}}

\theoremstyle{definition}
\newtheorem{theorem}{Theorem}[section]

\newtheorem{mainthm}{Theorem} 
\newtheorem{lemma*}[mainthm]{Lemma}
\newtheorem{proposition*}[mainthm]{Proposition}
\newtheorem{maincor}[mainthm]{Corollary}

\theoremstyle{definition}
\newtheorem{definition}[theorem]{Definition}
\newtheorem{lemma}[theorem]{Lemma}
\newtheorem{corollary}[theorem]{Corollary}
\newtheorem{proposition}[theorem]{Proposition}

\newtheorem*{notation}{Notation}

\newtheorem{maindef}[mainthm]{Definition}
\newtheorem{remark}[theorem]{Remark}
\newtheorem{example}[theorem]{Example}
\newtheorem{remark*}[mainthm]{Remark}
\newtheorem{definition*}[mainthm]{Definition}

\theoremstyle{remark}
\newtheorem{sublemma}[theorem]{Sublemma}

\title[log de Rham--Witt]{Logarithmic de Rham--Witt complexes \\via the D\'ecalage operator}
 
\author{Zijian Yao}
\address{Department of Mathematics, Harvard University.}
\email{zijian.yao.math@gmail.com}

\begin{document}

\maketitle

\begin{abstract}
We provide a new formalism of de Rham--Witt complexes in the logarithmic setting. This construction generalizes a result of Bhatt--Lurie--Mathew, and agrees with those of Hyodo--Kato and Matsuue for log-smooth schemes of log-Cartier type. We then apply our formalism to obtain a more direct proof of the log crystalline comparison of $\Ainf$-cohomology in the case of semistable reduction, which is established by Cesnavicius--Koshiwara. 
\end{abstract}

\tableofcontents

\section{Introduction} 

Recently Bhatt--Lurie--Mathew \cite{BLM} gave a relatively elementary construction of the de Rham--Witt complex (hence the crystalline cohomology) for smooth varieties over a perfect field of characteristic $p$. The goal of this paper is to extend this construction to the logarithmic setting, and to apply this new formalism to study the $\Ainf$-cohomology (more precisely its log crystalline specialization) in the case of semistable reduction. 
 
\subsection{Main Results} \label{ss:main_results} Let $\ul X$ be a quasi-coherent log scheme over a perfect field $\ul k$ of characteristic $p$, and $\ul R = (R, M)$ be a local chart for its log structure.  

To contextualize our construction, we first define the category $\textup{DA}^{\log}$ of log Dieudonn\'e algebras. Roughly speaking, a log Dieudonn\'e algebra is a commutative differential graded algebra (cdga) $A^*$ together with a Frobenius $F$, a log structure $L$ and a log derivation $\delta$ in a compatible manner. In particular, the data are required to satisfy $d F = p F d$ and $\delta F = p F \delta$. 
We also need a certain full subcategory  $\textup{DA}^{\log}_{\mathfrak{str}}$ of \textit{strict} log Dieudonn\'e algebras, which essentially consists of log Dieudonn\'e algebras that admit the Verschiebung maps $V$ and are complete with respect to the $V$-filtration. We formulate this using the log version of d\'ecalage operator $\eta_p$, which on $p$-torsion free cochain complexes is given by 
$$\eta_p A^i :=\{x \in p^i A^i  | dx \in p^{i+1} A^{i+1}\}.$$ For each $p$-torsion free object ${A}^*$ in $\textup{DA}^{\log}$, there is a morphism $\phi_F: {A}^* \ra \eta_p (A^*)$ of log Dieudonn\'e algebras. The existence of $V$ is closely related to  $\phi_F$ being an isomorphism.  This allows us to give the following primitive\footnote{Strictly speaking, we need to restrict to subcategories  $\ADpc$ of $\ADc$, spanned by objects for which the log Frobenius is given by multiplication by $p$. We ignore this issue in the rest of the introduction for the sake of exposition.} 
definition.  

\begin{maindef} \label{maindef:de_Rham_Witt}
The saturated log de Rham--Witt complex $\mW \omega^*_{\ul{R}/\ul{k}}$ of $\ul R$ is a strict log Dieudonn\'e algebra,  such that, for any strict log Dieudonn\'e algebra $A^*$, the following canonical map is a bijection  $$\Hom_{\textup{DA}^{\log}_{\mathfrak{str}}} (\mW \omega^*_{\ul{R}/\ul{k}}, A^*) \isom \Hom_{\textup{Alg}_{\ul{k}}^{\log}}(\ul R, \ul A^0/ \im V).$$
\end{maindef}

 In other words, the functor sending $\ul R$ to $\mW \omega^*_{\ul R/\ul k}$ (if exists) provides a left adjoint of $\iota: \ADc \ra \textup{Alg}_{\ul{k}}$, which sends $A^*$ to $\ul A^0/ \im V.$  
Our first main result states that 

\begin{mainthm}
\label{mthm:log_dR_is_DA} 
$\mW \omega^*_{\ul{R}/\ul{k}}$ exists and glues to a sheaf $\mW \omega^*_{\ul X/\ul k}$ on the $\etale$ site of $X$. 
\end{mainthm} 

Our version of the log de Rham--Witt complexes agrees with the existing ones in \cite{HK}, \cite{LZ} and \cite{Ma} under the additional assumption that $\ul X$ is log-smooth of log-Cartier type.\footnote{In general, the saturated complex $\mW \omega^*_{\ul R/\ul k}$ differs from the crystalline construction of Hyodo--Kato. We will not touch upon this distinction in our article, which is already present in the case of ordinary schemes and has little to do with log structures. We instead refer the reader to Section 6.3 of \cite{BLM} where the authors compute an explicit example of a cusp.} We summarize this in the following
\begin{mainthm} 
Let  $\ul{X}$ be a log scheme over $\ul k$ which is log-smooth of log-Cartier type, then there exists isomorphisms of sheaves of complexes
$$ \mW \omega^*_{\ul X/\ul k} \isom W \Lambda^*_{\ul X/\ul k} \cong W^{\textup{HK}} \omega^*_{\ul X/\ul k},$$ 
compatible with Frobenius operators, where 
\be
\item[\textup{(Ma)}] $W \Lambda^*_{\ul X/\ul k} $ is the complex constructed by Matsuue \cite{Ma} as the initial object in the category $\mC_{F\!V}$ of ``log F-V procomplexes''.
\item[\textup{(HK)}] $W^{\textup{\tiny HK}}\omega^*_{\ul X/\ul k}$ is the log de Rham--Witt complexes constructed by Hyodo--Kato \cite{HK} using the log crystalline site.
\ee 
(cf. Section \ref{section:comparison}). In other words, in this case 
all three versions of ``log de Rham--Witt complexes'' agree,  and compute the log crystalline cohomology $R \Gamma_{\logcris} (\ul X/W(\ul k))$. 
\end{mainthm} 
 
One notable feature of this construction is that we start with the Frobenius and produce the Verschiebung and Restriction maps  ``along the way''. Another  feature is that our formulation makes essential use of the d\'ecalage operator $\eta_p$. To be more precise, we turn the isomorphism  
$$\phi_F: W \omega^*_{\ul X/\ul k} \cong \eta_p W \omega^*_{\ul X/\ul k}$$ which is typically an output
\footnote{This isomorphism is due to Illusie for (ordinary) de Rham--Witt complexes and to Hyodo--Kato in the  logarithmic setting. See Definition \ref{def:usual_DA_2} and Lemma \ref{lemma:HK_complex_is_ADc_2}.} 
of the construction of the log de Rham--Witt complex, as part of the input. As a result, our complex is characterized by a different universal property, which makes it easier to compare to cdga's equipped with Frobenius. This is essentially what makes our formalism suitable for proving the log crystalline comparison of $\Ainf$-cohomology. To formulate this application,  let $R\Gamma_{\Ainf} (\fX)$ be the $\Ainf$-cohomology of a formal scheme $\fX$ over $\mO_C$ with semistable reduction\footnote{More generally, we may assume that $\fX$ has generalized semistable reduction over $\mO_C$ as specified in Subsection \ref{ss:semistable_setup}.}, where $C$ is a completed algebraic closure of $W(k)[\frac{1}{p}]$. 
Write $\ul \fX$ for the formal log scheme with the divisorial log-structure from its mod $p$ fiber $\fX_{\mO_C/p}$.  

\begin{mainthm} \label{mthm:log_crys_compare_global} 
There is a Frobenius compatible quasi-isomorphism 
$$R \Gamma_{\Ainf} (\fX) \otimes^\L_{\Ainf}  W(k) \cong R \Gamma_{\log\textup{-cris}} (\ul{\fX_k}/W(k)),$$
relating the specialization of $R\Gamma_{\Ainf} (\fX)$ to $W(k)$ with the log crystalline cohomology of the special fiber of $\fX$ over $k$. 
\end{mainthm}

This result has recently been obtained by Cesnavicius and Koshiwara  \cite{Kestutis} from an ``absolute" crystalline comparison over $\Acris$.  In fact, $R\Gamma_{\Ainf} (\fX)$ also interpolates the $p$-adic $\etale$ cohomology and (log) de Rham cohomology of $\fX$, but the most delicate part is the log crystalline comparison, which we give a new proof of. In our approach, we avoid the comparison over the larger ring $\Acris$, and directly analyze the base change to $W(k)$. 
When $\fX$ has good reduction over $\mO_C$, the theorem is proven in \cite{BLM} as a major application of the saturated de Rham--Witt complexes without log structures. Our theorem is a generalization of their result. However, here in the semistable case, the nontrivial log structures impose additional difficulties, as it is not \textit{a priori} clear how to directly read off the logarithmic data from the $\Ainf$-cohomology.

\bbr In a forthcoming paper \cite{ZY} of the author, we give another proof of Theorem \ref{mthm:log_crys_compare_global},  by constructing a certain functorial map from the absolute log crystalline cohomology $R \Gamma_{\logcris} (\ul \fX_{\mO_C/p})$ to the $\Acris$-base change of $R\Gamma_{\Ainf} (\fX)$. We remark that our map goes in the opposite direction compared to the one used in \cite{Kestutis}. The existence of such a map and its compatibility with the $\Bdr$ comparison will be enough for us to recover most of the results in \cite{Kestutis}, including the $\Bst$ comparison theorem (in fact in a slightly more general setting). 
\eer

In the rest of the introduction we explain the idea of the construction as well as proofs of our main results in slightly more detail, highlighting new phenomenon that appears in the logarithmic setting.  

\subsection{The construction}  \label{ss:intro_dRW}
The construction of the saturated log de Rham--Witt complex is a direct generalization of \cite{BLM} and only involves elementary commutative algebra. Let us assume that $R$ is reduced so $W(R)$ is $p$-torsion free for simplicity. We first consider $W(\ul R) = (W(R), M)$, where the log structure is obtained from the Teichmuller lift. The log de Rham complex $\omega^*_{W(\ul{R})/\ul{W}}$ admits a canonical structure of a log Dieudonn\'e algebra.  We then construct a certain left adjoint 
$W(\:\:)_{\textup{sat}}:  \textup{DA}^{\log} \longrightarrow  \textup{DA}^{\log}_{\mathfrak{str}}$ of the inclusion functor $ \textup{DA}^{\log}_{\mathfrak{str}} \hookrightarrow  \textup{DA}^{\log}$, 
and define the saturated log de Rham--Witt complex of $\ul{R}/\ul{k}$ as 
$$\mW \omega^*_{\ul{R}/\ul{k}} := W \; \big( (\omega^*_{W(\ul{R})/\ul{W}})_{\mathfrak{sat}} \big). $$  
Finally, for a quasi-coherent log scheme $\ul X$ over $\ul k$, we show that the construction $\mW \omega^i_{-/\ul k}$ satisfies \'etale descent for each 
$i \ge 0$, and hence globalizes to a sheaf $\mW \omega_{\ul X /\ul k}^*$ of log Dieudonn\'e algebras on the $\etale$ site $X_{\ett}.$

\bbr 
Although it might be possible to develop the global theory directly by working with sheaves of log Dieudonn\'e algebras, we prefer to work locally, and then globalize only in the last step. This is especially convenient for our application to $\Ainf$-cohomology, where local charts are used in an essential way.
\eer  

We then prove the following comparison theorems, which ultimately amounts to the existence of the Cartier isomorphism in the logarithmic setting. 

\begin{mainthm} 
\label{mainthm:comparison_with_all}  

Suppose that $\ul{R}$ is log-smooth over $\ul k$ of log-Cartier type. Then $\mW \omega^*_{\ul{R}/\ul{k}}$ satisfies the following ``de Rham comparisons''.
 
\noindent $-$ {\textup{(mod $V$).}} There is a canonical isomorphism  of cochain complexes $$\omega_{\ul R/\ul k}^* \isom \mW_1 \omega_{\ul R/\ul k}^*.$$  

\noindent $-$ {\textup{(mod $p$).}} There is a canonical quasi-isomorphism
$$ \mW \omega^*_{\ul{R}/\ul{k}} / p \mW \omega^*_{\ul{R}/\ul{k}}  \isom \mW_1 \omega_{\ul R/\ul k}^* \cong  \omega^*_{\ul{R}/\ul{k}}.$$

\noindent  $-$ {\textup{(with Frobenius lifts).}} 
Suppose that there exists a log-Frobenius lift $(\ul{A}, \varphi)$ of $\ul{R}$ over $W(\ul{k})$ in the sense of Subsection \ref{sss:via_Frob_liftings}, then there is a quasi-isomorphism
$$\iota_{\varphi}: \widehat{\omega}^*_{\ul{A}/W(\ul{k})} \isom \mW \omega^*_{\ul{R}/\ul{k}}. $$
As notation suggests, the map $\iota_\varphi$ depends on the choice of $\varphi$ (compare with the canonical isomorphism in Proposition \ref{prop:log_dRW_another_construction}).
\end{mainthm}

\bbr These comparisons typically fail for general log rings $\ul R$ without any additional assumptions. For example, the $\F_p$-algebra $\mW_1 \omega^0_{\ul R/\ul k}$ is always reduced by Remark \ref{remark:completion_of_A_is_saturated}, but $\omega^0_{\ul R/\ul k} = R$ could be arbitrary. In particular, this implies the following interesting corollary. 
\begin{maincor} 
Let $\ul X$ be a fine log scheme over $\ul k$ that is log-smooth of log-Cartier type, then the underlying scheme $X$ is reduced. 
\end{maincor}
This recovers a result of Tsuji in \cite{Tsuji}. Note that, unlike the non-logarithmic case, a log-smooth scheme over a log point may fail to be reduced. For a simple example, let $n$ be a positive integer coprime to $p$. The map $$(\F_p, \N: 1 \mapsto 0) \longrightarrow (\F_p[x]/x^n, \N: 1 \mapsto x)$$ of fine log algebras is log-smooth, with the map of monoids $\N \ra \N$ given by $1 \mapsto n$. The ring $\F_p[x]/x^n$ is obviously non-reduced when $n \ge 2$. 
\eer

\bbr For a log scheme of ``generalized semistable type'' (in the sense of Subsection \ref{sec:revisit_HK}) over $(k, \N)$, our formalism also provides a convenient framework to construct the monodromy operator $N$ on its log crystalline cohomology. 
\eer 

\begin{proposition*}  
\label{thm:main_monodromy}
Let $\ul X/\ul k$ be a log scheme of generalized semistable type over the standard log point $\ul k = (k, \N)$.  Let $k^\circ = (k, 0)$ be the trivial log point. Denote the saturated log de Rham--Witt complex of $\ul X$ over the trivial log point by $\mW \sq \omega^*_{\ul X/\ul k} : = \mW \omega^*_{\ul X/k^\circ}.$ Then there is a short exact sequence of cochain complexes 
$$ 0 \ra \mW \omega^*_{\ul X/\ul k} [-1] \longrightarrow \mW \sq \omega^*_{\ul X/\ul k} \longrightarrow \mW \omega^*_{\ul X/\ul k} \ra 0.$$
The connecting homomorphism on cohomology 
$$N: H_{\log\textup{-cris}}^*(\ul X/W(\ul k)) = H^*(X, \mW \omega^*_{\ul X/\ul k}) \longrightarrow H_{\log\textup{-cris}}^*(\ul X/W(\ul k))$$
satisfies $N \varphi = p \varphi N$, where $\varphi$ is the functorial Frobenius on the log crystalline cohomology, and agrees with the monodromy operator of Hyodo--Kato. 
\end{proposition*}

\subsection{Log crystalline comparison of $\Ainf$-cohomology} \label{ss:intro_Ainf_coh}

In the rest of the introduction we explain the key ideas for proving Theorem \ref{mthm:log_crys_compare_global}. 
Let $\fX$ be a formal scheme as before  and $\fX_C$ be its adic generic fiber over $C$. Recall  from \cite{BMS} the definition of $A \Omega_{\fX}:= L \eta_{\mu} (R \nu_* \Ainfx)$, where $L \eta_\mu$ is the derived d\'ecalage operator 
with $\mu:= [\epsilon] - 1 \in \Ainf$ (cf. Section \ref{sec:Ainf_intro}), $\nu: (\fX_{C})_{\text{pro\'et}} \ra \fX_{\ett}$ is the nearby cycle map, and $\Ainfx=  W(\widehat \mO_{\fX_C}^{+, \flat})$ is the basic period sheaf  constructed in \cite{Scholze}.  Note that the sheaf $\Ainfx$ is equipped with a natural Frobenius, which by functoriality induces a Frobenius $\varphi$ on $A \Omega_{\fX}$ and the $\Ainf$-cohomology $R\Gamma_{\Ainf}(\fX) := R \Gamma (\fX_{\ett}, A\Omega_{\fX})$.

We pause to record an ingredient that goes into our proof. Consider the category $\textup{DC}$ of Dieudonn\'e complexes, 
consisting of cochain complexes equipped with Frobenius. As before we have a full subcategory $\textup{DC}_{\str}$ consisting of strict objects. It is shown in \cite{BLM} that the canonical functor 
$$ \textup{DC}_{\mathfrak{str}} \longrightarrow \widehat D(\Z_p)^{L\eta_p} $$
from strict Dieudonn\'e complexes to the $L\eta_p$-fixed points of the derived $p$-complete $\Z$-modules induces an equivalence of categories (cf. Theorem \ref{thm:BLM_L_eta_p_fixed_points}).  In particular, the $L \eta_p$ fixed points of $\widehat D(\Z_p)$ admit canonical cochain representatives.

We will prove Theorem \ref{mthm:log_crys_compare_global} by restricting to small enough affine opens $\spf S \in \fX_{\ett}$ which admit semistable coordinates $\square$ in the sense of Subsection \ref{ss:semistable_setup}. For such $S$, define 
$A \Omega_S:= R \Gamma (\spf S, A \Omega_{\fX}) \in D(\Ainf).$ The key local statement we want to prove is the following

\begin{mainthm} \label{mainthm:comparison_local} There is a canonical quasi-isomorphism 
$$ \mW \omega^*_{\ul{S_k}} \isom A \Omega_S \widehat \otimes^\L_{\Ainf} W(k)$$ 
compatible with Frobenius on both sides, where $\ul{S_k} = \ul S \otimes_{\mO_C} k$.  
\end{mainthm} 

To proceed we first construct a quasi-isomorphism 
$$ A \Omega_S \widehat \otimes^\L_{\Ainf} W(k) \cong L\eta_p (A \Omega_S \widehat \otimes^\L_{\Ainf} W(k)). $$
This puts $ A \Omega_S \widehat \otimes^\L_{\Ainf} W(k) $ in the category $ \widehat D(\Z_p)^{L\eta_p}$, and gives rise to a strict Dieudonn\'e algebra $A_S^*$. 
To upgrade $A_S^*$ to a log Dieudonn\'e algebra is more subtle compared to the non-logarithmic setting, since it is difficult to read off log structures from the derived category (for example we do not know a version of Theorem \ref{thm:BLM_L_eta_p_fixed_points} that takes log structures into account). To remedy this, we first fix a choice of semistable coordinates $\square$, from which we obtain a chart $(S, M^\square)$ for the log structure where $M^\square$ is particularly simple. Using explicit computations for this log algebra, we are able to define log structures and log derivations 
\begin{align*}
\alpha_r^\square: & M^\square \ra A_R^0/ \im(V^r + d V^r) \\ 
\delta_r^\square: & M^\square \ra  A_R^1/ \im(V^r + d V^r) 
\end{align*}
for each $r$, which fit together to give rise to logarithmic data $\alpha^\square$ and $\delta^\square$ on $A_S^*$.  This way, $A_S^*$ becomes a strict log Dieudonn\'e algebra. From the universal property of saturated log de Rham--Witt complexes, we are able to construct  
$$\tau^\square: \mW \omega^*_{\ul{S_k}/\ul k} \longrightarrow A_S^*,$$
which \textit{a priori} depends on the choice of coordinates. The upshot is that, even though the logarithmic data depend on coordinates, the morphism $\tau^\square$ on the underlying complex does not depend on this choice. In other words, we have 

\begin{proposition*}
As a morphism of Dieudonn\'e algebras, the map $\tau^\square$ is independent of the choice of semistable coordinates. In particular, we get a canonical map $\tau: \mW \omega^*_{\ul{S_k}/\ul k} \longrightarrow A_S^*$ of Dieudonn\'e algebras, which is an isomorphism.
\end{proposition*}

\bbr 
The final step in the argument replies on the Hodge-Tate comparison (Theorem \ref{thm:HT_dR_specialization}). In fact, in both this article and \cite{ZY}, we only need to analyze $A \Omega_S$ along the Hodge--Tate specialization $\Ainf \ra \mO_C$, instead of objects such as $A \Omega_S \otimes^\L A_{\textup{cris}}^{(m)}$ in \cite{Kestutis}. Moreover, our approach sees no difference in treating formal schemes which are of semistable type or generalized semistable type.   
\eer 

\subsection{Outline of the paper} 

Section \ref{sec:BLM} is largely expository, where we summarize relevant results from \cite{BLM}.  In Section \ref{sec:log_AD}, we define log Dieudonn\'e algebras, and its subcategories consisting of $p$-compatible strict objects. The saturated log de Rham--Witt complex lives in the latter subcategory. Section \ref{sec:log_dR_Witt}  consists of the technical core of the construction of $\mW \omega^*_{\ul R/\ul k}$ and its globalization to $X_{\ett}$. In this section we also compare $\mW \omega^*_{\ul R/\ul k}$ to the (completed) de Rham complex of a log-smooth ``log-Frobenius lifting''.   
In section \ref{section:comparison}, we show that our construction agrees with the existing ones due to Hyodo--Kato and Matsuue in the case when the log scheme is sufficiently smooth. We then reconstruct the monodromy operator on the log crystalline cohomology for semistable log schemes. Finally, in the last two sections of the article we discuss the log crystalline comparison of $\Ainf$-cohomology. We recall the definition of $\Ainf$-cohomology with toy examples in Section \ref{sec:Ainf_intro} and give our proof of the log crystalline comparison in Section \ref{section:log crystalline_specialization}. We also include Appendix \ref{sec:log} to briefly review the necessary background on log schemes. 

\subsection{Conventions}  \label{ss:conventions}
We fix a prime $p$ once and for all. We say that a cochain complex $M^*$ is $p$-torsion free if each $M^i$ is $p$-torsion free.  By a cdga we mean a commutative differential graded aglebra $(A^* = \oplus_{i \ge 0} A^i, d)$. In particular, the differential operator $d$ increases grading by $+1$ and satisfies $\textup{d}(ab) = (\textup{d}a) b + (-1)^{k} a (\textup{d}b)$ for $a \in A^k$; commutativity requires that $ab = (-1)^{kl} ba$ for $a \in A^k, b \in A^l$ and that $a^2 = 0$ for all $a \in A^{2j+1}$. The latter condition is redundant unless $p = 2$. 

For the convention on log schemes we refer to Appendix \ref{sec:log}, where we mostly follow \cite{Kato} except that we denote a log scheme $(X, M_X)$ by $\ul X$.  In addition, by a log algebra $(R, L)$ we mean an algebra $R$ together with a monoid morphism $L \ra R$, this is the same data as giving a pre-log scheme $\ul Y = (\spec R, L)$ with constant pre-log structure $\beta: L_Y \ra \mO_Y$. We often denote by $\mL^a$ the associated log structure of $L$ on $Y = \spec R$, which is the pushout of $\beta^{-1} (\mO_Y^\times) \ra L_Y$ along $\beta^{-1} (\mO_Y^\times) \hookrightarrow \mO_Y$ on the $\etale$ site of $Y$, and denote its global section by $L^{\textup{sh}} = \Gamma(Y, \mL^a)$. 

Throughout the article, we use $k$ to denote a perfect field of characteristic $p$, which is further required to be algebraically closed in Section  \ref{sec:Ainf_intro} and \ref{section:log crystalline_specialization} for simplicity. Moreover, we denote by $\ul{W}$ a log algebra $(W(k), N)$ where $N$ could be arbitrary. We reserve the notation $W(\ul k)$ for $(W(k), [\alpha]: N \ra W(k))$, where the log structure comes from the Teichmuller lift of a log point $\ul k = (k, N)$. The reader is welcome to take $\ul W = W(\ul k)$ for convenience.  

Finally, in Section  \ref{sec:Ainf_intro} and \ref{section:log crystalline_specialization}, we use the language of derived $\infty$-categories. We mostly follow the conventions of Section 10 of \cite{BLM}, in particular, for a commutative ring $A$, we use $\mD(A)$ to denote the derived $\infty$-category of $A$-modules and identify the usual (triangulated) derived category $D(A)$ with the homotopy category of $\mD(A)$. 

\subsection*{Acknowledgements}
I am heartily grateful to Luc Illusie for his interests and various helpful correspondences throughout the project. I would like to thank Lin Chen, Mark Kisin, Matthew Morrow, Sasha Petrov, Lynnelle Ye and Allen Yuan for their help and suggestions while I was preparing the preprint. Additionally, as it will be clear to the reader, I am greatly indebted to Bhargav Bhatt, Jacob Lurie and Akhil Mathew for their work \cite{BLM}. Finally, this paper grows out of a seminar at Harvard in summer 2017, as an attempt to understand the work of Bhatt--Morrow--Scholze \cite{BMS}. I would like to thank all the participants for numerous helpful conversations or correspondences.

\newpage 

 
\section{Saturated de Rham--Witt complexes} \label{sec:BLM}

In this section, we summarize the relevant results of Bhatt--Lurie--Mathew \cite{BLM}, which we generalize in the rest of the article.  
 
\subsection{Dieudonn\'e algebras, saturation and $V$-completion} 

\subsubsection{The operator $\eta$}  

We first recall the d\'ecalage operator 
\footnote{The d\'ecalage operator was first introduced in \cite{BO}. It was later used in \cite{RI} for the crystalline construction of the de Rham--Witt complex $W \Omega^*_{X/k}$ (and similarly in \cite{HK} to construct the log de Rham--Witt complexes from the log crystalline site). Recently it appeared in \cite{BMS} to define $\Ainf$-cohomology.}.  

\bd \label{def:decalage}
Let $R$ be a ring and $\mu \in R$ a nonzero divisor. Let $(M^{\ast}, d)$ be a cochain complex of $\mu$-torsion free $R$-modules, then $(\eta_\mu M)^{\ast} \subset M^*[\frac{1}{\mu}]$ is defined to be the sub-complex given by 
$$(\eta_\mu M)^i = \{ x \in \mu^i M^i : dx \in \mu^{i+1} M^{i+1} \}.$$
\ed

Note that $\eta_\mu$ kills $\mu$-torsion in the cohomology, more precisely, we have $H^i(\eta_\mu M^*) \cong H^i(M^*)/H^i(M^*)[\mu]$. In particular, $\eta_\mu$ descends to the derived category, inducing a (non-exact) functor $L\eta_\mu: D(R) \ra D(R)$.  For other properties of $\eta_\mu$ and $L\eta_\mu$ we refer the reader to $\S 6$ of \cite{BMS}.  
 
\subsubsection{Dieudonn\'e algebras and saturated Dieudonn\'e algebras} \indent

\bd \label{def:usual_DA} 
A Dieudonn\'e algebra is a triple $(A^{\ast}  = \underset{i \ge 0}{\bigoplus} A^i, d, F)$ where $(A^{\ast}, d)$ is a cdga, $F: A^{\ast} \ra A^{\ast}$ is a graded algebra map such that 
\bi
\item $F (x) \equiv x^p \mod p  $ for all $x \in A^0$; 
\item $d F (x) = p F (dx)$ for all $x \in A^*$. 
\ei
\ed 

\bd \label{def:usual_DA_2} 
Let $A^{\ast} =(A^{\ast}, d, F)$ be a  $p$-torsion free Dieudonn\'e algebra, then $F$ determines a map of cochain complexes
$$\phi_F: A^{\ast} \ra \eta_p A^{\ast}$$ by sending $x \mapsto p^n F(x)$ for $x \in A^n$. A Dieudonn\'e algebra $A^{\ast}$ is saturated if it is $p$-torsion free  and $\phi_F$ is an isomorphism. 
\ed

Morphisms between two Dieudonn\'e algebras $(A^{\ast}, d, F)$ and $(B^{\ast}, d', F')$ is a morphism $f: A^{\ast} \ra B^{\ast}$ between cdgas compatible with the Frobenius maps $F$ and $F'$.  The category of Dieudonn\'e algebras (resp. the full subcategory spanned by saturated algebras) is denoted by $\textup{DA}$ (resp. $\textup{DA}_{\mathfrak{sat}}$)

\subsubsection{Verschiebung} \label{sss:V} 
Consider $A^{\ast} \in \textup{DA}_{\mathfrak{sat}}$. For each $i \in \Z$, the map 
$$\phi_F: A^i  \xrightarrow{\: \: F \: \: } \{x \in A^i: dx \in p A^{i+1}\} \xrightarrow{ \times p^i} (\eta_p A)^i$$ is an isomorphism, hence $F$ is injective and $F(A^{\ast})$ contains $p A^{\ast}$. Therefore, for each $x \in A^n$, there is a unique element $V x$ such that $F(Vx) = px$. It is straightforward to check that 
$$ F V = V F = p, \quad F  d   V = d, \quad V  d = p \; d V, \quad x Vy = V(Fx \cdot y).$$

\subsubsection{Saturation} \label{sss:phi_F_in_AD} 

Suppose $A^{\ast} \in \text{AD}$ is $p$-torsion free, then $(\eta_p A)^{\ast}$ with its inherited differential and Frobenius structure is again a Dieudonn\'e algebra. The only thing to check is that for any $x \in (\eta_p A)^0$, $F(x) = x^p + p y$ for some $y \in  (\eta_p A)^0$ (not just in $A^0$). The map $\phi_F: A^\ast \ra (\eta_p A)^{\ast}$ is a morphism of Dieudonn\'e algebras.  The inclusion functor $\text{AD}_{\mathfrak{sat}} \hookrightarrow \text{AD}$ admits a left adjoint, $A^{\ast} \mapsto A_{\mathfrak{sat}}^{\ast} $, which we call the \textbf{saturation} of $A^*$, given as follows. We replace $A^{\ast} \in \text{AD}$ by its $p$-torsion free quotient if necessary and suppose that $A^*$ is $p$-torsion free, and then dfine $ A_{\mathfrak{sat}}^{\ast} $ to be the direct limit of
$$ A^{\ast} \xrightarrow{\phi_F} (\eta_p A)^{\ast}  \xrightarrow{\eta_p(\phi_F)} (\eta_p \eta_p A)^{\ast}  \xrightarrow{\eta_p^2(\phi_F)} (\eta_p^3 A)^{\ast} \ra \cdots $$
By the discussion above, the saturation $A^{\ast}_{\mathfrak{sat}}$ inherits the structure of a Dieudonn\'e algebra, and the natural map $A^{\ast} \ra A^{\ast}_{\mathfrak{sat}}$ is a morphism in $\text{AD}$. 

\subsubsection{Strict Dieudonn\'e algebras and $V$-completion} \label{sss:W_r_of_AD}

Let $A^{\ast}$ be a saturated Dieudonn\'e algebra.  For each $r \ge 1$, we form the quotient
$$W_r (A^{\ast}):= A^{\ast}/(V^r A^{\ast} + d V^{r} A^{\ast}),$$ which is a cdga as $V^r A^{\ast} + d V^{r} A^{\ast}$ is a differential graded ideal.  Next we define the $V$-completion of $A^{\ast}$ to be the limit $W (A^{\ast}) := \varprojlim W_r A^{\ast}$ along the natural projection maps $R: W_r A^{\ast} \ra W_{r-1} A^{\ast}$.  It is easy to check that the map $F$ on  $A^{\ast}$ induces $F: W_r (A^{\ast}) \ra W_{r-1} (A^{\ast})$ on the quotients. Similarly, we have $V: W_r (A^{\ast}) \ra W_{r+1} (A^{\ast})$. 

\bd A saturated Dieudonn\'e algebra $A^{\ast}$ is strict if the canonical map 
$ A^{\ast} \longrightarrow W(A^{\ast})$ is an isomorphism.  The full subcategory of $\textup{DA}_{\mathfrak{sat}}$ spanned by strict algebras is denoted by $\textup{DA}_{\mathfrak{str}}$. 
\ed 

\br \label{remark:completion_of_A_is_saturated} Let $A^*$ be a saturated Dieudonn\'e algebra as above. 
\bi
\item $A^0/V A^0$ is a reduced $\F_p$-algebra by Lemma 3.6.1 of \cite{BLM}. 
\item $W(A^*)$ in the definition above becomes a Dieudonn\'e algebra with $F$ the inverse limit of $F:  W_r (A^{\ast}) \ra W_{r-1} (A^{\ast})$. It is in fact still saturated. One needs to check that (i). $W(A^*)$ is $p$-torsion free, (ii). if $x \in A^i :=  (W(A^{\ast}))^i$ is an element with $d x \in p W^{i+1}$, then $x  \in \im (F)$, and (iii). the inverse limit $F$ on $W(A^*)$ satisfies $F(x) \equiv x^p \mod p$, for the detail we refer to Section 2.6 and 3.5 of \cite{BLM}. 
\ei
\er

\br \label{remark:W_r_induces_isom}
The $V$-completion $W(A^*)$ for $A^* \in \textup{DA}_{\mathfrak{sat}}$ is strict. More precisely, for each $r \ge 1$, the canonical map $W_r (A^*) \ra W_r (W(A^*))$ is an isomorphism (by \cite{BLM} Proposition 2.7.5). The completion functor $A^{\ast} \mapsto W(A^{\ast})$ provides a left adjoint of the inclusion $\textup{DA}_{\mathfrak{str}} \hookrightarrow \textup{DA}_{\mathfrak{sat}}$. 
\er 

\br \label{lemma:A_0_is_Frob} 
For a strict Dieudonn\'e algebra $A^* \in \textup{DA}_{\mathfrak{str}}$, each $A^i$ is $p$-adically complete. Moreover, $A^0$ is the ring of Witt-vectors of $W_1(A)^0 = A^0/V A^0$. In other words, there is a unique isomorphism  $\mu: A^0 \ra W(A^0/VA^0)$ such that the composition $A^0 \ra W(A^0/VA^0) \xrightarrow{\textup{can}} A^0/V A^0$ is the projection map. Under the map $\mu$ the Frobenius $F$ on $A^0$ corresponds to the usual Frobenius of Witt vectors. 
\er
 
\subsection{The Cartier criterion}  \label{ss:Cartier_criterion}

Let $A \in \textup{DA}$ be a $p$-torsion free Dieudonn\'e algbera. Consider the cochain complex $(H^*(A^*/p A^*), \beta)$ where $\beta$ is the Bockstein differential induced from $0 \ra A^*/p \xrightarrow{p} A^*/p^2 \ra A^*/p \ra 0$. Then we have the following commutative diagrams where each arrow is a map of cochain complexes: 
\[
\begin{tikzcd}[column sep=0.2em]
 (A^*/p A^*, d)  \arrow[rr, "F"]  \arrow[rd, "\phi_F"] & &  (H^*(A^*/p A^*), \beta) \\
 & (\eta_p (A^*)/p, d)  \arrow[ru, "\gamma"]
\end{tikzcd}
\] 
and $\gamma$ is defined by sending $x \in (\eta_p A)^i$ to $x/p^i$. Suppose that $A^*$ is in addition saturated, then the composition $F$ in the diagram factors through the quotient cochain complex $W_1(A^*)$, as in the diagram below: 
\[
\begin{tikzcd}[column sep=0.2em]
& (W_1 (A^*), d) \arrow[rd, "F_1"]\\
 (A^*/p A^*, d) \arrow[ru, "q"]  \arrow[rr, "F"] & &  (H^*(A^*/p A^*), \beta) 
\end{tikzcd}
\] 

\bl \label{lemma:F_1_is_isom_if_saturated} Let $A^*$ be a $p$-torsion free Dieudonn\'e algebra, then 
\be
\item $\gamma: \eta_p(A^*)/p \ra H^*(A^*/p)$ is a quasi-isomorphism of cochain complexes; 
\item Suppose that $A^*$ is saturated, then $F_1: W_1(A^*) \ra H^*(A^*/p)$ is an isomorphism of cochain complexes. 
\ee 
\el

\bproof  (1).  The map $\gamma$ is clearly surjective, its kernel is acyclic by Proposition 2.4.4 of \cite{BLM} (also compare to Proposition 1.3.4 of \cite{Deligne}). (2). This is a special case of Proposition 2.7.1 of \cite{BLM}.
\eproof 

\bd \label{def:Cartier_criterion} We say that a Dieudonn\'e algebra $A^*  \in \textup{DA}$ satisfies the Cartier criterion (or is of Cartier type) if 
$$F: (A^*/p, d) \ra (H^*(A^*/p), \beta)$$ is an isomorphism of cochain complexes. 
\ed 

In the case in the above definition, we denote the map $F$ by $C^{-1}$. The following is immediate.  
\bc \label{cor:Cartier_implies_qi} Let $A^* \in \textup{DA}$. 
\be
\item If $A^*$ satisfies the Cartier criterion, then $\phi_F: A^*/p \ra \eta_p(A^*)/p$ is a quasi-isomorphism, therefore we have a quasi-isomorphism 
$$ A^{\ast}/p A^{\ast} \longrightarrow A_{\mathfrak{sat}}^{\ast}/p A_{\mathfrak{sat}}^{\ast} $$
\item If $A^*$ is saturated, then $q: A^*/p \ra W_1(A^*)$ is a quasi-isomorphism of cochain complexes. 
\item  Let $f: A^* \ra B^*$ be a morphism of saturated Dieudonn\'e algebras, then $f: A^*/p \ra B^*/p$ is a quasi-isomorphism if and only if $f: W_1 (A^*) \ra W_1(B^*)$ is an isomorphism, if and only if $f: W (A^*) \ra W (B^*)$ is an isomorphism. 
\ee 
\ec
\bproof 
The last part follows from Corollary 2.7.4 in \cite{BLM}.
\eproof 

\bc \label{cor:Cartier_implies_qi_2}
Let $A^*$ be a $p$-complete and $p$-torsion free Dieudonn\'e algebra satisfying the ``Cartier criterion'', then the canonical map 
$$A^{\ast}  \longrightarrow W(A_{\mathfrak{sat}}^{\ast})$$
is a quasi-isomorphism. 
\ec

\bproof 
It suffices to show that $A^*/p \ra A_{\mathfrak{sat}}^{\ast}/p \ra W(A_{\mathfrak{sat}}^{\ast})/p$ is a quasi-isomorphism since each $A^i$ is p-adically complete and $p$-torsion free. The first map is a quasi-isomorphism by Corollary \ref{cor:Cartier_implies_qi} (1), and the second map is a quasi-isomorphism by  Remark \ref{remark:W_r_induces_isom} and  Corollary \ref{cor:Cartier_implies_qi} (3).
\eproof 

\subsection{The saturated de Rham--Witt complexes} 

Given the preparations above, the saturated de Rham--Witt complexes $\mW \Omega_R^{*}$ for an $\F_p$ algebra $R$ is constructed in \cite{BLM} as follows: first equip $\Omega_{W(R_{\textup{red}})}^*$ with a Frobenius map $F$ under which it becomes a Dieudonn\'e algebra, then define $\mW\Omega_R^* := W(\Omega_{W(R_{\textup{red}})}^*)_{\mathfrak{sat}}$ (resp. define $\mW_n \Omega^*_{R} = W_n(\Omega_{W(R_{\textup{red}})}^*)_{\mathfrak{sat}}$). The construction $R \mapsto \mW \Omega_R^*$ is clearly functorial and provides a left adjoint of the functor $\textup{DA}_{\mathfrak{str}} \longrightarrow \textup{Alg}_{\F_p} $ that sends  $A^{\ast} \mapsto A^0/V A^0.$
 
We list some basic properties of $ \mW \Omega_R^*$ proved in \cite{BLM}. Suppose there exists a $p$-torsion free algebra $A$ lifting $R$, equipped with an endomorphism $\varphi: A \ra A$ satisfying $\varphi(a) \equiv a^p \mod p$. Then there exists a natural map 
$$\mu: \widehat \Omega_{A}^* \ra \mW \Omega^*_{R}.$$ The map $\mu$ induces an isomorphism $W(\widehat \Omega_{A}^*)_{\mathfrak{sat}} \isom \mW \Omega_{R}^*$. Suppose in addition that $R$ is smooth over a perfect extension of $\F_p$, then $\mu$ is a quasi-isomorphism.  Now further suppose that $R$ is a smooth $\F_p$-algebra, then there exists a canonical isomorphism 
$$\nu: \Omega^*_{R} \isom \mW_1 \Omega^*_{R}.$$   
Finally, if we write $W \Omega_{R}^*$ for the de Rham--Witt complexes constructed by Illusie in \cite{Illusie}, then there exists a natural isomorphism 
$$\gamma: W \Omega^*_R \isom \mW \Omega_R^*.$$

\subsection{$L\eta_p$-fixed points of the $p$-complete derived category}  \label{ss:L_eta_fixed_point}

The goal of this subsection is to state Theorem \ref{thm:BLM_L_eta_p_fixed_points}, which is crucial for our application to $\Ainf$-cohomology in Section \ref{section:log crystalline_specialization}. 

To start, we first relax the notion of Dieudonn\'e algebras to Dieudonn\'e complexes.  A Dieudonn\'e complex $(M^*, d, F)$ is a cochain complex $(M^*, d)$ of abelian groups equipped with a Frobenius operation $F: M^* \ra M^*$ satisfying $dF = p Fd$. They form a category $\textup{DC}$, where morphisms are maps of cochain complexes compatible with the Frobenius structure. The full subcategory $\textup{DC}_{\mathfrak{sat}}$ (resp.  $\textup{DC}_{\mathfrak{str}}$) of saturated (resp. strict) Dieudonn\'e complexes is similarly defined, as well as the saturation (resp. $V$-completion) functor, for detail we refer the interested readers to section 2 of \cite{BLM}.  Let $\widehat{D} (\Z)_p$ be the full subcategory of $D(\Z)$ (viewed as a triangulated category) generated by derived $p$-adically complete objects. For applications later, we also need the derived $\infty$-category enhancement $\mD (\Z)$ of its homotopy category $D(\Z)$, and similarly define $\widehat \mD (\Z)_p \subset \mD(\Z)$. Denote by $\widehat \mD (\Z)_p^{L\eta_p}$ (resp. $\widehat D(\Z)_p^{L\eta_p}$) the $L\eta_p$ fixed points of $\widehat \mD (\Z)_p$ (resp. $\widehat D (\Z)_p$), which consists of an object $X$ in the category equipped with an isomorphism $X \isom L\eta_p X$.  From definition, there is a natural functor 
$$ \textup{DC}_{\mathfrak{str}} \longrightarrow \widehat D (\Z)_p^{L\eta_p} $$

\bt [ \cite{BLM} Theorem 7.3.4 and 7.4.8]\label{thm:BLM_L_eta_p_fixed_points}\footnote{This equivalence is already hinted in \cite{RI} and \cite{HK} for the log version. For example, compare the inverse functor described below and the Hyodo--Kato construction in Subsection \ref{sss:construction_HK}.}  
The functor above factors through the following sequences of equivalences of categories: 
$$\textup{DC}_{\mathfrak{str}}  \isom  \widehat \mD (\Z)_p^{L\eta_p}  \isom  \widehat D (\Z)_p^{L\eta_p}. $$
\et 
We refer the interested readers to Section 7 of \cite{BLM} for more detail. 

\br \label{remark:BLM_L_eta_p_fixed_points} For the application later, we describe a quasi-inverse functor (which is not explicitly described in \cite{BLM} but not difficult to obtain). Let 
$$(X, \psi: X \isom L\eta_p X) \in \widehat{D} (\Z)_p^{L\eta_p}$$ be an $L \eta_p$ fixed point of the category $\widehat{D(\Z)}_p$. For each $r \ge 1$, we define the cochain complexes $(X_n^*, \beta_n)$ by 
$$X_r^i := H^i (X \otimes^{\L}\Z/p^r),$$ where the differentials $\beta_r$ are given by the Bockstein differentials (associated to $X \otimes^{\L}Z/p^r \ra X \otimes^{\L}Z/p^{2r} \ra X \otimes^{\L}Z/p^r$). Next define maps 
$$\mu_r: H^k (X\otimes^{\L}\Z/p^{r})  \longrightarrow H^k ( L \eta_p X\otimes^{\L}\Z/p^{r-1})$$ by sending $y \in H^k (X\otimes^{\L}\Z/p^{r})$ to $p^k y \in H^k(L\eta_p X\otimes^{\L}\Z/p^{r-1})$. This leads to the ``restriction maps'' $R_r: = X_{r}^* \longrightarrow X_{r-1}^* $ by setting $R_r := \psi^{-1} \circ \mu_r$, namely the composition  
\[\begin{tikzcd}
 H^*(X\otimes^{\L}\Z/p^{r}) \arrow[r, "\mu_r"] \arrow[rd, swap, "R_r"] & H^*(L\eta_p X\otimes^{\L}\Z/p^{r-1}) \arrow[d, "\psi^{-1}"]  \\
 &   H^*( X\otimes^{\L}\Z/p^{r-1})
\end{tikzcd} \]
We build a cochain complex  $X^* := \varprojlim_{R_r} X_r^*$ by taking the inverse limit along the restriction maps $R_r$. Equipped with Bockstein differentials, $X^*$ becomes a Dieudonn\'e complex as follows: the canonical projection $\Z/p^{r} \ra \Z/p^{r-1}$ and map $\Z/p^{r-1} \xrightarrow{\times p} \Z/p^{r}$ respectively induce the maps 
$$F: X_{r}^* \ra X_{r-1}^*, \qquad  V: X_{r-1}^* \ra X_{r}^*.$$
They are compatible with $R_r$ and therefore induce operators $F, V$ on $X^*$. 
\er

\newpage 

\section{Log Dieudonn\'e algebra} \indent  \label{sec:log_AD} 

In this section we describe an enhancement of Dieudonn\'e algebras that incorporates log structures, which serves as a basis for our theory of saturated log de Rham--Witt complexes on local charts. The definitions in this section are mostly analogous to those in \cite{BLM}, with minor differences. In particular, we prove that, for a strict Dieudonn\'e algebra with a $p$-compatible log structure, its log structure is ``valued in'' Teichmuller representatives (Proposition \ref{prop:image_of_log}). For notations on log geometry we refer to Subsection \ref{ss:conventions} and Appendix \ref{sec:log}. In this section and in the construction of saturated log de Rham--Witt complexes, we do not require the log structures to satisfy additional properties such as being integral or coherent. 
  
\subsection{Log Dieudonn\'e algebra}  
\bd \label{def:log_DA} Let $\ul W = (W(k), N)$ be a log algebra. A log Dieudonn\'e $\ul W$-algebra is a tuple $(A^*, L, d, \delta, F, F_L),$ consisting of the following data: 
\bi
\item $(A^{\ast} = \oplus_{i \ge 0} A^i, d)$ is a cdga over $W$, 
\item $\ul{A}^0 = (A^0, L)$ is a log algebra over $\ul W$, 
\item $\delta: L \ra A^1$ is a map of monoids, and 
\item $F: A^* \ra A^*$ a graded algebra homomorphism,
\ei  satisfying the following requirements:
\be
\item $(A^*, d, F)$ is a Dieudonn\'e algebra in the sense of Definition \ref{def:usual_DA}. 
\item $(d: A^0 \ra A^1, \delta: L \ra A^1)$ is a log derivation of $\ul{A}^0/\ul{W}$,  where we further require that the composition  $  L \xrightarrow{\delta} A^1 \xrightarrow{d} A^2 $ is $0$. 
\item $\delta$ is compatible with $F: A^1 \ra A^1$, in other words it satisfies 
$$\delta F_L = p F \delta.$$
\item $F_L: L \ra L $ is a monoid homomorphism satisfying 
$$F \circ \alpha = \alpha \circ F_L .$$
\ee    
\ed

\begin{notation} In what follows we often suppress notations and just write $A^*$  for a log Dieudonn\'e algebra. 
\end{notation}

\br
In contrast to the ordinary (= non-logarithmic) setting, for the definition of log Dieudonn\'e algebras we need to specify a base $\ul W$ to start with. In the ordinary case, $\Omega_{k/\F_p}  = 0$ (and $\Omega_{W(k)/\Z_p} = 0$), so the base is to some extent irrelevant (as one can see from the construction in Subsection \ref{ss:construction_via_Frob}, for example). If we have two distinct log structures $N_1, N_2$ on $k$, on the other hand, $\Omega_{(k, N_1)/(k, N_2)}$ is typically nonzero. Therefore it is necessary to specify the base $\ul W$ for log Dieudonn\'e algebras, in particular its log structure. 
\er 

\bd \label{def:log_DA_morphism} 
A morphism between two log Dieudonn\'e $\ul W$-algebras ${A}^*$ and ${B}^*$ is a pair $(f, \psi)$, where $f: (A^*, d, F) \ra (B^*, d', F')$ is a morphism between Dieudonna algebras over $W$, and $\psi: L_A \ra L_B$ is a monoid morphism over $N$ which is compatible with $F_L$, $\alpha$ and $\delta$. More precisely, 
\be
\item $f \circ d = d' \circ f$ and $f \circ F = F' \circ f$ on each $A^i$.  
\item The structure map $N \ra L_B$ agrees with $N \ra L_A \xrightarrow{\psi} L_B$ (resp. $W \ra B^0$ agrees with $W \ra A^0 \xrightarrow{f|_{A^0}} B^0$). Moreover $f \circ \alpha_A = \alpha_B \circ \psi$ on $L_A$. This is simply saying that $(f, \psi)|_{(A^0, L_A)}$ is a morphism of log $\ul W$-algebras. 
\item $\psi$ satisfies $\psi \circ F_{L_A} = F_{L_B} \circ \psi$,  and $f \circ \delta_A = \delta_B \circ \psi$. 
\ee
\ed

\begin{notation}
We denote the category of log Dieudonn\'e algebras by $\textup{DA}_{/ \ul W}^{\log}$, and will often suppress notations to write $\textup{DA}^{\log} =  \textup{DA}_{/ \ul W}^{\log}$ when $\ul W$ is understood. 
\end{notation}

\bl  \label{lemma:isom_in_AD_log}
A morphism $(f, \psi): A^* \ra B^*$ of log Dieudonn\'e algebras over $\ul W$ is an isomorphism if and only if both $f$ and $\psi$ are isomorphisms. 
\el

\bproof 
Suppose that $f$ has an inverse $g: B^* \ra A^*$ as Dieudonn\'e algebras over $W$ and that $\psi$ has an inverse $\chi: L_B \ra L_A$, the lemma claims that $\chi$ satisfies the compatibilities (2) and (3) in Definition \ref{def:log_DA_morphism}. This is straight-forward. For example, to check (2) we observe that for any $m \in L_B$
$$g \circ \alpha_B (m) = g (\alpha_B (\psi \circ \chi (m))) = g (f(\alpha_A ( \chi m))) = \alpha_A \circ \chi (m).$$ 
The remaining relations are similarly verified. 
\eproof

\subsection{$p$-compatible log Dieudonn\'e algebras} 
 
 \bd
 We say that a log Dieudonn\'e algebra  $A^*$  over $\ul W$ is $p$-compatible if the Frobenius on the log structure $L$ is given by $F_L = p$. Let  $\ADp$ denote the full subcategory of $\AD$ spanned by $p$-compatible objects. 
 \ed 

\br \label{remark:why_taking_honest_subcategory}  
In this article we are mostly interested in $p$-compatible log Dieudonn\'e algebras (see Proposition \ref{thm:omega_of_Witt_lifting_AD}). The $p$-compatibility condition imposes some restrictions on what log structures are allowed. For example, if $\fX$ is a scheme over $\Z_p$ and $\mL$ is a log structure on $\fX$, then for an affine open $U = \spec A \subset \fX$, the log algebra $(A, \Gamma (U, \mL))$ is usually not $p$-compatible. A typical example that arises is $(\Z_p[x, y]/xy, L = \N^2)$ with the pre-log structure $(a, b) \mapsto x^a y^b$ (see Example \ref{sss:examples}) and Frobenius $F_{L}$ given by multiplication by $p$, while on the (sections of the) associated log structure $L^a = \Z_p^\times \oplus \N^2$ the induced Frobenius is $F_{L^a} = \textup{id} \oplus F_L$. For this reason we are usually only work with pre-log structures on schemes. 
\er

\subsection{Saturated log Dieudonn\'e algebras} In this subsection and next we discuss saturation and $V$-completion in the logarithmic setting. 
 
\bd \label{def:log_decalage} Let $(A^*,L \xrightarrow{\alpha} A^0, d, \delta, F, F_L)$ be a log Dieudonn\'e algebra where $A^*$ is $p$-torsion free. We define the log Dieudonn\'e algebra $\eta_p ({A}^*)$ to be the tuple 
$$ \Big(\eta_p (A^*), L \xrightarrow{\eta_p (\alpha)} \eta_p A^0, d, \eta_p(\delta), F, F_L \Big) \in \AD$$
where 
\bi 
\item $\eta_p A^{i} = \{ x \in p^i A^i : d x \in p^{i+1} A^{i + 1} \}$ as in Definition \ref{def:decalage} and is equipped with $\phi_F: A^* \ra \eta_p A^*$, 
\item the log algebra $(\eta_p A^0, L)$ is given by $\eta_p(\alpha): L \xrightarrow{F_L} F_L (L) \xrightarrow{\alpha} \eta_p A^0$,
\item the log derivation $\eta_p (\delta)$ is given by $\eta_p (\delta): L \xrightarrow{F_L} F_L (L) \xrightarrow{\delta} \eta_p A^1$,  \footnote{Equivalently, $\eta_p (\delta)$ is given by the composition $L \xrightarrow{\delta} A^1 \xrightarrow{p F} \eta_p A^1$. Likewise, $\eta_p (\alpha)$ is given by the composition $L \xrightarrow{\alpha} A^0 \xrightarrow{F} \eta_p A ^0$.}
\item we use the same symbols $(d, F, F_L)$ to denote their restrictions to the corresponding (sub-)object. 
\ei
\ed
 
 Both $\eta_p (\alpha)$ and $\eta_p (\delta)$ are well-defined by the relations $d F =  p F d$ and $\delta F = p F \delta$.   By Remark 3.1.7 in \cite{BLM}, we know that $\eta_p (A^*)$ is a Dieudonn\'e algebra. It is then straightforward to verify that $\eta_p ({A}^*)$ is indeed a log Dieudonn\'e $\ul W$-algebra, with its $\ul W$-structure given by $W \ra A^0 \ra (\eta_p A)^0$ and $N \ra L$.
Now let $\textup{DA}^{\log}_{p\text{-tf}}$ be the full subcategory of $\AD$ spanned by $p$-torsion free complexes, then $\eta_p$ induces a functor $ \textup{DA}^{\log}_{p\text{-tf}} \longrightarrow  \textup{DA}^{\log}_{p\text{-tf}}$ which we call the log-d\'ecalage operator. For a morphism $f = (f, \psi): {A}^* \ra {B}^*$ between log Dieudonn\'e algebras, one obtains a morphism $\eta_p (f): \eta_p ({A}^*) \ra \eta_p({B}^*)$ in a functorial way. We only need to check that $f(\eta_p A)^i \subset (\eta_p B)^i$, which follows directly by unwinding definitions. The operator $\eta_p$ clearly preserves the $p$-compatible objects. 
Moreover, for a $p$-torsion free object ${A}^*$ in the category $\AD$ (resp. in $\ADp$), we have a natural map 
$${\phi}_F = (\phi_F, \textup{id}): {A}^* \longrightarrow \eta_p ({A}^*)$$ of log Dieudonn\'e algebras, where $\phi_F$ is the map $p^i F$ on $A^i$, which by  \ref{sss:phi_F_in_AD} is a morphism between Dieudonn\'e algebras. 

\bd 
A log Dieudonn\'e algebra ${A}^* \in  \AD$ is saturated if $A^*$ is $p$-torsion free and ${\phi}_F$ is an isomorphism, in other words, if the underlying Dieudonn\'e algebra is saturated. 
The full subcategory of $\AD$ spanned by saturated objects (resp. $p$-compatible saturated objects) is denoted by $\ADs$ (resp. $\ADps$). 
\ed

\br \label{remark:constructing_V} 
The definition of being saturated has nothing to do with the log structure $L$. Moreover, by the same discussion in Subsection \ref{sss:V}, a saturated log Dieudonn\'e algebra admits the Verschiebung operator $V: A^i \ra A^i$ satisfying the relations listed there. 
\er

Since the log structure is irrelevant in the condition of being saturated, the construction of the saturation of a log Dieudonn\'e algebra is essentially the same as in the non-logarithmic setup. First we replace $\ul A^* = (A^*, L)$ by $(A^*/A^*[p^\infty], L)$ to reduce to the $p$-torsion free case while keeping the monoid $L$.  For a $p$-torsion free log Dieudonn\'e algebra $ A^*$, we can now apply the log d\'ecalage operator $\eta_p$ to the morphism $\phi_F = (\phi_F, \textup{id}): {A}^* \ra \eta_p ( A^*)$  repeatedly and obtain 
$$ ( A^*)_{\mathfrak{sat}} := \textup{colim} \big( A^* \xrightarrow{\phi_F} \eta_p( A^*) \xrightarrow{\phi_F} \eta_p \eta_p( A^*) \ra \cdots  \: \big)$$ 
More explicitly, we may take $L_{\mathfrak{sat}} := \textup{colim} \Big( L \xrightarrow{F_L } F_L (L) \xrightarrow{F_L } F_L F_L (L) \ra  \cdots \Big) $
in the category of monoids ({colimit exists in the category of monoids, as both coproducts and coequalizers exist}). This is compatible with the colimit $\varinjlim_{n, \phi_F} \eta_p^{n} (A^*)$, hence we get natural monoid maps $  \alpha_{\mathfrak{sat}}: L_{\mathfrak{sat}} \ra (A_{\mathfrak{sat}})^0$ and  $\delta_{\mathfrak{sat}}: L_{\mathfrak{sat}} \ra (A_{\mathfrak{sat}})^1. $ 
The saturation of the log Dieudonn\'e algebra is then 
$$ \Big((A^*)_{\mathfrak{sat}}, L \ra L_{\mathfrak{sat}} \xrightarrow{\alpha_{\mathfrak{sat}}} (A_{\mathfrak{sat}})^0,  L \ra L_{\mathfrak{sat}} \xrightarrow{\delta_{\mathfrak{sat}}} (A_{\mathfrak{sat}})^1\Big).$$
The log structure and log derivation above can be equivalently described as the  compositions $ L \xrightarrow{\alpha} A^0 \ra  (A_{\mathfrak{sat}})^0$ and $L \xrightarrow{\delta} A^1 \ra  (A_{\mathfrak{sat}})^1, $ using the natural map $A^* \ra (A^*)_{\mathfrak{sat}}$ of Dieudonn\'e algebras. From this perspective it is clear that $(A^*)_{\mathfrak{sat}}$ forms a log Dieudonn\'e algebra over $\ul W$. 
The proof of the following lemma is omitted as it is completely analogous to the non-logarithmic situation. 
\bl \label{lemma:left_adjoint_1}
The saturation functor $\mathfrak{sat}: \AD \ra \ADs $ is the left adjoint of the inclusion $\ADs \subset \AD.$ The same is true for $\ADps$. 
\el

\subsection{Strict log Dieudonn\'e algebras and $V$-completion}   \label{ss:strict_completion}

Let $A^* \in \AD$ be a saturated object, recall that we can form the quotients $W_r (A^*) = A^*/\textup{Fil}^r$ where $\textup{Fil}^r$ is the $V$-filtration on $A^*$ given by the differential graded ideal  
$$ \textup{Fil}^r = V^r (A^*) + \textup{d} V^r (A^*) \subset A^*.$$
We then form the inverse limit $W(A^*) = \varprojlim W_r(A^*)$ along the restriction maps $R: W_r(A^*) \ra W_{r-1} (A^*).$   The Dieudonn\'e algebra $W(A^*)$ inherits a log structure $L \ra A^0 \ra W(A)^0$ and a log derivation $L \ra A^1 \ra W(A)^1$ from $A^*$, which makes it a log Dieudonn\'e algebra. The canonical map $\rho: A^* \ra W(A^*)$ of saturated Dieudonn\'e algebras upgrades to  $(\rho, \text{id}): A^* \ra W(A^*)$ as a morphism in $\AD$, which we still denote by $\rho$. By Remark \ref{remark:completion_of_A_is_saturated}, we know that, if $A^* \in \AD$ is saturated (resp. saturated and $p$-compatible), then $W(A^*)$ is again saturated (resp. saturated and $p$-compatible).  

\bd 
A saturated log Dieudonn\'e algebra $A^* \in \AD$ is strict if $\rho: A^* \ra W(A^*)$ is an isomorphism, in other words, if its underlying Dieudonn\'e algebra is strict. The full subcategory of $\ADs$ (resp. $\ADps$) spanned by strict objects is denoted by $\ADc$ (resp. $\ADpc$). 
\ed

The following lemma is immediate from Remark \ref{remark:W_r_induces_isom}. 
\bl \label{lemma:left_adjoint_2}
The completion $W(A^*)$ of a saturated log Dieudonn\'e algebra is strict.
Moreover the completion functor $W:  \ADs \ra \ADc$ provides a left adjoint of the inclusion $\ADc \subset \ADs$, similarly for the $p$-compatible subcategory. 
\el

The condition of being strict and $p$-compatible imposes some rather strict restrictions on the log structures. For the setup, let $A^* \in \ADpc$, and recall from Section \ref{sec:BLM} that we have a canonical isomorphism $\mu: A^0 \isom W(A^0/ V A^0)$ respecting Frobenii on both sides. We use $\mu$ to identify $A^0$ with the Witt vectors of $A^0/V A^0$ in the following proposition, which is important for our construction of log de Rham--Witt complexes in Section \ref{sec:log_dR_Witt}. 

\bp \label{prop:image_of_log}
Let $A^* \in \textup{DA}_{\mathfrak{str}}$  be a strict Dieudonn\'e algebra. Let $\alpha: L \ra A^0$ be a log structure on $A^0$, such that $F(\alpha(m)) = \alpha(m)^p$. Then the image of $\alpha$ consists of Teichmuller lifts of elements of $A^0/V A^0$. In other words, consider the log algebra $\cl \alpha: L \ra A^0 \ra A^0/V A^0$ with log structure inherited from $\alpha$, then $$\alpha (m) = [\cl \alpha (m)]$$
for any $m \in L$. This in particular applies to the case of a strict $p$-compatible log Dieudonn\'e algebra. 
\ep

\bproof 
It suffices to show that $\alpha(m) - [\cl \alpha(m)] \in \textup{Fil}^s (A^0) = V^s (A^0)$ for every $s \in \Z_{\ge 1}$ since $A^0$ is complete with respect to the $V$-filtration. We proceed by induction. By assumption $\alpha(m) - [\cl \alpha(m)] \in \textup{Fil}^1$. We claim that if $\alpha(m) - [\cl \alpha(m)] \in \textup{Fil}^{r}$ for $r \ge 1$, then $\alpha(m) - [\cl \alpha(m)] \in \textup{Fil}^{2r}.$ Suppose that $\alpha(m) - [\cl \alpha(m)] = V^rb$, then 
$$ F^r( \alpha(m)) - F^r([\cl \alpha(m)]) = \alpha(m)^{p^r} - [\cl \alpha (m)]^{p^r} = p^r \cdot b.$$ 
In the first equality we use that $A^*$ is $p$-compatible. On the other hand 
\begin{align*} \alpha(m)^{p^r} & = \Big( [\cl \alpha (m)] + V^rb \Big)^{p^r} \\ 
& =  [\cl \alpha (m)]^{p^r} + p^r  [\cl \alpha (m)]^{p^r -1} \cdot V^rb + \sum_{k = 2}^{p^r} {p^r \choose k}  [\cl \alpha (m)]^{p^r - k} \cdot \big( V^rb \big)^k 
\end{align*}
Note that since $x V^r y = V^r ((F^rx) y)$, we have
$$(V^rb)^2 = V^r((F^r V^r b)b) = V^r (p^r b^2) = p^r V^r(b^2).$$
Therefore
$$ p^r b =  p^r  [\cl \alpha (m)]^{p^r -1} \cdot V^rb + p^r V^r(b^2) \cdot \sum_{k = 2}^{p^r} {p^r \choose k}  [\cl \alpha (m)]^{p^r - k} \cdot \big( V^rb \big)^{k-2} $$
Since $A^0$ is $p$-torsion free, we conclude that $b \in V^r (A^0)$. 
\eproof 

\bc \label{cor:image_of_log}
Retain notations from above. Let $\ul{A^0}/V \ul A^0$ denote the log algebra $(A^0/V A^0, \cl \alpha: L \ra A^0/ V \ul A^0)$, and let $W(\ul A^0/V \ul A^0)$ be the log algebra defined as in Example \ref{sss:examples}. Then $\mu: A^0 \isom W(A^0/ V A^0)$ induces an isomorphism of log algebras 
$$\ul \mu= (\mu, \text{id}): \ul A^0 = (A^0, L) \ra W(\ul A^0/V \ul A^0).$$
\ec
\bproof 
This is immediate from Proposition \ref{prop:image_of_log}. 
\eproof

\newpage
 

\section{Saturated log de Rham--Witt complexes} \label{sec:log_dR_Witt} 

The goal of this section is to construct the saturated log de Rham--Witt complexes $\mW \omega^*_{\ul X/\ul k}$ for a log scheme $\ul X$ over $\ul k$, which is a sheaf on $X_{\ett}$ valued in log Dieudonn\'e algebras. We will construct this object locally on charts of the form $(R, M)$ in Subsections \ref{ss:construct_on_charts}, and compare it with de Rham complexes of $\ul R/\ul k$ (and its liftings to $\ul W$) in the subsequent subsections. In Subsection \ref{sec:etale_base_change} we glue the local construction on the $\etale$ site of $X$. 

For this section, we fix the base log algebra $\ul W = W(\ul k)$ where $\ul k = (k, N)$ as before\footnote{This is less general than the previous section, namely we require that the log structure on $W(k)$ comes from the Teichmuller lift from log structure on $k$. Examples include the trivial log structure on $W(k)$, or $\N \xrightarrow{1 \mapsto 0} W(k)$, but excludes $\N \xrightarrow{1 \mapsto p} W(k)$.}. In Subsection \ref{ss:construct_on_charts}, we do not impose any condition on the log structures, but later in this section we gradually impose restrictions on the log algebra to prove various forms of comparison theorems. 

\subsection{Log de Rham--Witt complexes for log algebras}  \label{ss:construct_on_charts}

Let $\textup{Alg}_{\ul{k}}^{\log}$ be the category of log algebras over $\ul k$. An element $\ul{R} = (R, M) \in \textup{Alg}_{\ul{k}}^{\log}$ can be equipped with a pair of Frobenius $(F, F_M)$ where $F(x) = x^p$ for all $x \in R$ and $F_M (m) = pm$ for all $m \in M$. From the previous section, we have a functor  
$$ \ADpc \longrightarrow \textup{Alg}_{\ul{k}}^{\log}  \qquad  \text{given by } \:\:   A^* \mapsto  \ul A^0/ V(A^0), $$ where $\ul A^0/ V(A^0) $ consists of the data 
$ (A^0/V(A^0),  \alpha: L \ra A^0 \twoheadrightarrow A^0/V(A^0))$. Now we are ready to give the definition of saturated log de Rham--Witt complexes for log algebras (Definition \ref{maindef:de_Rham_Witt} in the introduction).  
 
\bd 
Let $\ul{R} = (R, M) \in \textup{Alg}_{\ul{k}}^{\log}$ be a log algebra over $\ul{k}$. The saturated log de Rham--Witt complex $\mW \omega^*_{\ul{R}/\ul{k}}$ of $\ul{R}$ is a strict $p$-compatible log Dieudonn\'e $W(\ul k)$-algebra such that for any $A^* \in \textup{DA}^{\log}_{\mathfrak{str}}$, the natural map
$$\Hom_{\ADc} (\mW \omega^*_{\ul{R}/\ul{k}}, A^*) \longrightarrow  \Hom_{\textup{Alg}_{\ul{k}}^{\log}} (\ul{R}, \ul{A}^0/ V(A^0))$$ is an isomorphism. 
\ed 

In other words, the association $\ul{R} \mapsto \mW \omega^*_{\ul{R}/\ul{k}}$ (if exists) would provide a left adjoint of the functor 
$$ \textup{DA}_{\mathfrak{str}_{/\ul W}}^{\log, p} \longrightarrow \textup{Alg}_{\ul{k}}^{\log}, \qquad A^* \longmapsto (A^0/V\!(A^0), L).$$

\br
The restriction to the subcategory of $p$-compatible objects seems necessary to obtain the desired universal property. In particular, we do not get a left adjoint of $\ADc \ra \textup{Alg}_{\ul{k}}^{\log}$ in general. 
\er 

\subsubsection{Preliminaries \textup{I}}

\[
\begin{tikzcd}[column sep=1em,row sep=0.6em]
&L_A \arrow[rr, blue, "\psi"] \arrow[dd] & & L_B \arrow[dd] \\
L_A \arrow[ur, "F_{L_A}"] \arrow[rr, blue, crossing over, near end,  "\psi"] \arrow[dd] & & L_B  \arrow[ur, near start, swap, "F_{L_B}"] \\
& A \arrow[rr, blue, near start, "f"] & & B  \\
A \arrow[rr, blue, "f"]  \arrow[ur, "\varphi_A"] & & B  \arrow[ur, "\varphi_B"] \arrow[from=uu, crossing over]\\
\end{tikzcd} 
\] 
\bd \label{definition:Frob_compatible_homomorphism}  

Let $\ul{A} = (A, L_A)$ and $\ul{B} = (B, L_B)$ be two log algebras over $\ul W$ equipped with Frobenius maps $(\varphi_A, F_{L_A})$ and $(\varphi_B, F_{L_B})$ as morphisms of log $\ul W$-algebras, where $\varphi_A$ (resp. $\varphi_B$) lifts the $p$-power Frobenius on $A/p$ (resp. $B/p$). 
We define $\Hom_{F} (\ul A, \ul B)$ to be the set of Frobenius compatible $\ul W$-morphisms between $\ul A$ and $\ul B$, namely pairs $(f, \psi)$ of morphisms between log $\ul W$-algebras making every square in the cube above commute. 
\ed 

\br Consider the following special case of the definition above. Let $\ul{R}, \ul{R'} \in \textup{Alg}_{\ul{k}}^{\log}$ be two log algebras over $\ul{k}$. Equip $\ul{R}$ with the a log Froebnius map $(\varphi, F_M)$ where $F_M = p$. Now consider $W(\ul R) = (W(R), [\alpha])$, where 
$$[\alpha] = [\;] \circ \alpha: M \ra R \ra W(R)$$ is the Teichmuller lift of the log structure on $R$ as usual. $W(\ul R)$ is equipped with the Frobenius map $(F, p)$ from the Frobenius $F$ on the Witt vectors. Similarly define  $W(\ul R')$. Then $\Hom_{F} (W(\ul{R}), W(\ul{R'}))$ consists precisely of the homomorphisms between log algebras compatible with Frobenius $F$ on $W(R)$ and $W(R')$ respectively. More precisely, an element in $\Hom_{F} (W(\ul{R}), W(\ul{R'}))$ is a map $(f, \psi): W(\ul R) \ra W(\ul R')$ of log $\ul W$-algebras, such that $f \circ F = F \circ f$. 
\er

\bl  \label{lemma:equality_mod_V} 
Let $\ul{R}, \ul{R'} \in \textup{Alg}_{\ul{k}}^{\log}$ be as in the remark above. Suppose that both $R$ and $R'$ are reduced, then
$$\Hom_{F} (W(\ul{R}), W(\ul{R'})) \cong \Hom_{\textup{Alg}_{\ul{k}}^{\log}} (\ul{R}, \ul{R}')$$
\el

\bproof 
The Witt vector formalism provides a lift for morphisms from $R \ra R'$. Let $(\cl f, \psi): \ul R \ra \ul R'$ be a map over $\ul k$. We invoke the following property of Witt vectors: since $W(R)$ is $p$-torsion free by the reducedness assumption, and $F: W(R) \ra W(R)$ lifts the Frobenius on $R$, the map 
$$W(R) \ra R  \xrightarrow{\cl f} R'$$ lifts uniquely \footnote{This is sometimes called Cartier's Dieudonn\'e-Dwork lemma. Another way to produce this lifting is to observe that the Witt vectors provide a right adjoint to the forgetful functor from the category of $\delta$-rings to rings, in the sense of \cite{Joyal}.} to a map $f: W(R) \ra W(R')$ which is compatible with Frobenius on both $W(R)$ and $W(R')$. The same argument provides structure morphisms $W(k) \ra W(R)$ and $W(k) \ra W(R')$ such that $f: W(R) \ra W(R')$ is a map of $W(k)$-algebras (by the uniqueness of the lifting). 

It remains to show the compatibility between $f$ and $\psi$. This is implied by the compatibility of $f$ with the Teichmuller lifts: namely that $f([x]) = [\cl f (x)]$ for every $x \in R$. Since $W(R')$ is complete with respect to the $V$-filtration, it suffices to show that  $f([x]) - [\cl f (x)] \in \textup{Fil}^s = V^s W(R')$ for every $s \in \N$. We proceed by induction. Clearly $f([x]) - [\cl f (x)]  \in \textup{Fil}^1$, now suppose that $f([x]) - [\cl f (x)]  \in \textup{Fil}^r$, then $f([x]) = [\cl f(x)] + V^r(y)$ for some $y \in W(R')$, so we have 
$$F(f[x]) = f([x])^p = \Big( [\cl f(x)] + V^r(y) \Big)^p \equiv F[\cl f(x)] \mod p V^r W(R')$$
which implies that, there exists some $y' \in W(R)$ such that 
$$ F(f[x]) = F[\cl f(x)] + p V^r (y') = F \Big( [\cl f(x)] + V^{r+1}(y') \Big).$$
But $F$ is injective on $W(R')$ since $R'$ is reduced (so in particular it is $p$-torsion free), hence $f[x] = [\cl f(x)]$. 

Applying the same argument again, we see that the map $(f, \psi)$ is a morphism over $W(\ul k)$, which is equipped with the log structure $N \ra k \xrightarrow{[\:]} W(k)$. 
\eproof 

\subsubsection{A key example of log Dieudonn\'e algebra} 

The construction of log de Rham--Witt complexes for $\ul{R} \in \textup{Alg}_{\ul{k}}^{\log}$ relies on the following proposition. 

\bp \label{thm:omega_of_Witt_lifting_AD}
Let $(A, \alpha:L \ra A)$ be a log algebra over $\ul W$. Suppose that $A$ is $p$-torsion free and that the pair $(A, L)$ is equipped with Frobenius morphisms $\varphi: A \ra A$ and $F_L: L \ra L$ such that 
\bi
\item $\varphi (x) \equiv x^p \mod p$ is a lifting of the absolute Frobenius mod $p$; 
\item $F_L = p$ on $L$ and is compatible with $\varphi$, that is, $\varphi \circ \alpha = \alpha \circ F_L$. 
\ei
Then there exists a unique graded ring homomorphism 
$$ F: \omega^*_{\ul{A}/\ul{W}} \longrightarrow  \omega^*_{\ul{A}/\ul{W}}  $$ which extends $\varphi$ on $A$, such that $$ F (d x) = x^{p-1} dx + d (\frac{\varphi(x) - x^p}{p}) \quad \text{for all } x \in A.$$ 
Moreover, this makes $ \ul{\omega^*} := (\omega^*_{\ul{A}/W(\ul{k})}, L, d, \delta, F, F_L)$ a $p$-compatible log Dieudonn\'e $\ul W$-algebra. 
\ep

\bproof 

For the first claim we need to show the existence of $F: \omega^1_{\ul{A}/W(\ul k)} \ra \omega^1_{\ul{A}/W(\ul k)}$. We regard the second copy of $ \omega^1_{\ul{A}/W(\ul k)}$ as an $A$-module where $A$ acts by a twist of $\varphi$, namely $a \cdot \lambda = \varphi (a) \lambda$ for all $\lambda \in \omega^1_{\ul{A}/W(\ul k)}$. By the universal property of $ \omega^1_{\ul{A}/W(\ul k)}$ and Remark \ref{remark:pre_log_to_log_derivation}, we only need to supply the structure of a log derivation from $(A, L)$ to $\omega^1_{\ul{A}/W(\ul k)}$ with this twisted $A$-module structure. For this we define 
$$\sq \delta: L \ra   \omega^1_{\ul{A}/W(\ul k)}, \quad \text{and } \:\: \sq d: A \ra \omega^1_{\ul{A}/W(\ul k)}  $$
respectively as follows: 
\bi
\item we set $\sq \delta (l) :=  \delta (l)$, and 
\item $\sq d (x) :=   
x^{p-1}dx + d \theta (x)$ where $\theta(x) = \displaystyle \frac{\varphi(x) - x^p}{p} $. \footnote{Remark on notation: $\theta(x)$ here defines what is usually called a $\delta$-structure in literature, which features for example in the theory of prismatic cohomology. In this article we reserve the letter $\delta$ for log derivations.}
\ei 
  
To complete the proof of (1) we need to show that $\sq \delta$ and $\sq d$ form a log derivation of $\ul A$ over $W(\ul k)$, which requires
\be
\item $\sq d (\gamma x) =  \gamma \cdot \sq d x := F(\gamma) \sq d x$ for all $  \gamma \in W(k),  x \in A$;
\item $\sq d (x + y) = \sq d (x) + \sq d(y)$ for all $x, y\in A $;
\item $\sq d (x y) = x \cdot \sq d y + y \cdot \sq d x := \varphi(x) \sq d y + \varphi(y) \sq d x$ for all $x, y \in A$;
\item $\alpha (l) \cdot \sq \delta (l) : =  \varphi(\alpha (l)) \sq \delta (l) =\sq d (\alpha (l))$ for all $l \in L$; 
\item $\sq \delta (\iota (n)) = 0 $ for all $n \in N$, where $\iota: N \ra L$ is map of monoids in $W(\ul k) \ra \ul{A}.$  
\ee
(1) follows from (3) since $ p \sq d (\lambda) = d (F (\lambda)) = 0$, where $F: W(k) \ra W(k)$ is the Frobenius map on Witt vectors. (2) and (3) follow from the explicit description of $\theta(x)$, and the detail have been spelt out in \cite{BLM} Proposition 3.2.1, so we do not repeat the argument here. For (4),  since $\theta(\alpha(l)) = 0$, we have that 
$$ \varphi (\alpha(l)) \sq \delta (l) = \alpha(l)^p  \delta (l)  = \alpha(l)^{p-1} d (\alpha(l))  =   \sq d(\alpha (l)). $$ 
(5) is clear since $\delta (\iota (n)) = 0 $ for all $n \in N$. Therefore, we get the desired map 
$$F: \omega^1_{\ul{A}/W(\ul k)} \ra \omega^1_{\ul{A}/W(\ul k)} $$
which satisfies $F(\delta(l)) = \sq \delta(l) = \delta (l)$ and $F d(x) = \sq d (x) = x^{p-1} dx + d \theta (x)$.  \\

To finish the proof of the proposition, we need to check that (i). $ d F = p F d$ on all of $ \omega^*_{\ul{A}/W(\ul k)}$ and (ii). $\delta F_L = p F \delta$ on $L$. For (i), note that $\omega^*_{\ul{A}/W(\ul k)}$ is generated over $A$ by $\omega^1_{\ul{A}/W(\ul k)}$, and it is clear that we only need to check the relation on $x, dx $ for all $x \in A$, and on $\delta(l)$ for all $l \in L$. The fact that $d F (x) = p F d (x)$ follows directly from the construction of $F$ on $\omega^1_{\ul{A}/W(\ul k)}$, since $F d = \sq d$ by construction; on $d x$ and $\delta (l),$ both $d F$ and $p F d$ evaluate to $0$ (since $d \delta = 0$ on $L$). Part (ii) is automatic by construction of $F$ again, since $F \delta =\delta$.  
\eproof 

We record a variant of the proposition above, which allows slightly more flexible Frobenius $F_L$ on the monoid $L$. 

\bp \label{thm:omega_of_Witt_lifting_AD_variant}
Let $(A, \alpha:L \ra A)$ be a log algebra over $\ul W$.  
Suppose that $A$ is $p$-torsion free and $p$-separated, and that $\omega^1_{\ul{A}/\ul{W}}$ is $p$-torsion free. Also suppose that the pair $(A, L)$ is equipped with Frobenius morphisms $\varphi: A \ra A$ and $F_L: L \ra L$ such that 
\bi
\item $\varphi (x) \equiv x^p \mod p$ is a lifting of the absolute Frobenius mod $p$; 
\item $F_L$ is compatible with $\varphi$.  
\ei
Then the same conclusions of Proposition \ref{thm:omega_of_Witt_lifting_AD} hold except for $p$-compatibility. 
\ep

\bproof  First we need the following sublemma: 
\begin{sublemma} Retain the notation and assumption of the proposition. For any $k \in \N$, if $p^k | \varphi(x)$, then  $p^k | \big(x^{d-1} dx + d \theta(x) \big).$ \end{sublemma}
\bproof[Proof of the sublemma] The case $k = 0$ is tautological. We proceed by induction. Suppose that $k \ge 1$ and that the sublemma has been verified up to $k - 1$. Now let $\varphi (x) = p^{k} y = x^p + p \theta (x)$, since $p | \varphi (x)$, we know that $p | x$, so $x = p z$ for some $z \in A$. Since $A$ is $p$-torsion free, $\varphi (x) = p \varphi(z) = p^{k} y$ implies that $\varphi (z) = p^{k - 1} y$, so $p^{k - 1} |  \big({z}^{p-1} dz + d \theta(z) \big) $ by induction hypothesis. Note that ${x}^{p-1} dx + d \theta(x) = (pz)^{p-1} d(pz) + d \theta (pz) = p \big( {z}^{p-1} dz + d \theta(z) \big)$, so the lemma follows.  \eproof
\noindent Now we proceed as in the proof of Proposition \ref{thm:omega_of_Witt_lifting_AD}, but replace $\sq \delta$ by $$\sq \delta (l) := \displaystyle {\delta (F_L(l))}/{p}.$$ 
We need to show that  $p | \delta (F_L (l))$ for any $l \in L$ for $\sq \delta$ to be well-defined. For this we need the assumption that $A$ is $p$-separated: write $x = \alpha(l)$, then there exists some $k \in \N$ such that $p^k | \varphi (x)$ but $p^{k + 1} \nmid  \varphi (x)$. By the sublemma above we know that $p^{k+1} | d (\varphi(x))$. Now since  $$ \varphi(x)\; \delta(F_L (l)) =  \alpha (F_L (l))\; \delta(F_L(l)) = d \alpha (F_L(l)) = d \varphi(x) $$ and that $p^{k+1} \nmid \varphi(x)$, $p^{k+1} | d \varphi(x)$, we know that $p | \delta (F_L (l))$ (here we use the assumption that $ \omega^1_{\ul{A}/W(\ul k)} $ is $p$-torsion free). The rest of the proof of part (1) is the same, except that to check $\alpha (l) \cdot \sq \delta (l) = \sq d(\alpha (l))$, we need 
\begin{align*} p \sq d (\alpha (l)) & = d \varphi (\alpha (l )) \\ & = d \alpha (F_L (l)) \\ & =  \alpha (F_L (l)) \delta(F_L (l)) \\ & = \varphi (\alpha(l)) \delta (F_L (l)) \\ & =  p  \varphi (\alpha(l)) \sq \delta (l) = p \alpha(l) \cdot \sq \delta(l).
\end{align*} 
The rest of the proof is identical to that of Proposition \ref{thm:omega_of_Witt_lifting_AD}. 
\eproof 

\subsubsection{Frobenius on finite level}

\br \label{remark:Frob_finite_level} 

Proposition \ref{thm:omega_of_Witt_lifting_AD} will be applied to $A = W(R)$ where $R$ is a reduced $\F_p$-algebra. For applications later, we remark that the graded ring homomorphism $F$ can be constructed on the finite length Witt vectors, namely there is a natural graded ring homomorphism
$$F: \omega^*_{W_n(\ul R)/W_n(\ul k)} \ra \omega^*_{W_{n-1}(\ul R)/W_{n-1}(\ul k)} $$ extending the Frobenius $F: W_{n} (R) \ra W_{n-1} (R)$. We prove the claim using the same strategy as Proposition \ref{thm:omega_of_Witt_lifting_AD}, namely we need to construct a log derivation $(\sq d, \sq \delta)$ of $W_n(\ul R)/W_n(\ul k)$ into $F_*\omega^1_{W_{n-1}(\ul R)}$, where $\omega^1_{W_{n-1}(\ul R)}$ is regarded as a $W_n(R)$ structure via $F: W_n(R) \ra W_{n-1} (R)$. We define $\sq d$ as follows: write $x \in W_n(R)$ uniquely as $x = [x_0] + V x'$ with $x' \in W_{n-1}(R)$ and $x_0 \in R$, then define 
$$\sq d x := [x_0]^{p-1} d[x_0] + d x'$$  
where $[x_0]$ now lives in $W_{n-1}(R)$. Then define $\sq \delta (m) := \delta (m)$.  

$\bullet$ Let us first check that $\sq d (x + y) = \sq d(x) + \sq d(y)$, for which we only need to show that if $[x_0] + [y_0] = [x_0+y_0] + V z'$, then 
$$[x_0]^{p-1} d[x_0] + [y_0]^{p-1} d[y_0] = [x_0+y_0]^{p-1} d[x_0+y_0] + dz'.$$
It suffices to check the equality above in $\omega^1_{W(\ul R)/W(\ul k)}$. Now for any element $z = [x] + Vy \in W(R)$, we claim that
$$ [x]^{p-1}d[x] + d y = z^{p-1} dz + d \Big(\frac{F(z) - z^p}{p} \Big).$$ 
This is straightforward by expanding the terms on the right hand side: 
$$\big([x] + Vy \big)^{p-1} d\big([x] + Vy \big) + d \Big(\frac{ [x]^p + py - \big([x] + Vy \big)^p}{p}\Big).$$ 
Now the desired additivity follows from the proof of Proposition \ref{thm:omega_of_Witt_lifting_AD}. 

$\bullet$ Next we have 
\begin{align*} 
\sq d (x y) & = [x_0 y_0]^{p-1} d \Big([x_0 y_0]\Big) + d \Big([x_0]^p y' + [y_0]^p x' + p x' y' \Big) \\
& =  \Big([x_0]^{p} + px' \Big)\Big([y_0]^{p-1} d[y_0] + d y' \Big) +  \Big([y_0]^{p} + py' \Big)\Big([x_0]^{p-1} d[x_0] + d x' \Big)  \\
& = F(x) \sq d y + F(y) \sq d x = x \cdot \sq d y + y \cdot \sq d x
\end{align*}

$\bullet$ Finally, $(\sq d, \sq \delta)$ form a log derivation of $W_n(\ul R)/W_n(\ul k)$ into $\omega^1_{W_{n-1}(\ul R)/W_{n-1} (\ul k)}$, inducing the desired map $F: \omega^1_{W_{n}(\ul R)/W_{n}(\ul k)} \ra \omega^1_{W_{n-1}(\ul R)/W_{n-1}(\ul k)}$, which extends to a graded ring homomorphism on $\omega^*_{W_{n}(\ul R)/W_n(\ul k)}$. 
\er

\subsubsection{Preliminaries \textup{II}} 

\bl \label{lemma:homomorphism_F_compatible_algebra}
Let $\ul{A}$ be a log algebra over $\ul W = W(\ul k)$ satisfying the conditions in either Proposition \ref{thm:omega_of_Witt_lifting_AD} or \ref{thm:omega_of_Witt_lifting_AD_variant}, and let $\ul{\omega^*} = \ul{\omega}^*_{\ul A/ \ul W}$ be the corresponding log Dieudonn\'e $\ul W$-algebra given by Proposition \ref{thm:omega_of_Witt_lifting_AD} or \ref{thm:omega_of_Witt_lifting_AD_variant}.  Let $B^*$ be a $p$-torsion free log Dieudonn\'e $\ul W$-algebra. Then we have a canonical isomorphism between 
$$ \Hom_{\textup{DA}^{\log}} (\ul{\omega^*}, B^*) \cong \Hom_{F} (\ul A, \ul B^0)$$
where $\Hom_F$ is as defined in Definition \ref{definition:Frob_compatible_homomorphism}. 
\el

\bproof 
Given a morphism $(f, \psi) \in \Hom_{F} (\ul A, \ul B^0)$, we need to extend it to a morphism $f: \ul{\omega^*} \ra \ul B^*$ of log Dieudonn\'e $\ul W$-algebras. As in the proof of Proposition \ref{thm:omega_of_Witt_lifting_AD}, we construct a map  $f: \omega^1_{\ul{A}/\ul{W}} \ra B^1$ and extend to a morphism of cdga's on $\omega^*_{\ul{A}/\ul{W}}$, and check that it is compatible with log differentials $\delta$ and Frobenius $F$. In particular, we need the following diagram to commute:  
\[ 
\begin{tikzcd} [column sep=1.35em,row sep=1.35em]
A \arrow{r}{d} \arrow{d}{f} & \omega^1_{\ul{A}/W(\ul{k})} \arrow[dashed]{d}{f} \\
B^0 \arrow{r}{d_B} & B^1   
\end{tikzcd}
\]
Define $d_f: A \longrightarrow B^1$ by $d_f = d_B \circ f$ and $\delta_f: L \longrightarrow B^1$ by $\delta_f  = \delta_B \circ \psi$.   One easily verifies that $d_f (x + y) = d_f (x) + d_f (y)$, $d_f (x y) = f(x) d_f (y) + f(y) d_f (x)$, and that $d_f (\lambda x) = \lambda d_f (x)$. So $d_f$ is a derivation of $A/W(k)$ into $B^1$. To verify that $(d_f, \delta_f) $ is a log derivation of $\ul A/\ul W$, we check that $\delta_f (\iota(n)) = 0$ and that 
\begin{align*}  
\alpha (l) \cdot \delta_f (l) & = f(\alpha(l)) \: \delta_B (\psi (l)) \\ 
& = \alpha_B (\psi (l)) \: \delta_B (\psi (l)) \\
& = d_B \: \alpha_B (\psi (l)) = d_B \: f (\alpha (l)) = \delta_f (\alpha (l)). 
\end{align*}
Therefore, we have the dotted map $f: \omega^1_{\ul{A}/{W(\ul k)}} \ra B^1$ in the diagram above. 

$\omega^*_{\ul{A}/\ul{W}} $ is generated over $A$ by $\omega^1_{\ul{A}/\ul{W}}$, so we extend $f$ to $\omega^*_{\ul{A}/\ul{W}} $ as a differential graded morphism. Note that $f$ is compatible with $\delta$ by construction, so it remains  to show that $f$ is compatible with $F$. For this we proceed as the proof of Proposition \ref{thm:omega_of_Witt_lifting_AD}, namely we check it on $x, dx $ for $x \in A$  and $\delta(l)$ for $l \in L$   since both $F $ and $f$ are algebra morphisms. We know that $F f (x) = f F (x)$ by assumption. For $d x$ (where $x \in A^0$), we have
\begin{align*}
p f F (dx) & = f d F(x) = d( f F(x)) \\ & = d F (f (x)) = p F (d f(x)) = p F f (d x)
\end{align*} 
and since $B^*$ is $p$-torsion free, we get $f F (dx) = F f (dx)$. Finally for $\delta (l)$, we have (in both setups) $\delta F_{L_A} = p F \delta$, therefore
\begin{align*} 
p f F (\delta (l)) & = f \delta (F_{L_A} (l)) = \delta_B \psi ( F_{L_A} (l)) \\ &  = \delta_B (F_{L_B} \psi (l)) = p F \delta_B (\psi (l)) = p F f (\delta (l))
\end{align*}
again this implies that $ f F (\delta (l)) = F f (\delta (l))$. 
\eproof 

\subsubsection{Existence of saturated log de Rham--Witt complexes}
 
Now we prove our first main theorem (Theorem \ref{mthm:log_dR_is_DA} in the introduction). We continue to assume that the base log point $\ul k$ has the form $(k, N)$ where $N \minus \{0\} \mapsto 0 \in k$. 

\bt \label{thm:log_dRW_complex_exist}
Let $\ul{R} = (R, M) \in \textup{Alg}_{\ul{k}}^{\log}$ be a log algebra over $\ul k = (k, N)$. The saturated log de Rham--Witt complex $\mW \omega^*_{\ul R/\ul k}$ of $\ul R/\ul k$ exists. Moreover, the association of $\ul{R} \mapsto  \mW \omega^*_{\ul R/\ul k} $ is functorial.
\et 

\bproof 

Without loss of generality assume that $R$ is reduced by replacing it with its reduced quotient if necessary (for an object $B^* \in \ADpc$, $B^0/V(B^0)$ is reduced by Remark \ref{remark:completion_of_A_is_saturated}). Let  $W(\ul R) = (W(R), L =M \xrightarrow{[\alpha]} W(R))$ be the Witt vector of $\ul R$. Since $R$ is reduced, $W(R)$ is $p$-torsion free. \footnote{and in fact $p$-separated (since $W(R) \hookrightarrow W(K) \hookrightarrow W(K^{\text{perf}})$ where $K = \text{Fr}(R)$), however $\omega^1_{W(R)/W(k)}$ is not p-torsion free.} By Proposition \ref{thm:omega_of_Witt_lifting_AD}, there exits a unique log Dieudonn\'e $\ul W$-algebra on the relative log de Rham complexes:
$$ \ul{\omega^*} = (\omega^*_{W(\ul{R})/W(\ul{k})}, L, d, \delta, F, F_L) $$ where $L = M, F_L = p$ and $F (dx) = x^{p-1} dx + d(\frac{\varphi(x) - x^p}{p})$. Define 
$$ \mW \omega^*_{\ul R/\ul k}  = W (\ul{\omega^*})_{\mathfrak{sat}}.$$
The construction is clearly functorial in $\ul{R}$. We claim that $\mW \omega^*_{\ul R/\ul k}$ is the log de Rham--Witt complex of $\ul R/\ul k$. To prove the claim, consider any strict $p$-compatible log Dieudonn\'e $W(\ul k)$-algebra $B^* \in \textup{DA}^{\log, p}_{\mathfrak{str}}$, then we have
\begin{align*}
\Hom_{\ADpc} (\mW \omega^*_{\ul R/\ul k}, B^*) 
         & =  \Hom_{\ADp}  (\ul{\omega^*}, B^*) &\text{by Lemma } \ref{lemma:left_adjoint_2}  \\
           & =   \Hom_F(W(\ul{R}), \ul{B}^0)   & \text{by Lemma } \ref{lemma:homomorphism_F_compatible_algebra}  \\
          & =   \Hom_F(W(\ul{R}), W(\ul{B}^0/V \ul{B}^0))   & \text{by Cor } \ref{cor:image_of_log}  \\
          & =  \Hom_{\textup{Alg}^{\log}_{\ul{k}}} (\ul R, \ul{B}^0/V \ul B^0)  & \text{by Lemma } \ref{lemma:equality_mod_V}       
\end{align*}
\eproof 

\begin{notation} 
We denote $ \mW_n \omega^*_{\ul R/\ul k} := W_n \big( (\ul{\omega^*})_{\mathfrak{sat}} \big)$. 
\end{notation}

\br \label{remark:replace_with_finite_level}
For each $n \ge 1$, write $\omega^*_n :=  \omega^*_{W_n(\ul R)/W_n (\ul k)}$, then in the construction above we may use $\varprojlim \omega^*_n$ instead of $\omega^*_{W(\ul R)/ W(\ul k)}$. In other words, consider $\varprojlim \omega^*_n$ as a $p$-compatible log Dieudonn\'e algebra (see for example Remark \ref{remark:Frob_finite_level}), then we may define $\mW \omega^*_{\ul R /\ul k}$ as $W (\varprojlim \omega^*_n)_{\sat}$. From the proof of Theorem \ref{thm:log_dRW_complex_exist}, we need to show that the canonical map $\omega^*_{W(\ul R)/W(\ul k)} \ra \varprojlim \omega^*n$ induces a bijection \footnote{This is slightly different from the setup in Lemma \ref{lemma:homomorphism_F_compatible_algebra_2}, where $\widehat \omega^*$ there denotes the $p$-adic completion, while our limit here is the ``V-completion."}
$$\Hom_{\ADp}  (\ul{\omega^*}, B^*)  \isom \Hom_{\ADp} (\varprojlim \omega^*_n, B^*).$$
For this, it suffices to show that any map $\omega^* \ra B^*$ of log Dieudonn\'e algebras induces a map $\omega^*_n \ra W_n(B^*)$, which is compatible with transition maps as well as Frobenius maps on both sides, therefore it gives rise to $\varprojlim \omega^*_n \ra B^*$ on the inverse limit where the original map factors through. Clearly, we have $\omega^0_n = W_n(R) \ra W_n B^0 = B^0 /V^n (B^0)$. Similar to the proof of Lemma \ref{lemma:homomorphism_F_compatible_algebra}, we want to show that the map $\omega^1_{W(\ul R)/W(\ul k)} \ra B^1$ factors through $\omega_n^1 \ra W_n B^1$. This follows from the description of $\omega^1_{W(\ul R)/W(\ul k)}$, which is given by the quotient $$\Omega^1_{W(R)/W(k)} \oplus (W(R) \otimes M^{\textup{gp}})/\sim$$ where relations are specified in Subsection \ref{sss:log_differentials}. This tells us that the kernel of $\omega^1_{W(\ul R)/W(\ul k)} \ra \omega^1_n$ is generated by $\im (V^n) + \im (d V^n)$ as desired. 
\er

\subsection{Choice of charts}

In forming the saturated log de Rham--Witt complexes for a saturated log algebra $(R, M)$, we might have different log structures on $R$ which give rise to the same affine log scheme $(\spec R, \mM^a)$, which amounts to different choices of charts for the log structure $\mM^a$. This choice is almost irrelevant in our discussion, since the underlying Dieudonn\'e algebras will be the same. Now we make this more precise. Let $(R, \alpha: M \ra R)$ be a log algebra over $\ul k$. Following our convention in Introduction \ref{ss:conventions}, let $\mM^a$ be the log structure associated to the constant pre-log structure $M$ on $\spec R$. Let $M^{\textup{sh}} = \Gamma (\spec R, \mM^a)$ and denote the log algebra $(R, M^{\textup{sh}})$ by $\ul{R}^{\textup{sh}}$. 

\bp \label{prop:choice_of_pre_log_structures}

The induced map $ \mW \omega^*_{\ul R/ \ul k} \ra  \mW \omega^*_{\ul{R}^{\textup{sh}}/ \ul k}$  of the saturated log de Rham--Witt complexes from $\ul R \ra \ul R^{\textup{sh}}$ is an isomorphism on the underlying log Dieudonn\'e algebras. In particular, if $(R, M)$ and $(R, M')$ are two charts for the affine log scheme $(\spec R, \mM^a)$, then the underlying Dieudonn\'e algebra of their saturated log de Rham--Witt complexes are isomorphic.  
\ep

\bproof 
It suffices to show that for each $n \ge 1$, the log algebras $\beta: M \ra W_n(R)$ and $\beta^{\textup{sh}}: M^{\textup{sh}} \ra W_n(R)$ give rise to the same log structures over $\spec W_n(R)$, since it implies that there is an isomorphism
$$\omega^*_{\big(W_n(R), M \big)/W_n(\ul k)} \isom  \omega^*_{\big(W_n(R), M^{\textup{sh}} \big)/W_n(\ul k)}.$$ The proposition then follows from Remark \ref{remark:replace_with_finite_level}. 

Now write $X = \spec R$. Let $\mO_X$ be the structure sheaf of the $\etale$ site of $X$, and $W_n (\mO_X)$ be the sheaf sending $U \mapsto W_n (\Gamma(U, \mO_X))$.
\footnote{It is straightforward to directly check that the pre-sheaf $W_n (\mO_X)$ is a Zariski sheaf, and satisfies the sheaf condition for an $\etale$ cover $R \ra R'$. It corresponds to the structure sheaf on $(\spec W_n(R))_{\ett}$ under the identification of $\etale$ topos $X_{\ett}^{\sim} \cong (\spec W_n (R))_{\ett}^{\sim}$ by topological invariance of the $\etale$ site and the fact that, for any $\etale$ morphism $R \ra R'$, the induced map $W_n(R) \ra W_n (R')$ is $\etale$.} We need to show that the following two pushouts 
\[ 
\begin{tikzcd}[column sep=1.35em]
 \beta^{-1} (\mO_{\spec W_n(R)}^\times ) \arrow{r}{} \arrow{d}{} & \mO_{\spec W_n(R)}^\times \\ M 
\end{tikzcd} \qquad
\begin{tikzcd}[column sep=1.35em]
 (\beta^{\textup{sh}})^{-1} (\mO_{\spec W_n(R)}^\times ) \arrow{r}{} \arrow{d}{} & \mO_{\spec W_n(R)}^\times \\ M^{\textup{sh}}  
\end{tikzcd}  
\] of $\etale$ sheaves of monoids over $(\spec W_n(R))_{\ett}$ 
are isomorphic. By topological invariance of $\etale$ sites, we may compute the pushout of sheaves over $X_{\ett}$. Over $X$, the pre-log structure $\beta$ corresponds to morphism of monoids
$$\beta: M \xrightarrow{\alpha} \mO_X \xrightarrow{\gamma}  W_n(\mO_X)$$
where $\gamma$ is given by Teichmuller liftings. Note that $\gamma^{-1} (W_n (\mO_X)^\times) = \mO_X^\times$ since for any algebra $R$ the intersection $W_n(R)^\times \cap [R]$ is precisely $[R^\times]$. We have $\beta^{-1} (W_n (\mO_X)^\times) = \alpha^{-1} (\mO_X^\times)$, so the first pushout can be computed by two steps, as in the diagram 
\[
\begin{tikzcd}[column sep=1.35em]
 \alpha^{-1} (\mO_X^\times) \arrow{r}{\alpha} \arrow{d}{} & \mO_X^\times \arrow[r, "\gamma"]  & W_n (\mO_X)^\times \\ M 
\end{tikzcd}. 
\]
The first pushout precisely computes $\mM^a$, which is isomorphic to $(\mM^{\textup{sh}})^a$ by Lemma \ref{lemma:pushout_agrees}, therefore $M$ and $M^{\textup{sh}}$ indeed induces the same log structures on $\spec W_n (R)$, hence the Proposition follows. 
 \eproof 

\br 
The proof we give of Proposition \ref{prop:choice_of_pre_log_structures} is somewhat indirect (for example it involves Remark $\ref{remark:replace_with_finite_level}$). One issue is that we do not know whether $M$ and $M^{\textup{sh}}$ give rise to the same log structures over $W(R)$. However, if we assume that $\ul R/\ul k$ is log smooth of log Cartier type, then the proposition immediately follows from Corollary \ref{cor:Cartier_implies_qi} and Proposition \ref{prop:comparison_with_mod_p_dR}, since we have the following commutative diagram of cdga's
\[
\begin{tikzcd} 
\omega^*_{\ul R/\ul k} \arrow[r, "\sim"] \arrow[d] & \mW_1 \omega^*_{\ul R/\ul k} \arrow[d] \\ 
\omega^*_{\ul R^{\textup{sh}}/\ul k} \arrow[r, "\sim"] & \mW_1 \omega^*_{\ul R^{\textup{sh}}/\ul k} 
\end{tikzcd}
\] where the left vertical arrow is an isomorphism by Lemma \ref{lemma:pushout_agrees}. 
\er

\subsection{Construction via log Frobenius liftings} \label{ss:construction_via_Frob}
In this subsection we discuss another construction of saturated log de Rham Witt complexes for log algebras which admit liftings to $W(k)$ together with lifts of Frobenius.  

\subsubsection{Log Frobenius liftings} We record a (variant of a) definition from \cite{LZ}.  

\bd   Let $(f, \psi):  \spec \underline{R} = (\spec R, M) \ra \spec \underline{k}$ be a morphism of affine pre-log schemes with integral pre-log structures.
\be
\item A Witt-lifting of $f$ consists of a system 
$$\big(\underline{A_n} = (A_n, M_n), \:\: \delta_n: W_n (\spec \underline{R}) \ra  (\spec A_n, M_n) \big)_{n \ge 1},$$ where $(A_1, M_1) = (R, M)$, and satisfies the following list of conditions: 
\bi 
\item  $\spec \underline{A_n} = (\spec A_n, M_n)$ is log-smooth over $W_n(\spec \underline{k})$;
\item For each $n$, the following diagram is cartesian
\[
\begin{tikzcd} [column sep=1.2em,row sep=1.5em]
\spec \ul{A}_n \arrow{d}{} \arrow[r, hook, "R"] \arrow[rd, phantom, near start, "\square"] & \spec \ul{A}_{{n+1}} \arrow{d}{}  \\
W_n(\spec \ul{k} ) \arrow[r, hook, "R"] & W_{n+1}(\spec \ul{k} )
\end{tikzcd}
\] where $R: A_{n+1} \ra A_{n}$ is the transition map of the inverse system. 
\item The maps $\delta_n$ are compatible with the system $(\spec \ul{A}_n)$ under the transition maps: namely, for each $n$, the following diagram commutes 
\[
\begin{tikzcd} [column sep=1.2em,row sep=1.5em]
W_n(\spec \ul{R} ) \arrow[r, hook, "R"] \arrow{d}{\delta_n} & W_{n+1}(\spec \ul{R} ) \arrow{d}{\delta_{n+1}} \\
\spec \ul{A_n}  \arrow[r, hook, "R"] & \spec \ul{A}_{{n+1}}  
\end{tikzcd}
\]
\ei
\item A log Frobenius lifting of $f$ consists of a projective system 
$$( \ul{A_n}, \delta_n, \varphi_n)_{n \ge 1}$$ 
where $(\ul{A_n}, \delta_n)$ is a Witt-lifting of $f$, and $\varphi_n: \spec \ul{A}_{n} \ra \spec \ul{A}_{n+1}$ is a collection of morphisms of pre-log schemes satisfying:
\bi
\item It is compatible with the Frobenius on $\spec \ul{R}$ given by the $p^{th}$ power Frobenius on $R$ and multiplication by $p$ on $M$.
\item It is compatible with the Frobenius of Witt vectors on both $\spec \ul{k}$ and $\spec \ul{R}$. In other words, the following diagrams commute: 
\[
\begin{tikzcd} [column sep=0.2em,row sep=1.2em]
& \spec \ul{A}_{n+1}  \arrow[dd] \arrow[rr, "R"] &&  W_{n+1} (\spec \ul{R}) \\
\spec \ul{A}_{n} \arrow[rr, crossing over, near end, "R"] \arrow[dd] \arrow[ru, "\varphi_n"]  &&W_{n}(\spec \ul{R} ) \arrow[ru, "F_n"]  \\ 
& W_{n+1} (\spec \ul{k})\\
W_{n} (\spec \ul{k} ) \arrow[ru, "F_{n}"]
\end{tikzcd}
\]
\ei
\item A log Frobenius lifting $(\ul{A}_n = (A_n, M_n), \delta_n, \varphi_n)$ of $f$ is called $p$-compatible, if the log-restriction map $R: M_{n+1} \ra M_{n}$ is identity for all $n \ge 1$ (so in particular $M_n = M$ for all $n$) and the log Frobenius map $\varphi_n: M \ra M$ is multiplication by $p$. 
\ee
Though we stated the definitions on log algebras following \cite{LZ}, it is clear that one can globalize and extend the definition to pre-log schemes in general.
\ed

\bl \label{lemma:log_Frob_lifting_exist}
Suppose that $\ul R$ is log-smooth over $\ul k$, then there exists a $p$-compatible log Frobenius lifting $(\ul{A}_n, \delta_n, \varphi_n)$ of $\ul R$.  
  More generally, let $f: \ul{X} \ra \spec\ul{k}$ be a log-smooth morphism of fine (pre-)log schemes, then a $p$-compatible log Frobenius liftings exist $\etale$ locally.  
\el 
 
\bproof 
This follows from Proposition 3.2 in \cite{LZ} and Lemma 5.5 of \cite{Ma}.  The morphism $(k, N) \ra (R, M)$ factors through the morphism $(k, N) \ra (k \otimes_{\Z[N]} \Z [M], M) \ra (R, M)$, where the second map is $\etale$ on the underlying rings. The first arrow admits a $p$-compatible log Frobenius lifting given by $(T_n, M)$ where $T_n := W_n (k)\otimes_{\Z[N]} \Z[M]$ and $\varphi_n$ is given by $a \otimes b \mapsto F(a) \otimes b^p \in T_{n-1}$.
Then we need to lift the $\etale$ morphism $T_1 \ra R$ along $\cdots T_n \ra T_{n-1} \ra  \cdots \ra T_1$, which exists by \cite{LZ}  Proposition 3.2. Note that the log structure on each $A_n$ is given by $M \ra T_n \ra A_n$ (while $\varphi_n |M$ is still multiplication by $p$), so in particular we have constructed a $p$-compatible log Frobenius lifting. 
\eproof

\subsubsection{Log de Rham--Witt complex via log Frobenius liftings} \label{sss:via_Frob_liftings}   We make the following assumptions for the rest of this Subsection. 

\begin{framed}  \noindent Suppose that $\ul R/\ul k$ is integral and admits a $p$-compatible log Frobenius lifting $\ul A = (\ul A_n, \delta_n, \varphi_n)$.
\end{framed}
\noindent Typical examples include log algebras $\ul R$ that are log-smooth and integral over $\ul k$ by Lemma \ref{lemma:log_Frob_lifting_exist}. Note that our definition of log-smooth requires the monoid $M$ to be integral and coherent, but this is different from requiring the morphism $\spec \ul R \ra \spec \ul k$ to be integral, see Appendix \ref{sss:integral_morphism}.

Write $\ul{\widehat A} = (\widehat A, M)$ where $\widehat A = \varprojlim A_n$ and it is equipped with the log structure from the filtered inverse limit. $\widehat A$ defines an affine formal scheme with the $p$-adic topology.  Then we may form the completed log de Rham complex 
$$ \widehat{\omega}^*_{\ul{\widehat A} /W(\ul k)} (= {\widehat \omega}^*) := \varprojlim_n \omega^*_{\ul A_n/W_n (\ul k)} \cong \varprojlim_n \big( {\omega}^*_{\ul{\widehat A}/W(\ul k)}/p^n \big). $$ 
From Subsection \ref{sss:integral_morphism}, $f$ being integral means that the map $\Z[N] \ra \Z[M]$ induced from the map of monoids is flat, therefore for each $n$, $\spec \ul A_n \ra \spec W(\ul k)$ is log-smooth and integral, hence flat. Now we observe that the restriction map $R: A_{n+1} \ra A_n$ is surjective, so the inverse limit $\widehat A$ is also flat over $W(k)$, hence $p$-torsion free. This allows us to apply Proposition \ref{thm:omega_of_Witt_lifting_AD} and get a log Dieudonn\'e algebra $\omega^*_{\ul{\widehat A}/W(\ul k)}$. The completed log de Rham complex 
$\widehat \omega^*_{\ul{\widehat A}/W(\ul k)}$ has a unique cdga structure such that the canonical map  $\omega^*_{\ul{\widehat A}/W(\ul k)} \ra \widehat \omega^*_{\ul{\widehat A}/W(\ul k)}$ is a map of cdga's. Moreover, $\widehat \omega^*_{\ul{\widehat A}/W(\ul k)}$ inherits the Frobenius structure from $\omega^*_{\ul{\widehat A}/W(\ul k)}$ which makes it a $p$-compatible log Dieudonn\'e algebra. 

\subsubsection{Preliminaries \textup{III}} The following lemma is a variant of Lemma \ref{lemma:homomorphism_F_compatible_algebra}. In this lemma we do not need to assume that $\ul A$ comes from log Frobenius liftings.  

\bl \label{lemma:homomorphism_F_compatible_algebra_2} 
Let $\ul A$ be a log algebra over $W(\ul k)$ satisfying the conditions in either Proposition \ref{thm:omega_of_Witt_lifting_AD} or Proposition \ref{thm:omega_of_Witt_lifting_AD_variant}, and let ${\omega^*} = \ul{\omega}^*_{\ul A/W(\ul k)}$ be the log Deiudonne algebra constructed there. Let ${\widehat \omega}^*$ be the $p$-completion of $\omega^*$, which is a $p$-complete log Dieudonn\'e $W(\ul k)$-algebra. Then for any $p$-torsion free and $p$-complete log Dieudonn\'e $W(\ul k)$-algebra $B^*$, the canonical map
$$ \Hom_{\textup{DA}^{\log}} (\widehat \omega^*, B^*) \cong \Hom_{F} ({\ul A}, \ul B^0)$$
is a bijection. Suppose that $B^*$ is in addition strict and $p$-compatible, then there is in fact a canonical bijection between 
$$ \Hom_{\AD} (\widehat \omega^*, B^*) \cong \Hom (\ul A/p, \ul B^0/V \ul B^0)$$
where the second set denotes homomorphisms of log algebras over $\ul k$.  
\el 

\bproof 
The first assertion follows directly from Lemma \ref{lemma:homomorphism_F_compatible_algebra}. For the second assertion, we appeal to the Cartier-Dieudonn\'e-Dwork lemma, which in our setup says that, since $A$ is equipped with a lift of Frobenius $\varphi$  satisfying  $\varphi (a) \equiv a^p$ for all $a \in A,$ a homomorphism $\cl h: A \ra B^0/V B^0$ has a unique lift to a homomorphism 
$$h: A \ra W(B^0/V B^0) \xrightarrow{\mu^{-1}} B^0$$  such that $F \circ h = h \circ \varphi$. We still need to check that, given a morphism of log algebras $(\cl h, \psi)$, its lift $(h, \psi): \ul A \ra \ul B^0$ is a morphism of log algebras, namely the top square in the left diagram commutes:  
\[
\begin{tikzcd} [row sep=1.2em]
L  \arrow[d]  \arrow[r, "\psi"] & L_B \arrow[d] \\
A \arrow[rd, swap, "\cl h"] \arrow[r, dashed, "h"] & B^0 \arrow[d] \\
& B^0/V B^0
\end{tikzcd} \quad  \qquad 
\begin{tikzcd} [row sep=1.2em]
L  \arrow[dd]  \arrow[r, "\psi"] & L_B \arrow[d] \\
 & B^0/V B^0 \arrow[d, "{[\:\:]}"] \\
A \arrow[r, dashed, swap, "h"] \arrow[ru, "\cl h"] & B^0
\end{tikzcd}
\]
This is indeed the case, since the log structure on $\cl B^0$ factors through the Teichmuller lifts by Corollary \ref{cor:image_of_log}, as the diagram on the right indicates (where the outer square is the top square in the left diagram). Moreover, this also shows that any map $(h, \psi): \ul A \ra \ul B^0$ of log algebras comes from the lifting of a pair $(\cl h, \psi)$. Finally, since $p$ is $0$ in $B^0 / V B^0$, we get  $\Hom_{F} ({\ul A}, \ul B^0) = \Hom ({\ul A}/p, \ul B^0/ V \ul B^0)$ as desired. 
\eproof 

\subsubsection{Comparison with log de Rham--Witt complexes} Now let us apply the lemma above to the setup in Subsection \ref{sss:via_Frob_liftings}, namely when $\ul{\widehat A}$ comes from a log Frobenius lifting of $\ul R$. We have 
$$ \Hom_{\textup{DA}^{\log, p}} (\widehat \omega^*_{\ul{\widehat A}/W(\ul k)}, \mW \omega^*_{\ul R/\ul k}) \cong \Hom ({\ul R}, \mW_1 \omega^0_{\ul R/\ul k})$$
Let $e: R \ra \mW_1 \omega_{\ul R/\ul k}^0$ be the tautological map from the definition of the log de Rham--Witt complex $\mW \omega^*_{\ul R/\ul k}$. Then $(e, \textup{id}): {\ul R} \ra \mW_1 \omega^0_{\ul R/\ul k}$ gives rise to a map 
$$ \upsilon: \widehat \omega^*_{\ul{\widehat A}/W(\ul k)} \ra \mW \omega^*_{\ul R/\ul k}$$ in the category $\ADp$ of $p$-compatible log Dieudonn\'e $\ul W$-algebras. 

\bp \label{prop:log_dRW_another_construction}
The map $\upsilon$ induces a canonical isomorphism
$$ W \big(\widehat \omega^*_{\ul{\widehat A}/W(\ul k)} \big)_{\mathfrak{sat}} \isom \mW \omega^*_{\ul R/\ul k}. $$  
We denote this isomorphism again by $\upsilon$. 
\ep

\bproof 
By construction $W \big(\widehat \omega^*_{\ul{\widehat A}/W(\ul k)} \big)_{\mathfrak{sat}} \in \ADpc$, so it suffices to show that for any $B^* \in \ADpc$, the map $\upsilon: W \big(\widehat \omega^*_{\ul{\widehat A}/W(\ul k)} \big)_{\mathfrak{sat}} \ra \mW \omega^*_{\ul R/\ul k}$ induces a bijection between 
$$\Hom_{\ADpc} (\mW \omega^*_{\ul R/\ul k}, B^*) \stackrel{?}{=} \Hom_{\ADpc}(W \big(\widehat \omega^*_{\ul{\widehat A}/W(\ul k)} \big)_{\mathfrak{sat}} , B^*).$$ 
This follows from the following commutative diagram
\[
\begin{tikzcd}
\Hom_{\ADpc} (\mW \omega^*_{\ul R/\ul k}, B^*)  \arrow[r, "\upsilon"] \arrow[dd, equal, "\textup{Thm}. \ref{thm:log_dRW_complex_exist}"] & \Hom_{\ADpc} (W \big(\widehat \omega^*_{\ul{\widehat A}/W(\ul k)} \big)_{\mathfrak{sat}}, B^*) 
\arrow[d, equal] \\
& \Hom_{\ADp} ( \widehat \omega^*_{\ul{\widehat A}/W(\ul k)}, B^*) 
\arrow[d, equal, "\textup{Lemma } \ref{lemma:homomorphism_F_compatible_algebra_2}"]  \\
 \Hom_{\textup{Alg}^{\log}_{\ul{k}}} (\ul R, \ul{B}^0/V \ul B^0) \arrow[r, equal]  &  \Hom_{\textup{Alg}^{\log}_{\ul{k}}} (\ul R, \ul{B}^0/V \ul B^0) 
\end{tikzcd}
\]
\eproof 
Therefore, if $\ul R/ \ul k$ is integral and admits a $p$-compatible log Frobenius lifting, the saturated the log de Rham--Witt complex of $\ul R/ \ul k$ can be constructed by taking $\widehat \omega^*_{\ul{\widehat A}/W(\ul k)}$, and then taking its associated $V$-complete saturation (cf. Subsection \ref{ss:strict_completion}).

\subsection{Comparison with more general log Frobenius liftings} The goal of this Subsection is to prove a more general form of part (2) of Theorem \ref{mainthm:comparison_with_all}. In general, we cannot construct the log de Rham--Witt complexes from a non-$p$-compatible Frobenius lifting of $\ul R/\ul k$ as we did in Subsection \ref{ss:construction_via_Frob}. However, in many cases we can still compare the completed log de Rham complexes of Frobenius liftings with $\mW \omega^*_{\ul R/\ul k}$ in the derived category. This crucially relies on the notion of log-Cartier type (cf. Subsection \ref{sss:log_Cartier_type}). 
 
In this subsection, we consider the following setup 
\begin{framed} \noindent Suppose that $\ul R/\ul k$ is log-smooth of log-Cartier type (which is in particular integral), and $(\ul A_n, \varphi_n, \delta_n)$ is any log Frobenius lifting such that $\widehat A = \varprojlim A_n$ is $p$-torsion free (but the lifting may not be $p$-compatible).
\end{framed} 

\noindent By construction $\widehat A$ is $p$-complete. Moreover, each $ \omega^1_{{\ul R}_n/W_n(\ul k)}$ is a finitely generated locally free module over $A_n = \widehat A/p^n$, by log smoothness of $\ul A_n / W_n (\ul k)$. Therefore,  $\widehat \omega^1_{\ul{\widehat A}/W(\ul k)} = \varprojlim \omega^1_{{\ul R}_n/W_n(\ul k)}$ is locally free over $\widehat A$, and is in particular $p$-torsion free. Then $\ul{\widehat A}$ satisfies the condition in Proposition \ref{thm:omega_of_Witt_lifting_AD_variant}, with $\widehat \omega^1_{\ul{\widehat A}/W(\ul k)}$ replacing $\omega^1_{\ul{\widehat A}/W(\ul k)}$ in the requirement there. Following the same proof of Proposition \ref{thm:omega_of_Witt_lifting_AD_variant}, where instead of constructing $\sq \delta: \sq M \ra \omega^1_{\ul{\widehat A}/W(\ul k)}$ we only need to construct  $\sq \delta: \sq M \ra \widehat \omega^1_{\ul{\widehat A}/W(\ul k)}$, the same conclusion of Proposition \ref{thm:omega_of_Witt_lifting_AD_variant} holds. More precisely, $\widehat \omega^*_{\ul{\widehat A}/W(\ul k)}$ is equipped with a Frobenius extending $\varphi = \varprojlim \varphi_n$, which makes it a log Dieudonn\'e $W(\ul k)$-algebra,\footnote{Though we construct the log Dieudonn\'e algebra $\widehat \omega^*_{\ul{\widehat A}/W(\ul k)}$ directly using universal properties of $\widehat \omega^1_{\ul{\widehat A}/W(\ul k)}$, Lemma \ref{lemma:homomorphism_F_compatible_algebra_2} still holds by inspecting its proof (as well as the proof of Lemma \ref{lemma:homomorphism_F_compatible_algebra}).} and there is a bijection 
$$ \Hom_{\AD} (\widehat \omega^*_{\ul{\widehat A}/W(\ul k)}, \mW \omega^*_{\ul R/\ul k}) \cong \Hom (\ul R, \mW_1 \omega_{\ul R/\ul k}^0).$$
The tautological map $e: \ul R \ra \mW_1 \omega_{\ul R/\ul k}^0$ again induces a map 
$$\upsilon: \widehat \omega^*_{\ul{\widehat A}/W(\ul k)} \ra \mW \omega^*_{\ul R/\ul k}.$$ 

\br Unlike the previous subsection, in general $\upsilon$ will not induce an isomorphism from the complete saturation $W\big((\widehat \omega^*_{\ul{\widehat A}/W(\ul k)})_{\mathfrak{sat}} \big)$ of $\widehat \omega^*_{\ul{\widehat A}/W(\ul k)}$ to the saturated log de Rham--Witt complex, unless $\widehat \omega^*_{\ul{\widehat A}/W(\ul k)}$  is $p$-compatible.  
\er 

\bt \label{thm:compare_with_Frob_liftings}
Let $\ul R$ and $\widehat A$ be as in the setup above, then the map 
$$\upsilon: \widehat \omega^*_{\ul{\widehat A}/W(\ul k)} \ra \mW \omega^*_{\ul R/\ul k} $$ induces a quasi-isomorphism on the underlying cochain complexes. 
\et 

\bproof (1). First suppose that the Frobenius lift is $p$-compatible, then by Proposition  \ref{prop:log_dRW_another_construction}, we have  
$$\upsilon: \widehat \omega^*_{\ul{\widehat A}/W(\ul k)} \ra W(\widehat \omega^*_{\ul{\widehat A}/W(\ul k)})_{\mathfrak{sat}} \isom \mW \omega^*_{\ul R/\ul k}.$$
Since both $\widehat \omega^*_{\ul{\widehat A}/W(\ul k)}$ and $\mW \omega^*_{\ul R/\ul k}$ are $p$-torsion free and $p$-complete, it suffices to show that $\upsilon$ is a quasi-isomorphism after mod $p$. In other words, we need to show that 
$$ \widehat \omega^*_{\ul{\widehat A}/W(\ul k)}/p  \ra (\widehat \omega^*_{\ul{\widehat A}/W(\ul k)})_{\mathfrak{sat}} /p \ra W(\widehat \omega^*_{\ul{\widehat A}/W(\ul k)})_{\mathfrak{sat}} /p \isom \mW \omega^*_{\ul R/\ul k}/p$$
is a quasi-isomorphism. The rest of the proof is similar to the proof of Corollary \ref{cor:Cartier_implies_qi_2}. The key point is that the underlying Dieudonn\'e algebra $ \widehat \omega^*_{\ul{\widehat A}/W(\ul k)} $ satisfies the Cartier criterion (cf. Definition \ref{def:Cartier_criterion}) 
$$F: \widehat \omega^*_{\ul{\widehat A}/W(\ul k)}/p \cong \omega^*_{\ul R/\ul k} \isom H^*(\omega^{\bullet}_{\ul R/\ul k}) \cong H^*(\widehat \omega^{\bullet}_{\ul{\widehat A}/W(\ul k)}/p)$$
by Proposition \ref{prop:HK_Cartier_criterion}. To finish the proof, note that  both of the maps of cochain complexes  $(\widehat \omega^*_{\ul{\widehat A}/W(\ul k)})_{\mathfrak{sat}} /p \ra W(\widehat \omega^*_{\ul{\widehat A}/W(\ul k)})_{\mathfrak{sat}} /p$ and $ \widehat \omega^*_{\ul{\widehat A}/W(\ul k)}/p  \ra (\widehat \omega^*_{\ul{\widehat A}/W(\ul k)})_{\mathfrak{sat}} /p $ are quasi-isomorphisms by Corollary \ref{cor:Cartier_implies_qi}. 

(2). Now we treat the general case.  Again we need to show that 
$$\cl \upsilon: \widehat \omega^*_{\ul{\widehat A}/W(\ul k)}/p \cong \omega^*_{\ul R/\ul k} \xrightarrow{ \upsilon \!\!\mod p } \mW \omega^*_{\ul R/\ul k} /p $$ 
is a quasi-isomorphism.  The map $\cl \upsilon$ in fact depends on the lifting $\widehat A$ of $R$ and $\varphi$ of Frobenius, but the induced map in the derived category is independent of the lifting. More precisely, let $\ul{\widehat A'} / \ul W$ be a $p$-compatible formal lift of $\ul R/\ul k$ (which exists thanks to Lemma \ref{lemma:log_Frob_lifting_exist}),  then we consider the following diagram
\[
\begin{tikzcd} 
  \widehat \omega^*_{\ul{\widehat A}/W(\ul k)} \arrow[r] & \omega_{\ul R/\ul k}^* \arrow[r, "\cl \upsilon"]  \arrow[d, equal]&   \mW \omega^*_{\ul R/\ul k} /p  \arrow[r, "\textup{pr}"] & \mW_1 \omega^*_{\ul R /\ul k}   \arrow[d, equal] \\
  \widehat \omega^*_{\ul{\widehat A'} /W(\ul k)}  \arrow[r] & \omega_{\ul R/\ul k}^* \arrow[r, "\cl \upsilon' "] &   \mW \omega^*_{\ul R/\ul k} /p  \arrow[r, "\textup{pr}"] & \mW_1 \omega^*_{\ul R /\ul k}  
\end{tikzcd}
\]
where $\textup{pr}: \mW \omega^*_{\ul R/\ul k} /p  \twoheadrightarrow \mW \omega^*_{\ul R/\ul k} / \textup{Fil}^1 =  \mW_1 \omega^*_{\ul R/\ul k} $ is the canonical projection map. Unwinding definitions, we see that $\textup{pr} \circ \cl \upsilon = \textup{pr} \circ \cl \upsilon'$, namely that the square in the diagram commutes. From step (1), we know that $\cl \upsilon'$ is a quasi-isomorphism. By Proposition \ref{prop:comparison_with_mod_p_dR}, $ \textup{pr} \circ \cl \upsilon'$ is an isomorphism, therefore $\textup{pr}$ is a quasi-isomorphism. This in turn (using Proposition  \ref{prop:comparison_with_mod_p_dR} again) implies that $\cl \upsilon$ is a quasi-isomorphism. This concludes the theorem.  
\eproof 
 
\br  We give two remarks on the proof. 

$\bullet$ The proof of the general case (part (2)) uses Proposition \ref{prop:comparison_with_mod_p_dR}. The proof of Proposition \ref{prop:comparison_with_mod_p_dR} builds on the conclusion of the special case of the theorem where there is a $p$-compatible log Frobenius lift (i.e. part (1) in the proof), but does not use part (2) of this proof. Therefore there is no circular argument. 
 
$\bullet$ In general $\cl \upsilon$ does not agree with $\cl \upsilon'$ when $R$ is non-perfect. To see this, observe that the map $\upsilon$ comes from the algebra $f_{\varphi}: \widehat A \ra W(R)$, obtained by applying the Cartier-Dieudonn\'e-Dwork lemma (see the proof of Lemma \ref{lemma:homomorphism_F_compatible_algebra_2}). The map $f_\varphi$ is the composition $\widehat A \ra W(\widehat A) \ra W(R)$, where the first map can be determined on its ghost coordinates by 
$$\widehat A \longrightarrow  W(\widehat A) \xhookrightarrow{\textup{ gh }} {\widehat A }^{\;\N}, \qquad \sq x \mapsto (\sq x, \varphi (\sq x), \varphi^2 (\sq x), ... ).$$ 
On Witt coordinates, $f_\varphi(\sq x) = (\sq x, \theta(\sq x), ... )$ where $\theta(\sq x) = (\varphi (\sq x) - {\sq x}^p)/p \mod p \in R$. The reduction mod $p$ of $f_\varphi$ gives $\cl \upsilon|_{R}: R = \widehat A /p \ra W(R)/p$, which depends on $\widehat A $ and $\varphi$. This is not surprising since without the liftings, we cannot produce a map from $R = W(R)/V \ra W(R)/p$ when $R$ is non-perfect. 
\er

\subsection{Comparison with de Rham complexes in characteristic $p$}  

Let $\ul R/\ul k$ be a log algebra and let $\mW \omega^*_{\ul R/\ul k}$ be its saturated log de Rham--Witt complex.  We have a map $(e, \textup{id}): \ul R \ra \mW_1 \omega^0_{\ul R/\ul k}$, which extends uniquely to a map $\nu: \omega^*_{\ul R/\ul k} \ra \mW_1 \omega^*_{\ul R/\ul k}$ of cdga's that is compatible with log structures and log differentials (In fact, $\nu = \textup{pr} \circ \cl \upsilon =  \textup{pr} \circ \cl \upsilon'$ in the proof of Theorem \ref{thm:compare_with_Frob_liftings}). We have already used the following proposition to prove part (2) of Theorem \ref{thm:compare_with_Frob_liftings}.

\bp \label{prop:comparison_with_mod_p_dR}
Suppose that $\ul R /\ul k$ is log-smooth of log-Cartier type, then  
$$\nu: \omega^*_{\ul R/\ul k} \ra \mW_1 \omega^*_{\ul R/\ul k}$$ is an isomorphism of cdga's. 
\ep

\bproof 
The proof is identical to that of Proposition 4.3.2 of \cite{BLM}. Let $(\widehat A, M)$ be a p-compatible formal lift of $\ul R$ over $\ul W$, which exists by Lemma \ref{lemma:log_Frob_lifting_exist}. Then we have a commutative diagram of cochain complexes 
\[
\begin{tikzcd} [column sep=1em] 
\omega^*_{\ul R /\ul k} \arrow[rr, "C^{-1}"] \arrow[dd, "\cl \upsilon"] \arrow[rd, "\nu"] && H^*(\omega^*_{\ul R/\ul k}) = H^*(\widehat \omega^*_{\ul{\widehat A}/W(\ul k)}/p ) \arrow[dd, "H(\cl \upsilon)"]  \\
& \mW_1 \omega^*_{\ul R/\ul k} \arrow[rd, "F_1"] \\
\mW \omega^*_{\ul R/\ul k}/p \arrow[rr, "F"] \arrow[ru, "q"] && H^*(\mW \omega^*_{\ul R/\ul k}/p) 
\end{tikzcd}
\]
where both complexes on the right are equipped with the Bockstein differentials as in Subsection \ref{ss:Cartier_criterion}. The vertical maps are both induced by $\upsilon:  \widehat \omega^*_{\ul{\widehat A}/W(\ul k)} \ra \mW \omega^*_{\ul R/\ul k}$, while the bottom triangle is the triangle above Lemma \ref{lemma:F_1_is_isom_if_saturated} applied to the log de Rham--Witt complex.  

From the proof of part (1) of Theorem \ref{thm:compare_with_Frob_liftings},  $\cl \upsilon$ is a quasi-isomorphism, so $H(\cl \upsilon)$ is an isomorphism. The map $F_1$ is an isomorphism by part (2) of Lemma  \ref{lemma:F_1_is_isom_if_saturated}. Finally, by the assumption that $\ul R/\ul k$ is log-smooth of log-Cartier type, $C^{-1}$ is an isomorphism. Therefore, $\nu$ is an isomorphism. 
\eproof 

Consequently, we have
\bc \label{cor:comparison_with_mod_p_dR}
Let $\ul R /\ul k$ be log-smooth of log-Cartier type, then the projection $\mW \omega^*_{\ul R /\ul k} \ra \mW_1 \omega^*_{\ul R /\ul k} \xrightarrow{\nu^{-1}} \omega^*_{\ul R/\ul k}$ induces a quasi-isomorphism 
$$\textup{pr}: \mW \omega^*_{\ul R /\ul k} /p \mW \omega^*_{\ul R /\ul k}  \isom \omega^*_{\ul R /\ul k}$$
of cochain complexes.  
\ec 
In other words,  $\mW \omega^*_{\ul{R}/\ul{k}}$ is a deformation of $\omega^*_{\ul{R}/\ul{k}}$ in the derived category. 

\bproof 
This is immediate from the commutative diagram above, since $\cl \upsilon$ is an isomorphism. 
\eproof

\br Now suppose that $\ul R/\ul k$ satisfies the assumption in Proposition \ref{prop:comparison_with_mod_p_dR}, then it follows from $R \isom \mW_1 \omega^0_{\ul R/\ul k}$ and  by Remark \ref{remark:completion_of_A_is_saturated}  that $R$ is reduced. Therefore we obtain the follow corollary mentioned in the introduction.
\bc 
Let $\ul X$ be a fine log scheme over $\ul k$ that is log-smooth of log-Cartier type, then the underlying scheme $X$ is reduced. 
\ec

This recovers a known result in log geometry. Suppose in addition that $\ul X$ and $\spec \ul k$ are saturated log schemes as defined in \cite{Tsuji}. From \textit{loc.cit.} II.2.14 and II.3.1 (which is attributed to Kato) it follows that $\ul X$ is of log-Cartier type if and only if the morphism $\ul X \ra \ul k$ is saturated, and by Theorem II.4.2 of \textit{loc.cit.} the latter condition is equivalent to $\ul X$ being reduced.
\er

Our next remark is helpful for the discussion on monodromy operator. 

\br\label{remark:variant_of_comparison_in_char_p} 
If we look closely at the proof of Theorem \ref{thm:compare_with_Frob_liftings} and Proposition \ref{prop:comparison_with_mod_p_dR}, the assumption that $\ul R/\ul k$ is log-smooth of log-Cartier type is used to guarantee that the following two criterions hold:  
\begin{align*}
& \tag{1} \label{condition1}  \textup{There exists a lift }  \ul A/ W(\ul k) \textup{ of } \ul R \textup{ together with a lift of Frobenius}\\ &  \textup{such that } \omega^*_{\ul A/ W(\ul k)}  \textup{ is $p$-torsion free over  } W(k). \\
& \tag{2} \label{condition2}  \textup{The Cartier isomorphism } C^{-1}: \omega^{i}_{\ul R/\ul k} \isom H^i(\omega^{*}_{\ul R/\ul k}) \textup{ holds}. 
\end{align*}
Thus the conclusions of Theorem \ref{thm:compare_with_Frob_liftings}  and Proposition  \ref{prop:comparison_with_mod_p_dR} continue to hold for all log algebras $\ul R/\ul k$ that meet the two criterions above. 
\er

\subsection{$\Etale$ base change of  log de Rham--Witt complexes}  \label{sec:etale_base_change}

In this subsection, we show that the saturated log de Rham--Witt complex $\mW \omega^*_{\ul{R}/\ul{k}}$ behaves well with respect to $\etale$ base change, therefore, for any quasi-coherent log scheme $\ul{X}$ over $\ul k$, we obtain a sheaf $\mW \omega^*_{\ul X/\ul k}$ of log Dieudonn\'e algebras on the $\etale$ site $X_{\ett}$. 

To formulate the $\etale$ base change we consider the category of log cdga's, whose objects consist of $\ul A^* = (A^*, M, d, \delta)$ where $(A^*, d)$ is a cdga, $M \xrightarrow{\alpha} A^0$ a log algebra and $\delta: M \ra A^1$ a monoid morphism (with the additive monoid structure on $A^1$) satisfying $d \circ \delta = 0$ and $\alpha(m) \delta(m) = d (\alpha(m))$ for any $m \in M$. Morphisms between $\ul A^*, \ul B^*$ are pairs of maps $(f,\psi)$ where $f: A^* \ra B^*$ is a morphism of cdga's and $\psi: M_A \ra M_B$ is a morphism of monoids, compatible with $\alpha$ and $\delta$.

\bd \label{def:naively_etale}  \indent   

\bi 
\item A morphism $(f, \psi): (A, L_A) \ra (B, L_B)$ of log algebras is \textit{naively} $\etale$ if $f: A \ra B$ is $\etale$ and $\psi$ is an isomorphism of monoids. 
\item A morphism $(f,\psi): \ul A^* \ra \ul B^*$ between log cdga's is $\etale$ if its restriction on the log algebra $(f, \psi): (A^0, M_A) \ra (B^0, M_B)$ is naively $\etale$ and $f$ induces an isomorphism $A^* \otimes_{A^0} B^0 \ra B^*$ of graded algebras. 
\item If $\ul A^*, \ul B^* \in \ADc$ are strict log Dieudonn\'e algebras, then $(f, \psi): \ul A^* \ra \ul B^*$ is \textit{$V$-adically} $\etale$ if for each $n \ge 1$, the induced map of log cdga's $W_n (A^*) \ra W_n (B^*)$ is $\etale$. 
\ei
\ed 

If $\ul A^*$ is a log cdga, then we denote by $\Et_{\ul A^*}$ the category of log cdga's which are $\etale$ over $\ul A^*$. Similarly, if $\ul A^* \in \ADc$, then we denote by $\textup{V\'Et}_{\ul A^*}$ the category of $V$-adically $\etale$ $\ul A^*$-algebras in $\ADc$. Morphisms in $\Et_{\ul A^*}$ are morphisms of log cdga's (resp. morphisms in $\textup{V\'Et}_{\ul A^*}$ are morphisms of log Dieudonn\'e algebras). 

\bp   \indent 
 
(1). Let $\ul A^*$ be a log cdga, then the functor $\ul B^* \mapsto B^0$ induces an equivalence of categories between 
$$ \Et_{\ul A^*} \isom \Et_{A^0} = \{\etale A^0\textup{-algebras}\}.$$
 
(2). Let $\ul A^* \in \ADc$,  the functor that sends a V-adically $\etale$ $\ul A^*$-algebra $\ul B^* \in \ADc$ to the $ A^0/V (A^0)$-algebra $ B^0/V (B^0)$, induces an equivalence of categories 
$$ \textup{V\'Et}_{\ul A^*} \isom  \Et_{A^0/V (A^0)}. \quad \quad $$  
\ep

\bproof 
Note that the functor $\ul B^* \mapsto B^*$ forgetting the log  structure on $\ul B^*$ is an equivalence of categories between $\etale$ $\ul A^*$-algebras and $\etale$ $A^*$-algebras. Likewise $\ul B^* \mapsto B^*$ induces an equivalence between $\textup{V\'Et}_{\ul A^*}$ and $\textup{V\'Et}_{A^*}$. Therefore claim (1) and (2) follow from Proposition 5.3.2 and Theorem 5.3.4 in \cite{BLM} respectively. 
\eproof 

By the same proof of Corollary 5.3.5 of \cite{BLM}, we arrive at
\bc \label{cor:etale_bc}
Let $\ul R = (R, M)$ be a log algebra over $\ul k$, and let $R \ra S$ be an $\etale$ morphism of $k$-algebras. Let $\ul S = (S, M)$ be the naively $\etale$ algebra over $\ul R$. Then for any $n \ge 1$, the map $\mW_n \omega^*_{\ul R/\ul k} \ra \mW_n \omega^*_{\ul S/\ul k}$ is $\etale$. In other words, $\mW \omega^*_{\ul R/\ul k} \ra \mW \omega^*_{\ul S/\ul k}$ is V-adically $\etale$. Moreover, we have the following natural isomorphism of graded algebras
$$\mW_n \omega^*_{\ul R/\ul k} \otimes_{W_n (R)} W_n (S) \isom \mW_n \omega^*_{\ul S/\ul k}.$$ 
\ec

\br[Other approaches for $\etale$ base change] In the special case when $\ul R$ is log-smooth over $\ul k$ of log-Cartier type, we can also deduce the base change isomorphism $\mW_n \omega^*_{\ul R/\ul k} \otimes_{W_n (R)} W_n (S) \isom \mW_n \omega^*_{\ul S/\ul k}$ by comparing to the log de Rham--Witt complexes constructed by Hyodo--Kato (which is defined globally on the $\etale$ site) or Matsuue (using Proposition 3.6 in \cite{Ma}), once we prove the comparison theorems for such log algebras (Propositions \ref{thm:compare_with_HK_affine} and x\ref{prop:compare_to_Matsuue}). 
\er

Now for any quasi-coherent log scheme $\ul X$ over $\spec \ul k$, we construct the sheaf of log de Rham--Witt complexes $\mW \omega^*_{\ul X/\ul k}$ on $X_{\ett}$.  

\bd \label{def:X_affine_with_chart} Let $\ul X = (X, M_X)$ be a log scheme over $\ul k$. Define $X_{\ett, \textup{aff}}$ to be the site of ``small enough'' affine $\etale$ opens where the log structure of $\ul X$ admits charts. More precisely, the objects of  $X_{\ett, \textup{aff}}$ are $\etale$ morphisms $U \xrightarrow{h} X$ over $k$ where $U = \spec R$ is affine, such that there exists a constant log structure $L \ra R$ and an isomorphism $(\mL_U)^a \isom M_X|_{U} := h^* M_X$. Here $(\mL_U)^a $ is the associated log structure of the constant pre-log structure $L_U \ra \mO_{U}$. The topology on $X_{\ett, \textup{aff}}$ is given by $\etale$ coverings. 
\ed 
It is straightforward to check that $X_{\ett, \textup{aff}}$ indeed forms a site. Moreover, if the log structure $M_X$ admits charts $\etale$ locally, then these small affine opens form a basis for the $\etale$ topology on $X$. 

\bl \label{lemma:affine_spans_topology} Suppose that $\ul X$ is a quasi-coherent log scheme over $\ul k$, then the restriction of an $\etale$ sheaf from $X_{\ett}$ to $X_{\ett, \text{aff}}$ induces an equivalence of topoi. 
\el

\bproof 
The canonical functor $u: X_{\ett, \text{aff}} \longrightarrow X_{\ett}$ is clearly continuous, cocontinuous, and satisfies the property that any $\etale$ open $Y \in X_{\ett}$ admits a cover by $\{U_i\}$ in $X_{\ett, \text{aff}}$, as the log structure $M_X$ is quasi-coherent by assumption. Therefore $u$ induces an equivalence $\textup{Sh}(X_{\ett, \text{aff}}) \isom \textup{Sh} (X_{\ett})$\footnote{For example by \cite{SP} 03A0}.
\eproof

\bt \label{theorem:global_dRW}
Let $\ul X = (X, M_X)$ be a quasi-coherent log scheme over $\ul k$. There exists a unique sheaf $\mW \omega^*_{\ul X/\ul k}$ on $X_{\ett}$ valued in log Dieudonn\'e algebras, such that on each  $\etale$ local chart $ (U, L) = \spec (R, L)$ of $\ul X$, there is a canonical map of log Dieudonn\'e algebras  
$$ \mW \omega^*_{\ul R/\ul k} \longrightarrow \Gamma (U, \mW \omega^*_{\ul X/\ul k})$$ 
given by an isomorphism on the underlying Dieudonn\'e algebras, and $L \ra \Gamma (U, M_X)$ on the log structures. 
\et

\bproof 
We will construct $\mW \omega^*_{\ul X/\ul k}$ as the inverse limit of the tower
$$ \cdots \ra \mW_n \omega^*_{\ul X/\ul k} \ra \mW_{n-1} \omega^*_{\ul X/\ul k} \ra \cdots \ra \mW_1 \omega^*_{\ul X/\ul k},$$ so it suffices to construct the sheaf of abelian groups $\mW_n \omega^i_{\ul X/\ul k}$ for each $n\ge 1, i \ge 0$.  By Lemma \ref{lemma:affine_spans_topology} it suffices to construct each $\mW_n \omega^*_{\ul X/\ul k}$ on $X_{\ett, \textup{aff}}$. 

Now we define a presheaf of abelian groups
$\mW_n \omega^i_{\ul X/\ul k}: X_{\ett, \textup{aff}} \ra \textup{Ab}$  as follows.  For each  $U = \spec R \ra X$ in $X_{\ett, \textup{aff}},$ we write $M(U) :=  \Gamma (U, M_X)$, and define
$$\mW_n \omega^i_{\ul X/\ul k} (U \ra X) := \mW_n \omega^i_{ (R, M(U))/\ul k}.$$
This is clearly a presheaf.\footnote{We remark that the value of $\mW_n \omega^i_{\ul X/\ul k} (U \ra X)$ is independent of the choice of chart on $U$, by Proposition \ref{prop:choice_of_pre_log_structures}. In fact, the constant monoid $M(U)$ gives a (functorial) choice of chart on each $U \in X_{\ett, \textup{aff}}$, by Lemma \ref{lemma:pushout_agrees}.}  We need to show that this presheaf is a sheaf on $X_{\ett, \textup{aff}}$. For this it suffices to check that, on each $\etale$ opens $U = \spec R \in X_{\ett, \textup{aff}}$, the restriction $W_n \omega^i_{\ul X/\ul k} \;\vline_{ U_{\ett, \textup{aff}} } $ defines a sheaf on $U_{\ett, \textup{aff}} = X_{\ett, \textup{aff}} /U$, which consists precisely of all $\etale$ $R$-algebras. In other words, we may assume that $X = \spec R$ is a ``small enough'' affine log scheme which admits a chart for its log structure.  Now let $R \ra S$ be an $\etale$ map, write $V = \spec S$, then we have morphisms 
$$ \mW_n \omega^i_{(R, M(X))/\ul k} \xrightarrow{\gamma_1} \mW_n \omega^i_{(S, M(X))/\ul k} \xrightarrow{\gamma_2} \mW_n \omega^i_{(S, M(V))/\ul k}$$
by functoriality (of saturated log de Rham--Witt complexes). $\gamma_2$ is an isomorphism by Lemma \ref{lemma:etale_chart} and Proposition \ref{prop:choice_of_pre_log_structures}. By Corollary \ref{cor:etale_bc}, we know that 
$$ \mW_n \omega^i_{(S, M(V))/\ul k} \cong  \mW_n \omega^i_{(R, M(X))/\ul k} \otimes_{W_n(R)} W_n(S).$$
Now by the invariance of $\etale$ sites under nilpotent thickenings and Theorem 5.4.1 of \cite{BLM}, we know that the association $S \mapsto W_n (S)$ identifies the $\etale$ algebras over $R$ with $\etale$ algebras over $W_n(R)$. Therefore, on the affine $\etale$ site $\big(\spec W_n(R) \big)_{\ett, \textup{affine}}$, the presheaf $\mW_n \omega^i_{\ul X/\ul k}$ is in fact a sheaf associated to the $W_n(R)$-module $\mW_n \omega^i_{(R, M(X))/\ul k}$, by standard $\etale$ descent. It is straightforward to keep track of differentials, Frobenius operators, as well as the monoid morphisms $M_X \ra \mW \omega^0_{\ul X/\ul k}$ and $M_X \xrightarrow{\delta} \mW \omega^1_{\ul X/\ul k}$, hence the theorem follows. 
\eproof

\br \label{remark:direct_global_version}

It might be possible to develop the global theory directly by working with sheaves of log Dieudonn\'e algebras over $X_{\ett}$, by mimicking the steps used in the proof of Theorem \ref{thm:log_dRW_complex_exist}. However, we prefer to work locally and then globalize in the last step. Otherwise it seems that one needs to sheafify at each step and it becomes more involved to describe global sections.  In fact, for our application to $\Ainf$-cohomology of $\fX$ (see Introduction for the setup),  it is more convenient to construct the saturated log de Rham--Witt complexes in local charts.  For example, as discussed in the introduction, the log structures do not appear so transparently on $A \Omega_{\fX} \otimes^\L W(k)$, but we could get around this problem locally by choosing a certain coordinate. Moreover, the formal log scheme $\ul \fX$ is strictly speaking not log-smooth over the base $\ul{\mO_C}$, since the log structures are not coherent (finitely generated). This will not become an issue if we restrict to small enough charts, where we could pretend that $\ul \fX$ is log-smooth over $\ul{\mO_C}$ since all the information of the log structures over these charts are contained in a discrete subring $\mO \subset \mO_C$.  For the detail of these discussions see Subsection \ref{ss:semistable_setup}, \ref{ss:log_structure_after_coordinates}, and Remark \ref{remark:pullback_log_structure_1}. 

\er 

\newpage

\section{Comparison theorems} \label{section:comparison} In this section we show that our log de Rham--Witt complexes agree with the construction of Hyodo--Kato and Matsuue for log-smooth schemes of log-Cartier type. In particular, in this entire section our log schemes ($\ul X$ and $\ul Y$) are assumed to be log-smooth of log-Cartier type over $\ul k = (k, N)$, where $N \minus \{0\} \mapsto 0$, in particular, all log structures are assumed to be fine (and hence quasi-coherent) in this Section. 

\subsection{Comparison with the constructions of Hyodo--Kato}  \label{sss:construction_HK} 

In \cite{HK}, Hyodo--Kato constructs a log de Rham--Witt complex $W^{\textup{\tiny HK}} \omega^*_{\ul X/\ul k}$ using the log crystalline site, as a generalization of \cite{RI} \footnote{In particular, their construction is global, while our log de Rham--Witt complex (as well as Matsuue's version, reviewed in the next subsection) is constructed $\etale$ locally and then globalized.}.  Let $\ul X/\ul k$ be a log-smooth scheme over $\ul k$ of log-Cartier type, and let $(\ul X/\ul{W_n})_{\log\textup{-cris}}$ be the log crystalline site of $\ul X$ over $\ul{W_n}$, where $\ul{W_n}  = W_n(\ul k)$ is equipped with the standard PD structure.  
Let  $$u^{\log}_{\ul X/W_n}: (\ul X/\ul{W_n})_{\log\textup{-cris}} \ra X_{\ett} $$  be the canonical projection to the $\etale$ site. 
Now define  
$$ W_n^{\textup{\tiny HK}} \omega^i_{\ul X/\ul k} := R^i u^{\log}_{\ul X/W_n, *} (\mO_{\ul X/\ul{W_n}}) $$
where $\mO_{\ul X/\ul{W_n}}$ is the structure sheaf on the log crystalline site.
These objects are equipped with (collections of) operators $d, F, V$, and the canonical projections $ R$, which we now briefly recall (for detail, see 4.1 and 4.2 of \cite{HK}). Let $C_n^*$ be a crystalline complex (for a chosen embedding system, cf. 2.18-2.19 of  \textit{loc.cit}), which computes $ R u^{\log}_{\ul X/W_n, *} (\mO_{\ul X/\ul{W_n}})$ (cf. Proposition 2.20 of  \textit{loc.cit}),  
\bi 
\item Define $d: W_n^{\textup{\tiny HK}} \omega^i_{\ul X/\ul k} \ra W_n^{\textup{\tiny HK}} \omega^{i+1}_{\ul X/\ul k}$ to be the ``Bockstein differential'' induced from the exact sequence $0 \ra C_n^* \xrightarrow{p^n} C_{2n}^* \ra C_{n}^* \ra 0$.
\item Let $F: W_{n}^{\textup{\tiny HK}} \omega^i_{\ul X/\ul k} \ra W_{n-1}^{\textup{\tiny HK}} \omega^i_{\ul X/\ul k}$ be the map induced by $C_{n}^* \ra C_{n-1}^*$.
\item Let $V: W_{n-1}^{\textup{\tiny HK}} \omega^i_{\ul X/\ul k} \ra W_{n}^{\textup{\tiny HK}} \omega^i_{\ul X/\ul k}$ be the map induced by $C_{n-1}^* \xrightarrow{p} C_{n}^*$. 
\item Define $R_n: W_{n}^{\textup{\tiny HK}} \omega^i_{\ul X/\ul k} \ra W_{n-1}^{\textup{\tiny HK}} \omega^i_{\ul X/\ul k}$ as follows. Let $\mu_n: H^i(C_n^*) \mapsto H^i(\eta_p C_{n-1}^*)$ be the map induced by $x \mapsto p^i x$ on the level of cochain complexes \footnote{The complex ``$\eta_p C_{n-1}^*$'' needs justification, since $C_{n-1}^*$ is $p^{n-1}$-torsion. This is in fact the complex denoted by $E_{n-1}^*$ in \cite{HK}, where $E_{n-1}^i$ is defined to be ``$(\eta_p C_m)^i/p^{n-1}$'' for $m$ sufficiently large,  where $\eta_p C_m^i = \{x \in p^i C_{m}^i \: : \: dx \in p^{i+1} C_{m}^{i+1}\}$ as usual. For $m > n + i + 1$, the issue of $p$-torsion goes away and the quotient $(\eta_p C_m)^i/p^{n-1}$ is well-defined.}, 
and let $\psi: H^i(C_{n-1}^*) \isom H^i (\eta_p C_{n-1}^*)$ be the isomorphism given by Lemma 2.25 of \cite{HK}, then define $R_n : = \mu_n \circ \psi^{-1}$. 
\ei  
As mentioned in Subsection \ref{ss:L_eta_fixed_point}, this is essentially the same construction as the inverse functor in Remark \ref{remark:BLM_L_eta_p_fixed_points}. We then take the inverse limit along the restriction maps to form the  Hyodo--Kato complex $W^{\textup{\tiny HK}} \omega^*_{\ul X/\ul k} :=  \varprojlim_{R} W_n^{\textup{\tiny HK}} \omega^*_{\ul X/\ul k}$. 

\bt \label{thm:compare_with_HK}
Let $\ul X$ be a log scheme over $\ul k$ that is log-smooth of log-Cartier type. Then there is a natural isomorphism of sheaves of complexes $$ \mW \omega^*_{\ul X/\ul k} \isom  W^{\textup{\tiny HK}} \omega^*_{\ul X/\ul k}$$
compatible with Frobenius operators on both sides, which can be upgraded to an isomorphism of sheaves of log Dieudonn\'e algebras. 
\et

\br 
In particular, $\mW \omega^*_{\ul X/\ul k}$ computes the log crystalline cohomology for quasi-coherent log-smooth schemes of log-Cartier type.
\er

The proof of the Theorem will be given in the end of this Subsection. We first restrict to the affine case. Let $\ul Y = (\spec R, M)$ be an affine log scheme over $\ul k$, and consider the Hyodo--Kato complex in the affine setup
$$W^{\textup{\tiny HK}} \omega^*_{\ul R/\ul k} := \Gamma(Y, W^{\textup{\tiny HK}} \omega^*_{\ul Y/\ul k}) = \varprojlim_{R_n}  \Gamma(Y, W_n^{\textup{\tiny HK}} \omega^*_{\ul Y/\ul k}).$$  Unwinding definitions, we see that the canonical projections $R_n$ commute with operators $d$, $F$ and $V$ on each finite level. Therefore we obtain a cochain complex in the inverse limit, equipped with $F$ and $V$. By relations given in 4.1.1 of \cite{HK}, $W^{\textup{\tiny HK}} \omega^*_{\ul R/\ul k}$ is a Dieudonn\'e complex (cf. Subsection  \ref{ss:L_eta_fixed_point}).

\bl  \label{lemma:HK_complex_is_ADc_1}
The Frobenius $F$ on $W^{\textup{\tiny HK}} \omega^*_{\ul R/\ul k}$ is a graded ring homomorphism. 
\el
This is of course of no surprise and probably well-known to experts, but we cannot track down an explicit proof, so we record the proof here. 

\bproof 
For notational simplicity we will write $W^{\textup{\tiny HK}, *}$ (resp. $W_n^{\textup{\tiny HK}, *}$) for the complex $W^{\textup{\tiny HK}} \omega^*_{\ul R/\ul k}$ (resp. $W_n^{\textup{\tiny HK}} \omega^*_{\ul R/\ul k} := \Gamma(Y, W_n^{\textup{\tiny HK}} \omega^*_{\ul Y/\ul k})$) in this proof. We need to show that $F: W_n^{\textup{\tiny HK}, *} \ra W_{n-1}^{\textup{\tiny HK}, *}$ is a graded ring homomorphism.  By Proposition 4.7 of \cite{HK} and its proof,  there is a canonical surjection
$$\psi: \omega^*_{W_n(\ul R)/W_n(\ul k)} \twoheadrightarrow W_n^{\textup{\tiny HK}, *}$$ with kernel $I_n \subset  \omega^*_{W_n(\ul R)/W_n(\ul k)} $ a graded differential ideal generated by elements of the form $\{\eta_{i, j, a, m}, d \eta_{i, j, a, m} \}$ for $0 \le j \le i < n, a \in R, m \in M$, where 
$$ \eta_{i, j, a, m}:= V^i\big([a]\big) d V^j \big([\alpha(m)]\big) - V^i \big([a \alpha(m)^{p^{i-j}}]\big) \textup{d}\log (m).$$
To further ease notation, we write $\omega^*_{W_n(\ul R)} := \omega^*_{W_n(\ul R)/W_n(\ul k)} $ in the rest of the proof.  
The map $\psi$ is compatible under the canonical projection maps. By (the correction in Section 7 in \cite{Na} of) 4.9.1 of \cite{HK}, we have the following isomorphism of cochain complexes
$$s= \psi \circ t: (W_n^{\textup{\tiny HK}, *})'' \xrightarrow{ \:\: t \: \:}  \omega^*_{W_n(\ul R)}/I_n  \xrightarrow{\:\: \psi\:\: } W_n^{\textup{\tiny HK}, *}$$
where $(W^{\textup{\tiny HK}, i})''$ is defined as a quotient of \footnote{The complex $(W^{\textup{\tiny HK}, i})'$ which appears in Proposition 4.6 of \cite{HK} is incorrect (and the proof given there is incomplete), instead we should take the modification $(W^{\textup{\tiny HK}, i})''$ given in Section 7 of \cite{Na}. With this correction, the map $s$ is indeed an isomorphism, by Theorem 7.5 of \textit{loc. cit}.} 
$$ (W_n (R) \otimes \wedge^{i} M^{\textup{gp}}/N^{\textup{gp}}) \oplus  (W_n (R) \otimes \wedge^{i} M^{\textup{gp}}/N^{\textup{gp}}).$$ 

By Remark \ref{remark:Frob_finite_level}, we have a graded ring homomorphism  $F: \omega^*_{W_n(\ul R)} \ra \omega^*_{W_{n-1}(\ul R)}.$
From the construction given there, it sends the ideal $I_{n}$ into $I_{n-1}$, therefore inducing a graded ring homomorphism 
$$F: \omega^*_{W_n(\ul R)}/I_n \ra \omega^*_{W_{n-1}(\ul R)}/I_{n-1}.$$
We claim that the following diagram (on the right) commutes 
\[
\begin{tikzcd}[column sep=1.2em]
 (W_n^{\textup{\tiny HK}, *})'' \arrow[r, "t"] \arrow[d, "F"] & \omega^*_{W_n(\ul R)}/I_n \arrow[r, "\psi"] \arrow[d, "F"] & W_n^{\textup{\tiny HK}, *}\arrow[d, "F"] \\
  (W_{n-1}^{\textup{\tiny HK}, *})'' \arrow[r, "t"] & \omega^*_{W_{n-1}(\ul R)}/I_{n-1} \arrow[r, "\psi"] & W_{n-1}^{\textup{\tiny HK}, *} 
\end{tikzcd}
\]
Where $F: (W_{n}^{\textup{\tiny HK}, i})'' \ra (W_{n-1}^{\textup{\tiny HK}, i})''$ is defined as follows: for any $x \in W_n(R)$ and $y \in W_n(R)^\times$, write $y = [y_0] + V y'$, and send
\begin{align*}
(x \otimes m_1 \wedge \cdots m_i, 0) & \mapsto (F(x) \otimes m_1 \wedge \cdots m_i, 0)   \\
(0, y \otimes m_2 \wedge \cdots m_i) & \mapsto ([y_0]^p \otimes \beta^{-1}([y_0]) \wedge m_2 \wedge \cdots \wedge m_i, 0) \\
& \qquad \qquad \qquad  \qquad \qquad  + (0, y' \otimes m_2 \wedge \cdots \wedge m_i).
\end{align*}
$\beta$ here denotes the log  structure $M \xrightarrow{\alpha} R \xrightarrow{[\;]} W(R)$, where the first arrow $\alpha$ is a log structure by assumption, so there exists a unique element $\beta^{-1} ([y_0]) \in M$ which is sent to $[y_0]$. 

The commutativity of the right square follows from the commutativity of the outer square (by Remark 7.7.4 of \cite{Na}) and the commutativity of the left square, which we now check by hand on the generators of $ (W_{n}^{\textup{\tiny HK}, i})''$. Clearly $F \circ t = t \circ F$ on $(x \otimes m_1 \wedge \cdots m_i, 0)$; for $(0, y \otimes m_2 \wedge \cdots m_i) $, we simply observe that 
$$ F d y = F d ([y_0] + Vy') = [y_0]^p \textup{d}\log [y_0] + d y'.$$
This finishes the proof. 
\eproof 

\bl  \label{lemma:HK_complex_is_ADc_2}
$W^{\textup{\tiny HK}} \omega^*_{\ul R/\ul k}$ is a saturated Dieudonn\'e algebra. 
\el

\bproof 
Keep the notations from the proof of Lemma \ref{lemma:HK_complex_is_ADc_1}. 
By Corollary 4.5 (2) of \cite{HK}, $W^{\textup{\tiny HK}} \omega^*$ is $p$-torsion free. Since $\ul R/\ul k$ is log-smooth of log-Cartier type, we have the Cartier isomorphism (cf. Proposition \ref{prop:HK_Cartier_criterion}), which implies (by the same proof of I.3.11.4 of \cite{Illusie} in the case $n = 0$) that 
$$ \ker (d: W^{\textup{\tiny HK}}_1 \omega^i \ra W^{\textup{\tiny HK}}_1 \omega^{i+1}) = \im \big(F: W^{\textup{\tiny HK}}_2\omega^i \ra W^{\textup{\tiny HK}}_1 \omega^i  \big).$$
This in turn implies that 
$$d^{-1} (p \cdot W^{\textup{\tiny HK}} \omega^*) = F \big(W^{\textup{\tiny HK}} \omega^* \big).$$
To see this, suppose that $dx \in p  \cdot W^{\textup{\tiny HK}} \omega^*$, then by the equality above we have $x = F x_1' + V x_2' + d V y_2' = F x_1 + V x_2'$ where $x_1 = x_1' + d V^2 (y_2')$. This implies that $d V x_2' \in p \cdot W^{\textup{\tiny HK}} \omega^*$, so $d x_2' = F d V x_2' \in p \cdot W^{\textup{\tiny HK}} \omega^*$ and we may write $x_2' = F x_2 + V x_3'$. By repeating this procedure we may write $x = F(x_1 + V x_2 + V^2 x_3 + \cdots)$. This together with the injectivity of $F$ (since $W^{\textup{\tiny HK}} \omega^*$ is $p$-torsion free) implies that $W^{\textup{\tiny HK}} \omega^*$ is saturated (as a Dieudonn\'e complex).  

To finish the proof of the Lemma, it remains to show that $F (x) \equiv x^p \mod p$ for all $x \in W^{\textup{\tiny HK}} \omega^0$. Since $W^{\textup{\tiny HK}} \omega^*$ is saturated, in fact it suffices to check that $F (x) \equiv x^p \mod V$ for all $x \in W^{\textup{\tiny HK}} \omega^0$.  But this follows from the commutative diagram in the proof of Lemma \ref{lemma:HK_complex_is_ADc_1}, since the middle vertical map is constructed to extend the Witt vector Frobenius and the map $\psi$ preserves the $V$-filtration (by Proposition 7.1 and diagram 7.5.6 in \cite{Na}). 
\eproof

\bc \label{prop:HK_complex_is_ADc}
$W^{\textup{\tiny HK}} \omega^*_{\ul R/\ul k}$ admits a natural structure of a log Dieudonn\'e algebra, which is strict (in particular saturated) and $p$-compatible.  
\ec 

\bproof 
For each $n$, the log structure $[\alpha]_n: M \xrightarrow{\alpha} R \xrightarrow{[\;]} W_n(R)$ induces a log derivation $(d, \delta_n: M \ra \omega^1_{W_n(\ul R)})$. By composing with projection we get a monoid map $\delta_n: M \ra \omega^1_{W_n(\ul R)}/I_n$. Both $[\alpha_n]$ and $\delta_n$ are compatible under the canonical projections. Now the map $\psi$ in  
$$(W_n^{\textup{\tiny HK}, *})'' \xrightarrow{ \:\: t \: \:}  \omega^*_{W_n(\ul R)}/I_n  \xrightarrow{\:\: \psi\:\: } W_n^{\textup{\tiny HK}, *}$$
(from the proof of Lemma \ref{lemma:HK_complex_is_ADc_1}) is also compatible with the canonical projections. Therefore we have monoid morphisms (which we still denote by the same symbols): $ [\alpha]_n: M \ra W_n^{\textup{\tiny HK}} \omega^0 $ and $\delta_n:  M \ra W_n^{\textup{\tiny HK}} \omega^1.$
In the inverse limit they give rise to 
$$[\alpha]:  M \ra W^{\textup{\tiny HK}} \omega^0,  \quad \delta:  M \ra W^{\textup{\tiny HK}} \omega^1.$$
It is straightforward to check that this makes $W^{\textup{\tiny HK}} \omega^*$ a $p$-compatible log Dieudonn\'e algebra. It is saturated by Lemma \ref{lemma:HK_complex_is_ADc_2};  and it is strict since for each $n \ge 1$ and any $k \ge 0$, we have the following exact sequence 
$$ 0 \ra \im \Big(W_{k}^{\textup{\tiny HK}} \omega^* \xrightarrow{V^n + d V^n} W_{n+k}^{\textup{\tiny HK}} \omega^* \Big) \longrightarrow W_{n+k}^{\textup{\tiny HK}} \omega^*  \xrightarrow{\textup{pr}} W_n^{\textup{\tiny HK}} \omega^* \ra 0 $$
by Section 4.9 of \cite{HK}. 
\eproof 

With the preparations above, the identity map 
$$ \ul R \longrightarrow W_1^{\textup{\tiny HK}} \omega^0_{\ul R/\ul k} \cong \ul R $$
induces a map of strict log Dieudonn\'e algebras 
$\gamma: \mW \omega^*_{\ul R/\ul k} \longrightarrow W^{\textup{\tiny HK}} \omega^*_{\ul R/\ul k}.$ This essentially proves Theorem \ref{thm:compare_with_HK} in the affine case. 

\bp \label{thm:compare_with_HK_affine}
$\gamma$ is an isomorphism of log Dieudonn\'e algebras 
$$\gamma: \mW \omega^*_{\ul R/\ul k} \isom W^{\textup{\tiny HK}} \omega^*_{\ul R/\ul k}. $$ 
\ep 

\bproof 
The induced map $\cl \gamma: \mW_1 \omega^*_{\ul R/\ul k}  \ra W_1^{\textup{\tiny HK}} \omega^*_{\ul R/\ul k}$ is an isomorphism by Proposition \ref{prop:comparison_with_mod_p_dR}. The theorem then  follows from Corollary \ref{cor:Cartier_implies_qi} part (3). 
\eproof 

\bproof[Proof of Theorem \ref{thm:compare_with_HK}] 
The construction on affine log schemes gives rise to an isomorphism of sheaves $\gamma: \mW \omega^*_{\ul X/\ul k} \isom W^{\textup{\tiny HK}} \omega^*_{\ul X/\ul k}$ on the site $X_{\ett, \textup{aff}}$ (see Definition \ref{def:X_affine_with_chart}). The functorality follows from the commutative diagram  
\[
\begin{tikzcd}[row sep=1.2em]
R \arrow[r] \arrow[d] & W_1^{\textup{\tiny HK}} \omega^0_{\ul R/\ul k} \arrow[d] \\
S \arrow[r]  & W_1^{\textup{\tiny HK}} \omega^0_{\ul S/\ul k}
\end{tikzcd}
\]
for any morphism for any $\etale$ morphism $R \ra S$, where $U = \spec R \in X_{\ett, \textup{aff}}$. 
\eproof

\subsection{Comparison with the constructions of Matsuue}  \label{sss:construction_Matsuue} 

We start by recalling the definition of the category $\mC_{F\!V}$ of (R-framed) log F-V-procomplexes.  
\bd[\cite{Ma} Definition 3.4]  \label{def:log_pro_F_V} 
Let $\ul R = (R, M)$ be a log algebra over $\ul k = (k, N)$. A $R$-framed F-V-procomplex over $\ul R/\ul k$ is a projective system 
$$ \cdots \ra E_{m+1}^* \xrightarrow{R_m} E_{m}^* \xrightarrow{R_{m\text{-}1}} \cdots \ra E_1^* \ra E_0^* = 0 $$ where each $E_m^* = (E_m^*, d_m, \delta_m)$ is a log cdga over $W_m(\ul R)/W_m (\ul k)$, 
together with 
\bi
\item a collection $ F: E_{m+1}^* \ra E_{m}^*$ of graded ring homomorphisms;
\item a collection $V: E_m^* \ra E_{m+1}^*$ of graded abelian group homomorphisms. 
\ei
These data are required to satisfy the following conditions
\be[(a)]
\item $R=\{R_m\}$ is compatible with $\delta =\{\delta_m\}$, i.e. $\delta_m = R_m \circ \delta_{m+1}$ for $m \ge 0$. 
\item $R$ is compatible with the collection of maps $F$ and $V$.
\item The structure maps $\beta_m: W_m (R) \ra E_m^0$ are compatible with $F$ and $V$. 
\item The collections of $F, V, d, \delta$ satisfies 
\be[i)]
\item $F V = p$
\item $F d_{m+1} V = d_m$
\item $F d_{m+1} [x] = [x]^{p-1} d_m [x]$
\item $(Vx)y = V(x Fy)$
\item $F \delta_{m+1} = \delta_{m}$. 
\ee
\ee
Here by the structure maps, we mean the collection of ring homomorphisms  
$$\beta_m: W_m(R) \ra E_m^0$$ 
given as part of the data of the $W_m(\ul R)/W_m (\ul k)$-cdga $(E_m^*, d_m, \delta_m)$ where $(d_m \circ \beta_m, \delta_m: M \ra E_m^1)$ forms a log derivation of $W_m(\ul R)/W_m(\ul k)$ into $E^1$ and satisfies $d_m \circ \delta_m = 0$. 
\ed 

Denote the category of log F-V-procomplexes over $\ul R/\ul k$ by  $\mC_{FV, \ul R}$. In \cite{Ma}, Matsuue constructs an initial object $\{W_m \Lambda^*_{\ul R/\ul k}\}$ in $\mC_{FV, \ul R}$, where each $W_m \Lambda^*_{\ul R/\ul k}$ is constructed as a certain quotient of $\omega^*_{W_m(\ul R)/W_m(\ul k)}$. The complexes $ W \Lambda^*_{\ul R /\ul k } = \lim_{n} W_n \Lambda^*_{\ul R /\ul k }$ is then globalized to a sheaf $W \Lambda^*_{\ul X/\ul k }$ on $X_{\ett}$.  In fact, the construction in \cite{Ma} works more generally over non-perfect bases, but we have restricted to the case of perfect base (in fact a perfect field) to compare with our version of saturated log de Rham--Witt complexes. The following lemma is clear. 

\bl 
The log de Rham--Witt complexes $\mW \omega^*_{\ul R/\ul k}$  gives rise to an $R$-framed F-V-procomplex $\{W_m \omega^*_{\ul R/\ul k}\}$ over $\ul R/\ul k$. 
\el

Therefore by the universal property of $W \Lambda^*$ we have a map $$ \gamma': W \Lambda^*_{\ul R /\ul k } \longrightarrow \mW \omega^*_{\ul R/\ul k}$$ of cochain complexes which is compatible with $F, V$ and $R$. 

\bp \label{prop:compare_to_Matsuue}
$\gamma'$ is an isomorphism (under the assumption that $\ul R/\ul k$ is log-smooth of log-Cartier type). 
\ep 

\bproof 
The map $\gamma'$ induces commutative diagrams
\[
\begin{tikzcd}
0 \arrow[r]   & \textup{Fil}^n W_{n+1} \Lambda^*_{\ul R/\ul k} \arrow[r] \arrow[d]  & W_{n+1} \Lambda^*_{\ul R/\ul k}  \arrow[r]  \arrow[d]  & W_n \Lambda^*_{\ul R/\ul k} \arrow[r] \arrow[d] & 0 \\ 
0 \arrow[r]   & \textup{Fil}^n \mW_{n+1} \omega^*_{\ul R/\ul k} \arrow[r]   & \mW_{n+1} \omega^*_{\ul R/\ul k}  \arrow[r]   & \mW_n \omega^*_{\ul R/\ul k} \arrow[r] & 0 
\end{tikzcd} 
\]
between short exact sequences, where the filtrations are standard filtrations, given by $\textup{Fil}^n W_{n+1} \Lambda^*_{\ul R/\ul k} = V^n W_1 \Lambda^*_{\ul R/\ul k} +  \textup{d} V^n W_1 \Lambda^{*-1}_{\ul R/\ul k}$ and likewise for $ \mW_{n+1} \omega^*_{\ul R/\ul k}$. By Proposition  \ref{prop:comparison_with_mod_p_dR}, we know that $\gamma'$ induces an isomorphism $W_1 \Lambda^*_{\ul R/\ul k} \isom \mW_1 \omega^*_{\ul R/\ul k} $ and hence $\textup{Fil}^n W_{n+1} \Lambda^*_{\ul R/\ul k}  \isom \textup{Fil}^n \mW_{n+1} \omega^*_{\ul R/\ul k} $ for each $n$. The proposition therefore follows by induction. 
\eproof 

As in Subsection \ref{sss:construction_HK}, we have the following 
\bc For a quasi-coherent log scheme $\ul X$ which is log-smooth of log-Cartier type over $\ul k$, the isomorphism $\gamma'$ globalizes to an isomorphism of sheaves $\gamma':  W \Lambda^*_{\ul X /\ul k } \isom \mW \omega^*_{\ul X/\ul k}$.
\ec

\br  \label{remark:GL_comparison}
One possible way to obtain $\mW \omega_{\ul R/\ul k}^* \isom W \Lambda^*_{\ul R/\ul k}$ is to use the comparison with Hyodo--Kato complexes in Subsection \ref{sss:construction_HK} and prove that the Hyodo--Kato complexes agree with $W \Lambda_{\ul R/\ul k}^*$ constructed by Matsuue. The latter is claimed in Proposition 3.1 of \cite{GL}, which asserts that there is a natural isomorphism 
$$ W \Lambda^*_{\ul R/\ul k} \isom W^{\textup{\tiny HK}} \omega^*_{\ul R/\ul k}$$ 
identifying the log  structures, $F, V$ and $R$ on both sides. The key assertion used in their proof is that $W_m^{\tiny \textup{HK}} \omega^*_{\ul R/\ul k}$ form an $R$-framed F-V-procomplex, which is claimed without justification. It seems to the author that this would involve proving most of the statements in Subsection \ref{sss:construction_HK}, which is probably well known to experts. 
\er

\subsection{The monodromy operator} \label{sec:revisit_HK}

In this subsection we explain the construction of the monodromy operator on the log crystalline cohomology for log schemes of ``generalized semistable type'' in the following sense: let $\ul k = (k, \N)$ be the standard log point, $\ul X/\ul k$ a log scheme such that $\etale$ locally, $\ul X$ is naively $\etale$ over a finite fiber product of log schemes that are given by log algebras of the form 
\[\begin{tikzcd} [column sep=1.2em]
\ul R^\square \arrow[r,equal] & \Big(\N^{r+1} \arrow[r, "\alpha"] & k[X_0, ..., X_r, ..., X_{d}]/\prod_{0 \le i \le r} X_i \Big) \\
\ul k \arrow[r, equal]  &  \Big(\: \N  \: \: \: \arrow[r] \arrow[u, swap, "1 \mapsto {(1, ..., 1)}"] & \: \: \:  k \arrow[u] \: \Big)
\end{tikzcd} \]
where the log structure $\alpha$ on $R^\square$ is given by $e_i \mapsto X_i$ for $ 0 \le i \le r$, here $e_i = (0, ..., 1 , ..., 0)$ with $1$ at the $i^{th}$ position. More precisely, $\ul X$ can be covered by charts which are naively $\etale$ over log algebras of the form $\ul S^\square = (M, S^\square)$, where
$$ S^\square = k[X_{l,h}]_{\substack{1\le l \le s \\ 0 \le h \le d_l}} / \Big( \prod_{1 \le i \le d_1} X_{1, i}, \cdots,  \prod_{1 \le i \le d_s} X_{s, i} \Big)$$
and $M = \underset{1 \le l \le s}{\bigsqcup_\N} \N^{r_l + 1}$ is the pushout of the diagonal maps $\N \ra \N^{r_l + 1}$ for $ 1 \le l \le s$. Clearly $\ul R^\square$ (resp. $\ul S^\square$) is fine and log-smooth over $\ul k$ of log-Cartier type. 

\bl \label{lemma:W_1_compare_works_with_trivial_log}
Let $\ul R^\square$ be as above, let $k^\circ = (k, 0)$ be the trivial log point, then $\ul R^\square/ k^{\circ}$ satisfies conditions (1) and (2) in Remark \ref{remark:variant_of_comparison_in_char_p}. In particular, there is a canonical isomorphism 
$$\omega_{\ul R^\square /k^{\circ}} \isom \mW_1\omega_{\ul R^\square /k^{\circ}}.  $$ The same conclusion holds for $\ul S^\square$. 
\el

\bproof 
For notational simplicity we check the claim for $\ul R^\square$. For (1) we may take the lifting $ A^\square = W(k) [X_0, ..., X_d]/\prod_{0 \le i \le r} X_i$ with log structure $\N^{r+1} \xrightarrow{e_i \mapsto X_i} A^\square$. For the lift of Frobenius we take the canonical lifting $\sigma$ on $W(k)$ and $X_i \mapsto X_i^p$, while Frobenius on $\N^{r+1}$ is given by multiplication by $p$. The log differential $\omega^1_{\ul A^\square/\ul k}$ is free over $\ul A$ with generators $\textup{dlog} X_0, ..., \textup{dlog} X_r, \textup{d} X_{r+1}, ..., \textup{d} X_{d}$, using notation in Subsection \ref{sss:log_differentials}, hence both condition (1) and (2) in Remark \ref{remark:variant_of_comparison_in_char_p} are satisfied. 
\eproof

The canonical map $\ul R^\square /k^{\circ} \longrightarrow \ul R^\square /\ul k$ induces a morphism of saturated log de Rham--Witt complexes 
$$ \Theta: \mW \omega^*_{\ul R^\square /k^{\circ}} \longrightarrow \mW \omega^*_{\ul R^\square /\ul k}.$$
The cochain complex $\ker \Theta$ is concentrated in degree $\ge 1$. 

\bp \label{prop:local_monodromy}
The shifted cochain complex $(\ker \Theta)[1]$ admits a structure of a Dieudonn\'e complex, which is saturated and strict. Moreover, there is a short exact sequence of Dieudonn\'e complexes
$$ 0 \ra  \mW  \omega^*_{\ul R^\square/\ul k} [-1] \longrightarrow \mW \omega^*_{\ul R^\square/k^{\circ}} \xrightarrow{\; \Theta \;} \mW \omega^*_{\ul R^\square /\ul k} \ra 0.$$ 
The same conclusion holds when $\ul R^\square $ is replaced by $\ul S^\square$. In particular, $(\ker \Theta)[1]$ carries a (log) Dieudonn\'e algebra structure. 
\ep

\bproof 
Again we check the proposition for $\ul R^\square$ since the only complication for $\ul S^\square$ is notational.  Write $K^* = \ker \Theta$ which is a $p$-torsion free sub-Dieudonn\'e complex of $\mW \omega^*_{\ul R^\square/k^{\circ}}$ (here we ignore the algebra structure, as in Subsection \ref{ss:L_eta_fixed_point}). It is straightforward to check that  $$\phi_F: K^* \longrightarrow \eta_p (K^*), \quad x \in K^i \mapsto p^i F (x) $$ is an isomorphism,  therefore $K^*[1]$ is a saturated Dieudonn\'e complex as multiplication by $p$ gives an isomorphism $\eta_p (K^* [1]) \isom \eta_p (K^*) [1]$. By Lemma \ref{lemma:kernel_of_strict_is_again_strict}, we have a short exact sequence 
$$ 0 \ra W_m (K^*) \longrightarrow \mW_m \omega^*_{\ul R^\square/k^{\circ}}  \longrightarrow  \mW_m  \omega^*_{\ul R^\square /\ul k}  \ra 0$$
for each $m$. Taking the inverse limit we see that $K^*[1]$ is strict. 

Next we construct a map of cochain complexes $\omega^*_{\ul A^\square/\ul k} \ra \omega^*_{\ul A^\square/ k^\circ}[1]$ preserving Frobenius on both sides. For this we define an element $\theta \in  \omega^1_{\ul A^\square/ k^\circ}$ from the log differential $\textup{dlog}: \N^{r+1} \ra  \omega^1_{\ul A^\square/ k^\circ}$ by 
$$\theta:= \textup{dlog} (e_0) + \cdot +  \textup{dlog} (e_r).$$
(More generally, for the product $\ul S^\square = (S^\square, M)$, we have a map $\N \ra M$ coming from each diagonal map $\N \ra \N^{r_l + 1}$. Let $e \in M$ be the image of $1 \in \N$. Then $\theta$ is the element $\textup{dlog} (e)$). 
We abuse notation and use the same symbol to denote the image of $\theta$ in $ \mW \omega^*_{\ul R^\square/k^\circ}$, then clearly $\theta \in K^1$. Since $\omega^*_{\ul A^\square/\ul k}$ is generated over $A^\square$ by $(\wedge_{i}\textup{dlog} X_i )\wedge (\wedge_{j}\textup{d} X_j)$, we may a map of complexes $\Psi: \omega^*_{\ul A^\square/\ul k} \ra \omega^*_{\ul A^\square/k^\circ}[1]$ by 
$$(\wedge_{i}\textup{dlog} X_i )\wedge (\wedge_{j}\textup{d} X_j) \mapsto (\wedge_{i}\textup{dlog} X_i )\wedge (\wedge_{j}\textup{d} X_j)\wedge \theta.$$
The map is well defined, since 
$$\sum_{0 \le i \le r} \textup{dlog} X_i \mapsto (\sum_{0 \le i \le r} \textup{dlog} X_i) \wedge \theta = \theta \wedge \theta = 0.$$
Moreover, since $d (\theta) = 0 $ (by $d \circ \textup{dlog} = 0$) and $F(\theta) = \theta$, the map $\Psi$ preserves differential and Frobenius structures on both sides. 

By passing to the (p-adic) completion and strictification (complete saturation), we get a map $\Psi: \mW \omega^*_{\ul R^\square/\ul k} \ra \mW \omega^*_{\ul R^\square/k^{\circ}} [1]$. By the choice of $\theta$, the image of $\Psi$ lie in $K^*[1]$, so we have now constructed a map of strict  Dieudonn\'e complexes
$$\Psi: \mW \omega^*_{\ul R^\square/\ul k} \longrightarrow K^*[1].$$
To finish the proof of the proposition, it remains to prove that this map is an isomorphism. By Corollary 2.7.4 in \cite{BLM} \footnote{see also Corollary \ref{cor:Cartier_implies_qi}. Note that the proof we cite from \cite{BLM} is stated for Dieudonn\'e complexes}, it suffices to show that 
$$W_1(\Psi): \mW_1 \omega^*_{\ul R^\square/\ul k} \longrightarrow W_1(K^*[1]) $$
is an isomorphism. By Lemma \ref{lemma:kernel_of_strict_is_again_strict} and Lemma \ref{lemma:W_1_compare_works_with_trivial_log}, we need to check that 
$$\omega^*_{\ul R^\square/\ul k} \xrightarrow{\wedge \theta_1} \ker (\omega^*_{\ul R^\square/ k^\circ} \ra \omega^*_{\ul R^\square/\ul k})[1]$$ is an isomorphism of cochain complexes, which follows from their explicit descriptions.  
\eproof  

The following lemma is used in the proof above

\bl  \label{lemma:kernel_of_strict_is_again_strict}
Let $\Theta: A^* \longrightarrow B^*$ be a map between two saturated Dieudonn\'e complexes. Let $K^* = \ker \Theta$ be the sub Dieudonn\'e complex (which is saturated). Define the standard filtration on $A^*$ (resp. on $B^*$ and $K^*$) as in Subsection \ref{ss:strict_completion} by $\textup{Fil}^i(A^*) = V^i(A^*) + \textup{d} V^i (A^*)$. Then 
$$ \textup{Fil}^i(K^*) = \ker (\textup{Fil}^i(A^*) \xrightarrow{\; \Phi \;} \textup{Fil}^i(B^*)).$$
\el

\bproof  
We need to show that $ \ker (\textup{Fil}^i(A^*) \xrightarrow{\; \Phi \;} \textup{Fil}^i(B^*)) \subset \textup{Fil}^i(K^*)$. In other words, $x = V^i y + d V^i z$, and $\Phi(x) = 0$, we claim that $x \in \textup{Fil}^i(K^*)$. For notational simplicity we write $\cl x = \Phi (x)$ for the image of $x \in A^*$ in $B^*$.  By assumption we know that 
$$ d \cl z = F^i d V^i  \cl z= - F^i(V^i \cl y)  = - p^i \cl y \in B^*.$$
Since $B^*$ is saturated, we know that $\cl z = F^i \cl w$ for some $\cl w \in B^*$, so there exists $w \in A^*$ and $\alpha \in K^*$ such that 
$$x = V^i y + d V^i (F^i w + \alpha) = V^i  (y + F^i d w) + d V^i \alpha.$$
Therefore, it suffices to prove the claim for $x$ of the form $x = V^i y$, but this is clear since $V^i \cl y = 0$ implies that $\cl y = 0$ as $B^*$ is $p$-torsion free. 
\eproof 

As a corollary, we arrive at Theorem \ref{thm:main_monodromy} in the introduction. 

\bc  \label{cor:monodromy}
Let $\ul X/\ul k$ be a log scheme of ``generalized semistable type'', then there is a short exact sequence of cochain complexes 
$$ 0 \ra \mW \omega^*_{\ul X/\ul k} [-1] \longrightarrow \mW \omega^*_{\ul X/k^\circ} \longrightarrow \mW \omega^*_{\ul X/\ul k} \ra 0$$
that preserves the Frobenius structures. As a consequence, this gives rise to 
a connecting homomorphism on the log crystalline cohomology  $$N: H_{\log\textup{-cris}}^*(X/W(\ul k)) \longrightarrow H_{\log\textup{-cris}}^*(X/W(\ul k)).$$ This operator satisfies $N \varphi = p \varphi N$ and agrees with the monodromy operator constructed in \cite{HK}. 
\ec

\bproof 

Let $R$ be a $k$-algebra with an $\etale$ coordinate, namely an $\etale$ morphism $S^\square \ra R$ where $S^\square$ has the form described as in the beginning of this section.  Equip $R$ with the log structure from $M \ra S^\square$. Then for each $m$ we have 
$$ 0 \ra \mW_m \omega^*_{\ul R/\ul k} [-1] \longrightarrow \mW_m \omega^*_{\ul R/k^{\circ}} \longrightarrow \mW_m \omega^*_{\ul R/\ul k} \ra 0$$
from the $\etale$ base change of the corresponding sequence for $\mW_m \omega^*_{\ul S^\square/ \ul k}$. 
To globalize from these $\etale$ coordinates, we observe that the element $\theta \in \mW \omega^1_{\ul R/k^\circ}$ is independent of the choice of coordinates, since it can be described as the image of $\textup{dlog} (e)$ in $\omega^1_{W(\ul R)/W(k^\circ)} \longrightarrow W(\omega^*_{W(\ul R)/W(k^\circ)})^1_{\mathfrak{sat}} = : \mW \omega^1_{\ul R/k^\circ}$, so the sequences above glue to an exact sequence of sheaves 
$$ 0 \ra \mW_m  \omega^*_{\ul X/\ul k} [-1] \longrightarrow \mW_m \omega^*_{\ul X/k^\circ} \longrightarrow \mW_m \omega^*_{\ul X/\ul k} \ra 0.$$ Since there is no higher $\textup{Rlim}$ terms, taking the inverse limit we reach the conclusion. The relation $N \varphi = p \varphi N$ follows from $d F = p Fd $ since $\varphi$ is induced by $\phi_F: x \in \mW\omega^i_{\ul X /\ul k} \mapsto p^i F(x) \in \mW\omega^i_{\ul X /\ul k}$.  

Finally, the comparison with Hyodo--Kato's construction follows from Theorem \ref{thm:compare_with_HK} and the description of $(W_n \sq \omega^q_{Y})'$ in 4.20 in \cite{HK}. 
\eproof

\newpage 

\section{$\Ainf$-cohomology theory} \label{sec:Ainf_intro}
In the rest of the article we formulate and prove Theorem \ref{mthm:log_crys_compare_global} in the introduction. The goal of this section is to recall the theory of $\Ainf$-cohomology developed by Bhatt--Morrow--Scholze in \cite{BMS}. In particular, we need the Hodge--Tate comparison (Theorem \ref{thm:dR_comparison}) and the fact that $L\eta$ commutes with a certain form of global sections on small enough objects (Theorem \ref{prop:kestutis_3_20}), as generalized to the semistable case by \cite{Kestutis}. We also carry out some local computations needed for the log crystalline comparison. 

\begin{notation}
$\bullet $ For the rest of the article, let $k$ be an algebraically closed field over $\F_p$ and $C$ be a completed algebraic closure of $W(k)[\frac{1}{p}]$.  Let $\fm$ be the maximal ideal of its ring of integers $\mO_C$. For convenience we also fix a choice of $p^\Q \subset C$ (in particular a choice of compatible $p^{n}$th roots of unity $\{\zeta_{p^n}\} \in \mu_{p^n} (C)$). Consider the element $\epsilon = (1, \zeta_p, \zeta_{p^2}, ...)  \in \mO_C^{\flat}$, which determines 
$$\mu := [\epsilon] - 1, \quad  \xi := \frac{\mu}{\varphi^{-1} (\mu)} = \sum_{i = 0}^{p-1} [\epsilon]^{\frac{i}{p}}, \quad \sq \xi := \varphi(\xi)$$
in  $\Ainf = W(\mO_C^{\flat})$ where $\mO_C^{\flat}$ is the tilt of $\mO_C$. 
The ring $\Ainf$ is equipped with  the $(p, \mu)$-adic topology.  

$\bullet$ Let $\theta$ be the surjection $\Ainf \twoheadrightarrow \mO_C$ lifting the canonical morphism $\mO_C^{\flat} \ra \mO_{C}/p$, and define $\sq \theta: = \theta \circ \varphi^{-1}$. Therefore we have $ \ker (\theta) = (\xi)$ (resp. $\ker (\sq \theta ) = (\sq \xi)$).
We will call the map $\theta$ (resp. $\sq \theta$) the de Rham  (resp. Hodge-Tate) specialization map. 

$\bullet $ Let  $\sq \vartheta: \Ainf \ra W(k)$ be the unique lift of the quotient map $\Ainf \xrightarrow{\sq \theta} \mO_C \ra k$ and define $ \vartheta := \sq \vartheta \circ \varphi = F \circ \sq \vartheta$  where $F$ denotes the Witt vector Frobenius. The map $\sq \vartheta$ will be referred to as the crystalline specialization map  \footnote{ It is helpful to observe that under the crystalline specializations and its twist, we have 
$\vartheta ([\epsilon^{\frac{1}{p}}]) = \sq \vartheta([\epsilon^{\frac{1}{p}}])  = \vartheta ([\epsilon]) = \sq \vartheta([\epsilon]) = 1$ and $ \vartheta (\xi)  = \vartheta (\sq \xi) =\sq  \vartheta (\xi)  =  \sq \vartheta(\sq \xi) = p.$ }. Unless otherwise stated, all maps $\Ainf \ra W(k)$ in the rest of the article are assumed to be the crystalline specialization $\sq \vartheta$ by default. 

$\bullet$ Finally, we make the following notation changes from the previous sections. 
\bi
\item[$-$] We will let $k$ be an algebraically closed field over $\F_p$.
\item[$-$] In the rest of the article, we will use $R$ to denote $\mO_C$-algebras (instead of a $k$-algebra), its base change to $k$ is denoted by $R_k$. 
\ei

\end{notation}

\subsection{Semistabe coordinates and log structures} \label{ss:semistable_setup}

Let $\fX$ be a $p$-adic formal scheme over $\mO_C$ of generalized semistable reduction, in the sense that there exists a covering by affine opens $\fU = \spf R$ in the $\etale$ topology of $\fX$, such that each $\fU$ admits an $\etale$ morphism (called a \textit{semistable coordinate})
 $$ \square: \fU  \longrightarrow \spf R^\square = \spf (R_1^{\square} \widehat \otimes \cdots \widehat \otimes R_s^{\square})$$ over $\mO_C$,
where each $R_i^{\square}$ is of the form
$$\mO_C \gr{T_0, ..., T_r, T_{r+1}^{\pm 1}, ..., T_{d}^{\pm 1}}/ (T_0 \cdots T_r - \varpi)$$ 
with $0 \le r \le d$, where $\varpi \in \fm \minus \{0\}$ a non-unit of rational valuation.   Let $\fX_C := \fX_{\eta}^{\textup{ad}}$ denote the adic generic fiber of $\fX$, which is smooth over $C$ by the requirement above.  We equip $\mO_C$ (resp. $\mO_C/p$) with the log structure $\mO_C \minus \{0\} \hookrightarrow \mO_C$ (resp. its pullback log structure), and $\fX$ its ``divisorial log structure'' $M_{\fX}$ given by the sheafification of $(\mO_{\fX} [\frac{1}{p}])^\times \cap \mO_{\fX}$.  The log scheme $\ul{\fX}_{\mO_C/p}$ over $\spec \ul{\mO_C}/{p}$ and its special fiber $\ul{\fX}_{k}$ are equipped with the pullback log structures. Similarly define $\ul{\fU}$ and $\ul{\fU_k}$.  The semistable coordinate morphism $\fU = \spf R \ra \spf (R_1^{\square} \widehat \otimes \cdots \widehat \otimes R_s^{\square})$ is the $p$-adic completion of (the $\mO_C$-base change of) an $\etale$ $\mO$-morphism 
\begin{align*}
U \quad  \longrightarrow \quad  &  \spec  \Big( \mO [T_{h, i}, T_{h, j}^{\pm 1}]_{1 \le h \le s}/ \big( \prod T_{h, i} - \varpi_h \big)_{1 \le h \le s} \Big) \\
= & \spec \Big( \mO [T_{1, i}, T_{1, j}^{\pm 1}]/ (T_{1, 0}\!\cdot\!\cdot\!\cdot\! T_{1, r_1} \text{--} \varpi_1) \Big) \times \cdots \\   
& \qquad \qquad \qquad  \times \spec \Big( \mO [T_{s, i}, T_{s, j}^{\pm 1}]/ (T_{s, 0}\!\cdot\!\cdot\!\cdot\! T_{s, r_s} \text{--} \varpi_s) \Big)
\end{align*}
where $\mO \subset \mO_C$ is discretely valued, and $\varpi_1, ..., \varpi_s \in \mO$. We still call this a semistable coordinate and continue to denote it by $\square$. This determines a choice of charts for the (similarly defined) divisorial log structure on $U$. More precisely, let $\sq M^\square$ be the monoid defined by the following push-out diagram in monoids:  
\[
\begin{tikzcd}
\N^s \arrow[r, "\textup{diag}"] \arrow[d] 
& \N^{r_1 + 1} \oplus \cdots \oplus \N^{r_s+1} \arrow[d]\\
\mO \minus \{0\} \arrow[r] 
& \sq M^\square= \Big( \underset{1 \le h \le s}{\oplus} \N^{r_h+1} \Big) \sqcup_{\N^s} (\mO\minus \{0\})
\end{tikzcd}
\]
where the top horizontal map is the ``diagonal'' map, given by $1 \mapsto (1, ..., 1)$ on each $\N \ra \N^{r_h+1}$;  and the left vertical map $\N^s \ra \mO\minus \{0\}$ is 
$(m_1, ..., m_s) \longmapsto \varpi_1^{m_1} \cdots \varpi_s^{m_s}.$
Then the log structure $(\mO_{U}[\frac{1}{p}])^\times \cap \mO_{U}$ has a chart given by 
$\sq M^\square \longrightarrow \Gamma (U, \mO_U),$
where $(\{m_{h, i}\}_{\substack{1 \le h \le s \\0 \le i \le r_h}}, a)$ gets mapped to $\big(\prod (T_{h, i})^{m_{h, i}} \big) \cdot a \in \Gamma (U, \mO_U)$. \footnote{See Lemma \ref{lemma:etale_chart}.} Moreover, the pre-log structure described by the chart above is the base change of 
\[
\begin{tikzcd}
M^\square :=  \underset{1 \le h \le s}{\oplus} \N^{r_h+1} \arrow[r] &R_{\circ} : = \Gamma (U, \mO_U) \\
N^\square := \N^s \arrow[r] \arrow[u, "\textup{diag}"] & \mO \arrow[u]
\end{tikzcd}
\]
along $(\mO, N^\square = \N^s) \ra (\mO, \mO\minus \{0\})$, where the underlying map of rings is identity and the monoid map is $(m_1, ..., m_s) \longmapsto \varpi_1^{m_1} \cdots \varpi_s^{m_s}$. 

\br \label{remark:pullback_log_structure_1} From the discussion above, even though $\ul \fU = \spf \ul R$ is strictly speaking not log-smooth over $\ul{\mO_C}$ (the log structures are not fine), we may treat it as if it were so, since the map $\ul{\mO_C} \longrightarrow \ul R$ is the base change of the log-smooth morphism $ (\mO_C, N^\square) \longrightarrow (R, M^\square) $ of fine log rings along  $(\mO_C, N^\square) \longrightarrow (\mO_C, \mO\minus \{0\})$, which is defined as above with $\mO_C$ in place of $\mO$. Consequently, the continuous log differentials $ \omega^{*, \textup{ct}}_{\ul R/\ul{\mO_C}} \cong  \omega^{*, \textup{ct}}_{(R, M^\square)/(\mO_C, N^\square)}$ are isomorphic, and the $\mO_{\fX}$-module  $\omega^{*, \textup{ct}}_{\ul \fX/ \ul{\mO_C}}$ is locally free  of finite rank.  Now we turn to the special fiber $R_k = R \otimes_{\mO_C} k$ of $\fU = \spf R$. The $\etale$ local chart $(R_k, \sq M^\square)$ of $\ul \fX_k$ over $(k, \mO_C \minus \{0\})$ is the base change (this time along $(k, N^\square) \ra (k, \mO_C \minus \{0\})$) of $ \ul k^\square := (k, N^\square) \longrightarrow (R_k, M^\square),$
which is log-smooth of log-Cartier type. This in turn implies that 
$$\omega^*_{(R_k, M^\square)/\ul k^\square} \cong \omega^*_{(R_k, \sq M^\square)/(k, \mO_{C} \minus \{0\})} \cong \omega^{*, \textup{ct}}_{\ul{R}/\ul{\mO_C}} \otimes_{\mO_C} k $$
is locally free over $R_k$ and satisfies the Cartier isomorphism. More importantly, by Remark \ref{remark:variant_of_comparison_in_char_p}, the comparison with de Rham complexes (in particular the conclusion of Proposition  \ref{prop:comparison_with_mod_p_dR}) still applies to $\mW \omega^*_{(R_k, M^\square)/\ul k^\square}$ and the $\etale$ sheaf $\mW \omega^*_{\ul \fX_k/ \ul k}$, where $\ul k = (k, \mO \minus \{0\})$. In particular, the map $(R_k, M^\square) \ra \ul{R_k} = (R_k, \Gamma_{\ett}(\fU, M_{\fX}))$ induces an isomorphism on the underlying Dieudonn\'e algebras
$$ \mW \omega^*_{(R_k, M^\square)/\ul k^\square} \isom \mW \omega^*_{\ul{R_k}/\ul k}$$  of the corresponding log de Rham--Witt complexes, by Corollary \ref{cor:Cartier_implies_qi}. 
\er 

\br  \label{remark:pullback_log_structure_2}
In our setup, we still have an isomorphism $\mW \omega^*_{\ul{\fX_k}/\ul k} \isom W^{\textup{\tiny HK}} \omega^*_{\ul{\fX_k}/\ul k}$\footnote{Even though in the original work of \cite{HK}, all log schemes are assumed to be fine, we may still define the Hyodo--Kato complex by taking the inverse limit of $R^* u^{\log}_{\ul X/W_n, *} (\mO_{\ul X/\ul{W_n}})$ as in Subsection \ref{sss:construction_HK}. }, and both compute the log crystalline cohomology of $\ul{\fX_k}$ over $\ul k$. The first claim essentially follows from the same line of arguments as the proof of Theorem \ref{thm:compare_with_HK}. To justify this, we regard both complexes as sheaves on the $\etale$ site $(\fX_{k})'_{\ett, \textup{aff}}$ whose objects consisting of $\spec R_k$, where $\spf R$ is a formal affine scheme admitting semistable coordinates as above. On each $\spec R_k$, $W^{\textup{\tiny HK}} \omega^*_{\ul{R_k}/\ul k}$ is still a $p$-torsion free Dieudonn\'e complex (by the proof of Corollary 4.5 of \cite{HK}) and in fact saturated by the proof of Lemma \ref{lemma:HK_complex_is_ADc_2}, the key being the Cartier isomorphism (which still holds by the previous remark). Therefore, $W^{\textup{\tiny HK}} \omega^*_{\ul{\fX_k}/\ul k}$ takes values in ($p$-compatible) strict log de Rham--Witt complexes over $W(\ul k)$, and we have a map $\gamma: \mW \omega^*_{\ul{\fX_k}/\ul k} \ra W^{\textup{\tiny HK}} \omega^*_{\ul{\fX_k}/\ul k}$ coming from each $\ul{R_k} \ra W_1^{\textup{\tiny HK}} \omega^0_{\ul{R_k}/\ul k}$. Now apply the Cartier isomorphism one more time we get an isomorphism $W_1^{\textup{\tiny HK}} \omega^*_{\ul{R_k}/\ul k} \isom \omega^*_{\ul{R_k}/\ul k}$, and again by Remark  \ref{remark:pullback_log_structure_1}, this implies that $\gamma$ is an isomorphism. For the second claim in this remark, note that there is a canonical quasi-isomorphism $W_n^{\textup{\tiny HK}} \omega^*_{\ul{\fX_k}/\ul k} \isom R u^{\log}_{\ul X/W_n, *} (\mO_{\ul X/\ul{W_n}})$ for each $n \ge 1$, following the notation in the beginning of Subsection \ref{sss:construction_HK}. The proof in \cite{HK} (Theorem 4.19) goes through and amounts to the Cartier isomorphism once again.  
\er

\subsection{$\Ainf$-cohomology} \label{ss:local_analysis_sq_A_Omega}

Following \cite{BMS}, we define $ A\Omega_{\fX} := L \eta_{\mu} R \nu_{*} (\Ainfx)$, where  $\nu: \fX_{C, \pet} \ra \fX_{\ett}$ is the $p$-adic nearby cycle map and $\Ainfx:= W(\widehat \mO^{+, \flat}_{\fX_C})$ is the period sheaf on the pro-$\etale$ site of $\fX_C$ (cf. \cite{Scholze}). The Frobenius automorphism on $\widehat \mO^{+, \flat}_{\fX_C}$ lifts to a Frobenius on $\Ainfx$, which by functoriality induces an $\Ainf$-linear Frobenius map $ \varphi_{\fX}: A \Omega_{\fX} \ra \varphi_* A \Omega_{\fX}$. The map $\varphi_{\fX}$ is the composition \footnote{The isomorphism $\sq \varphi_{\fX}$ is called the divided Frobenius in \cite{BLM}.}
$$  A \Omega_{\fX} \xrightarrow{\sq \varphi_{\fX}} \varphi_* L \eta_{\varphi(\mu)} R \nu_* \Ainfx  =  \varphi_* L \eta_{\sq \xi} A \Omega_{\fX} \ra \varphi_* A \Omega_{\fX} .$$
Similar to Subsection \ref{sec:etale_base_change}, it is convenient to consider the smaller site $\fX_{\ett, \textup{aff}}$, consisting of open affine formal schemes $\fU \in \fX_{\ett}$ which admit $\etale$ semistable coordinates, where the topology is given by $\etale$ coverings. As in \cite{BMS}, there are two presheaves on $\fX_{\ett, \textup{aff}}$ that play important roles for us, which we define now. Let $\mD(\fX_{\ett, \textup{aff}}, \Ainf)$ be the derived $\infty$-category of sheaves of $\Ainf$-modules on $\fX_{\ett, \textup{aff}}$, and $\textup{Psh}(\fX_{\ett, \textup{aff}}, \mD(\Ainf))$ the $\infty$-category of presheaves valued in $\mD(\Ainf).$ \footnote{The homotopy category of $\mD(\fX_{\ett, \textup{aff}}, \Ainf)$ is identified with the usual triangulated derived category $\textup{D}(\fX_{\ett, \textup{aff}}, \Ainf)$. Similarly, the homotopy category of $\textup{Psh}(\fX_{\ett, \textup{aff}}, \mD(\Ainf))$ is $\textup{D}(\fX_{\ett, \textup{aff}}^{\textup{pre}}, \mD(\Ainf)) $ where $\fX_{\ett, \textup{aff}}^{\textup{pre}}$ is the site consisting of the same objects as $\fX_{\ett, \textup{aff}}$ but equipped with the indiscrete (= chaotic) Grothendieck topology.}
Let  $$\iota: \mD(\fX_{\ett, \textup{aff}}, \Ainf) \longrightarrow \textup{Psh}(\fX_{\ett, \textup{aff}}, \mD(\Ainf))$$
denote canonical (fully faithful) map which takes $\mF$ to $\iota(\mF): \fU \mapsto R \Gamma (\fU, \mF)$. The functor $\iota$ admits a left adjoint $j^{-1}$ which is the sheafification functor. We will interchangeably denote $\iota(\mF)$ by $\mF^{\textup{pre}}$ to emphasize that it is the presheaf associated to $\mF$ (viewed as a sheaf taking values in the derived $\infty$-category $\mD(\Ainf)$). Now we define the following two objects in $ \textup{Psh}(\fX_{\ett, \textup{aff}}, \mD(\Ainf))$: \footnote{The presheaf $A \Omega^{\textup{pre}}_{\fX}$ is denoted by $A \Omega^{\textup{sm}}_{\fX}$ in \cite{BLM} in the case of good reductions; the presheaf $A \Omega^{\textup{psh}}_{\fX}$ agrees with the notation used in \cite{Kestutis}, where it is considered as an object in the homotopy category.}
\bi
\item  $A \Omega^{\textup{pre}}_{\fX} : = \iota (A \Omega_{\fX}) = (L \eta_\mu R \nu_* \AAinf)^{\textup{pre}}$; 
\item $A \Omega^{\textup{psh}}_{\fX} := L \eta_\mu (\iota (R \nu_* \AAinf)) = L\eta_\mu (R \nu_* \AAinf)^{\textup{pre}}$. 
\ei
Unwinding definitions, we know that on each affine open $\fU = \spf R \in \fX_{\ett, \textup{aff}}$, the global sections of the  two presheaves are 
  $$A \Omega^{\textup{pre}}_{\fX} (\fU) = R \Gamma_{\ett}(\spf R, A \Omega_{\fX}); \qquad A \Omega^{\textup{psh}}_{\fX} (\fU) = L \eta_\mu R\Gamma_{\pet} (U, \AAinf)$$
where $U = \fU_{\eta}^{\textup{ad}}$ is the adic generic fiber of $\fU$. Both presheaves sheafify to $A\Omega_{\fX}$, so we have a natural map $A \Omega^{\textup{psh}}_{\fX} \longrightarrow A \Omega^{\textup{pre}}_{\fX}$ by adjunction. To simplify notation we will often write  
 $A\Omega_R: =A \Omega_{\fX}^{\textup{pre}} (\fU)$  and $A \Omega_R^{\textup{psh}} := A \Omega_{\fX}^{\textup{psh}} (\fU)$.

\subsection{Local perfectoid covers and group cohomology} \label{sss:local_setup_of_R_infty}
 In this subsection, we recall the local theory of $\Ainf$-cohomology for small enough $\spf R$ in $\fX_{\ett}$, mostly developed in \cite{Kestutis}, following the general strategy outlined in \cite{BMS} -- namely to compute both $A\Omega_R^{\textup{psh}}$ and $A \Omega_R$ (thus to relate them) using continuous group cohomology.  
Let us assume that $\spf R \in \fX_{\ett}$ admits an $\etale$ coordinate map  $\square: \spf R \ra \spf (R_1^\square \widehat \otimes \cdots \widehat \otimes R_s^\square)$, with each 
$$R_h^\square = \mO_C \gr{T_{h, i}, T_{h, j}^{\pm1}} / (\prod T_{h, i} - \varpi_h) $$ 
where the variables range over $0 \le i \le r_h, \; r_h+1 \le j \le  d_h$ and $\varpi_h = p^{q_h}$ for some $q_h \in \Q_{> 0}$ as in Subsection \ref{ss:semistable_setup}. In particular, our choice of $p^\Q \in \mO_C$ determines a choice of $\varpi_h^{1/p^m} \in \mO_C$  for each $m \ge 1$. 

\subsubsection{Construction of $R_\infty$}
 For each $m \ge 1$, we define 
$$R_{h,m}^\square = \mO_C \gr{T_{h, i}^{1/p^m}, T_{h, j}^{\pm 1/p^m}} / (\prod T_{h, i}^{1/p^m} - \varpi_h^{1/p^m}), $$
and then let $R_{m}^\square := R_{1, m}^\square \widehat \otimes \cdots \widehat \otimes R_{s, m}^\square$. In the direct limit we get ``perfectoid covers''
$$ R_{h, \infty}^\square :=  (\varinjlim R_{h, m}^\square)^{\widehat{\;}}, \quad \:\: R_\infty^\square : = R_{1, \infty}^\square \widehat \otimes \cdots \widehat \otimes R_{s, \infty}^\square.$$
These perfectoid $\mO_C$-algebras admit explicit decompositions as $\mO_C$-modules by
\begin{align*}
R_\infty^\square & =  \widehat \bigoplus_{a_{h, i}} \mO_C  T_{1, 0}^{a_{1, o}} \cdots T_{1, d_1}^{a_{1, d_1}}  \cdots   T_{s, 0}^{a_{s, 0}}  \cdots  T_{s, d_s}^{a_{s, d_s}}  \\
& = R^\square \oplus M^\square_\infty,
\end{align*}
where the completed direct sum ranges over  $a_{h, 0}, ..., a_{h, d_h} \in \Z[\frac{1}{p}]$ for $1 \le h \le s$, such that each $a_{h, i} \ge 0$ and $\prod_i a_{h, i} = 0$ for each $h$.  The $R^\square$-module decomposition $R_\infty^\square = R^\square \oplus M^\square_\infty$ follows the same notation used in \cite{Kestutis}, where the direct summand $M^\square_{\infty}$ denotes the completed direct sum of all $ \mO_C \cdot T_{1, 0}^{a_{1, o}}  \cdot \cdots \cdot  T_{s, d_s}^{a_{s, d_s}}$ such that there exists some non-integral supscript $a_{h, i}$ or $a_{h, j}$ (that belongs to $\Z[\frac{1}{p}] \minus \Z$). 
This allows us to obtain perfectoid $R$-algebras by setting $R_m := R_m^\square \otimes_{R^\square} R_m^\square, $ and then 
$$ R_\infty := (R \otimes_{R^\square} R_\infty^\square)^{\widehat{\;}} = R \oplus M_\infty.$$ 

We write $\Z_p (1) = (\varprojlim_{m \ge 0} \mu_{p^m} (\mO_C))$, 
and define for each $1 \le h \le s$ the following abelian group
$$\Delta_h := \Big\{ (\epsilon_{h, 0}, ..., \epsilon_{h, d_h}) \in \Z_p(1)^{\oplus (d_h+1)} \: \vline \: \epsilon_{h, 0} + \cdots + \epsilon_{h, d_h} = 0 \Big\}. $$ 
Finally we set $\Delta : = \Delta_1 \oplus \cdots \oplus \Delta_s,$ with each $\Delta_h$ isomorphic to $\Z_p^{\oplus d_h}$,  generated by 
\bi
\item $\gamma_{h, i} = (-1, 0, ..., 0, 1, 0, ..., 0)$ \; for $i = 1, ..., r_h$,\; \quad with $1$ at the $i^{th}$ entry;
\item $\gamma_{h, j} = (0, ..., 0, 1, 0, ..., 0)$\;  for $j = r_h +1, ..., d_h$, \;\: with $1$ at the $j^{th}$ entry.
\ei
Under this identification, the element $\gamma_{h, i}$ acts on $R_{h, m}^\square$ by sending 
$$T_{h, 0}^{1/p^m} \mapsto \zeta_{p^m}^{-1} T_{h, 0}^{1/p^m},\quad T_{h, i}^{1/p^m} \mapsto \zeta_{p^m} T_{h, i}^{1/p^m}, \quad T_{h, k}^{1/p^m} \mapsto T_{h, k}^{1/p^m} \textup{ for } k \ne 0, i$$ 
\big(resp. $\gamma_{h, j}$ acts by 
$T_{h, j}^{1/p^m} \mapsto \zeta_{p^m}  T_{h, j}^{1/p^m}, \: T_{h, k}^{1/p^m} \mapsto T_{h, k}^{1/p^m} \textup{ for } k \ne j \big)$. The product $\Delta$ therefore acts on $R^\square_m$, and continuously on $R_\infty^\square$ (resp. $R_\infty$).

\subsubsection{Construction of $A(R^\square)$}
We next observe that
$$(R^\square_{h, \infty})^\flat = \Big(\varinjlim_{m} \mO_C^\flat [U_{h, i}^{1/p^m}, U_{h, j}^{\pm 1/p^m}] / (\prod U_{h, i}^{1/p^m} - (\varpi_h^{\flat})^{1/p^m}) \Big)^{ \widehat{\;}} $$
where $U_{h, i}^{1/p^m}$ denotes the elements 
$$U_{h, i}^{1/p^m} = (T_{h, i}^{1/p^m}, T_{h, i}^{1/p^{m+1}},  T_{h, i}^{1/p^{m+2}}, ...) \in (R^\square_{h, \infty})^\flat = \varprojlim_{y \mapsto y^p} R^\square_{h, \infty} $$ for each $i = 0, ..., r_h$. We likewise define elements $U_{h, j}$ for $j = r_{h}+1, ..., d_h$.\footnote{We reserve the symbols $X_{h, i}$ for the Teichumuller representatives of $U_{h, i}$ in the Witt vectors.}\footnote{The tilt $(R_\infty)^\flat$ of $R_\infty$ can be identified with the $p^{\flat}$-adic completion of $R'$, where $(R_\infty^\square)^\flat \ra R'$ is any lift of the $\etale$ algebra $R_\infty^\square /p \ra R_\infty / p$ to an $\etale$ $(R_\infty^\square)^\flat$-algebra.} 
From the description of $(R_\infty^\square)^\flat$ and the decomposition $R_\infty^\square = R^\square \oplus M_\infty^\square$, we have 
\begin{align*}
\Ainf (R_\infty^\square) &  \cong \widehat \bigotimes_{h} \Big(\varinjlim_{m} \Ainf   \Big[X_{h, i}^{1/p^m}, X_{h, j}^{\pm 1/p^m} \Big] / (\prod X_{h, i}^{1/p^m} - [\varpi_h^{\flat}]^{1/p^m})  \Big)^{\widehat{\;}} \\
&  \cong  \widehat {\bigoplus}_{\substack{a_{h, 0}, ..., a_{h, d_h} \in \Z[\frac{1}{p}]_{_{\ge 0}} \\  \prod_i a_{h, i} = 0, \: 1 \le h \le s}} \:\:\: \Ainf \cdot X_{1, 0}^{a_{1, o}} \cdots X_{1, d_1}^{a_{1, d_1}} \cdots X_{s, d_s}^{a_{s, d_s}} \\
& \cong A(R^\square) \oplus N^\square_\infty \end{align*}
with $X_{h, i}^{1/p^m} = [U_{h, i}^{1/p^m}]$ (likewise for $X_{h, j}$), and all completions being $(p, \mu)$-adic. We pose to remark that in particular, we have $\sq \theta(X_{h, i}^{1/p^m}) = T_{h, i}^{1/p^{m+1}}$ under the map $\sq \theta: \Ainf(R_\infty^\square) \ra R_\infty^\square$. Here the subring $A(R^\square)$ is given by 
$$ \Ainf (R^\square)   \cong \widehat \bigotimes_{h} \Big( \Ainf   \Big[X_{h, i}, X_{h, j} \Big] / (\prod X_{h, i}  - [\varpi_h^{\flat}])  \Big)_{(p, \mu)}^{\widehat{\;}}  $$
We have similar decompositions for $\Ainf(R_\infty) = W(R_\infty^\flat) = A(R) \oplus N_\infty$. Moreover, $A(R)$ is (by construction) formally $\etale$ over $A(R^\square)$, and also $(p, \mu)$-adically complete. 
The group $\Delta = \oplus \Delta_h$ acts on $\Ainf (R_\infty^\square)$ (resp. $\Ainf (R_\infty)$) by functoriality, and the action respects the decomposition $A(R^\square) \oplus N_\infty^\square$ (resp. $A(R) \oplus N_\infty$) described above. Moreover, by the same analysis in \cite{BMS} (or \cite{Kestutis}), for each $i$ the (continuous) group cohomology $H^i(\Delta, N_\infty^\square)$ (resp. $H^i(\Delta, N_\infty)$)  is entirely $\mu$-torsion, hence $L \eta_\mu R\Gamma_{\textup{ct}} (\Delta, N_\infty) = 0$ and consequently $L\eta_\mu R\Gamma_{\textup{ct}} (\Delta, \Ainf(R_\infty))$ only involves the action of $\Delta$ on $A(R)$. 

\subsubsection{The computation of $\sq \Omega_R$ and $A \Omega_R$}

To ease notation we write $\widehat \mO^+$ for the pro-$\etale$ sheaf $\widehat \mO_{\fX_C}^+$. By the almost purity theorem (\cite{Scholze} Theorem 4.10), the perfectoid pro-$\etale$ $\Delta$-cover $\fU_{\eta, \infty}^{\textup{ad}} = \textup{Spa} (R_\infty[\frac{1}{p}], R_\infty)$ over $U = \fU_\eta^{\textup{ad}} := (\spf R)_{\eta}^{\textup{ad}}$ gives rise to an almost quasi-isomorphism $e:  R\Gamma_{\textup{ct}} (\Delta, R_\infty) \xrightarrow{ \: \sim \: }_a R\Gamma_{\pet} (U, \widehat \mO^+)$.\footnote{this, of course, depends on the choice of coordinates.} Let $\sq \Omega_{\fX} := L \eta_{\zeta_p - 1} R \nu_* \widehat \mO^+$ and define $\sq \Omega_{\fX}^{\textup{psh}}$ and $\sq \Omega_{\fX}^{\textup{pre}}$ as in Subsection \ref{ss:local_analysis_sq_A_Omega}.

\bt[\cite{Kestutis} 3.9, 3.20, 4.6]  \label{prop:kestutis_3_20}
We have quasi-isomorphisms 
$$L\eta_{\zeta_p - 1} (R\Gamma_{\textup{ct}} (\Delta, R_\infty)) \isom  L\eta_{\zeta_p - 1} (R\Gamma_{\pet} (U, \widehat \mO^+)) \isom R \Gamma(\fU, \sq \Omega_{\fX}) $$
where the the quasi-isomorphism is given by $L\eta_\mu (e)$  and the second comes from the canonical map $\sq \Omega_{\fX}^{\textup{psh}} \ra \sq \Omega_{\fX}^{\textup{pre}}$. Analogously, there are quasi-isomorphisms
$$ L \eta_\mu R\Gamma_{\textup{ct}} (\Delta, \Ainf(R_\infty)) \isom  L \eta_\mu R\Gamma_{\pet} (U, \AAinf) \isom R \Gamma(\fU, A\Omega_{\fX}).$$ In particular, the canonical map $A \Omega^{\textup{psh}}_{\fX} \isom A \Omega^{\textup{pre}}_{\fX}$ is an isomorphism of presheaves in $\textup{Psh}(\fX_{\ett, \textup{aff}}, \mD(\Ainf))$. 
\et

\br For simplicity, we denote the asserted quasi-isomorphisms by 
$$\sq \Omega_R^{\textup{gp}} \isom \sq \Omega_R^{\textup{psh}} \isom \sq \Omega_R, \quad \text{and } \:\:  A\Omega_R^{\textup{gp}} \isom  A\Omega_R^{\textup{psh}} \isom A\Omega_R. $$
(See also Subsection \ref{ss:local_analysis_sq_A_Omega}). The proof of $\sq \Omega_R^{\textup{gp}} \isom \sq \Omega_R^{\textup{psh}} \isom \sq \Omega_R$ is similar to \cite{BMS}, while the proof of $ A\Omega_R^{\textup{gp}} \isom  A\Omega_R^{\textup{psh}}$ uses Lemma 3.18 of \cite{Kestutis},\footnote{This is a variant of Lemma 8.11 and Proposition 9.12 of \cite{BMS}, which states that if a map $f: B \ra B'$ in $D(\Ainf)$ is an almost quasi-isomorphism in the sense that $W(\fm^\flat)$ kills $\textup{Cone}(f)$, and if $H^i(B \otimes^{\L} \Ainf/\mu)$ is $W(\fm^\flat)$-torsion free, then $L\eta_\mu$ turns $f$ into an actual quasi-isomorphism.} which reduces to the computation of $W(\fm^\flat)$-torsion of $H^i_{\textup{ct}}(\Delta, \Ainf(R_\infty)/\mu)$. This computation is carried out  in 3.15 -- 3.19 of \textit{loc.cit.} in the case when $s = 1$ in the product $R^\square = R_1^\square \widehat\otimes \cdots \widehat \otimes R_s^\square$, but the same argument goes through in our setting. Finally, the proof of the isomorphism  $A \Omega^{\textup{psh}}_{\fX} \isom A \Omega^{\textup{pre}}_{\fX}$  uses the following claim:
$$\textup{Claim: \: the natural map $ A \Omega_{\fX} \otimes^{\mathbb L}_{\sq \theta} \mO_C \isom \sq \Omega_{\fX}$ is a quasi-isomorphism.}$$  
The claim involves commuting $L \eta_\mu$ with certain tensor products, and relies on
\bl[\cite{Bhatt} 5.16]  \label{lemma:Bhargav_5_16}
Suppose $f, g$ are nonzerodivisors in a ring $A$, and $K \in D(A)$ satisfies that each $H^i(K \otimes^\L A/f)$ is $g$-torsion free, then 
$$ (L \eta_f K )\otimes^\L A/g \isom L\eta_{\cl f} (K \otimes^\L A/g)$$
is a quasi-isomorphism, where $\cl f \in A/g$ is the image of $f$. 
\el 
\noindent To justify the claim it suffices 
to show that for each object $\fU= \spf R \in \fX_{\ett, \textup{aff}}$, we have $A \Omega_R^{\textup{psh}} \otimes^\L_{\sq \theta} \mO_C \isom  \sq \Omega_R^{\textup{psh}}$. By what we have already outlined and Lemma \ref{lemma:Bhargav_5_16}, this amounts to showing that $H^i_{\textup{ct}}(\Delta, \Ainf(R_\infty)/\mu)$ is $p$-torsion free, which again follows from the same argument of Proposition 3.19 in \cite{Kestutis}. Returning to the proof of the quasi-isomorphism $A \Omega_R^{\textup{psh}} \isom A \Omega_R$, since both sides are derived $\sq \xi$-adically complete (see 4.1 and 4.6 of \cite{Kestutis}), and $\Ainf$ is moreover $\sq \xi$-torsion free, it suffices to check this after reducing mod $\sq \xi$, in other words, it suffices to show that $A \Omega_R^{\textup{psh}} \otimes^\L_{\sq \theta} \mO_C \isom A \Omega_R \otimes^\L_{\sq \theta} \mO_C$, by the previous discussion and the projection formula, this reduces to $\sq \Omega_R^{\textup{psh}} \isom \sq \Omega_R$. 
\er 

\subsubsection{The Hodge-Tate (and de Rham) specialization}  \label{ss:dR_comparison}

\bt[\cite{Kestutis}, 4.11] \label{thm:HT_dR_specialization}  \label{thm:dR_comparison} 
There is a (non-canonical) isomorphism of cdga's 
$$H^*(\sq \Omega_{\fX}) \isom \omega^*_{\fX/\mO_C},$$ where the differential on $H^* (\sq \Omega_{\fX}) = H^* ( A \Omega_{\fX} \otimes^\L \Ainf/\sq \xi)$ is the Bockstein differential.\footnote{To get a canonical isomorphism, one needs $\omega_{\fX/\mO}^i \{-i\}$ instead on the right hand side, where $\{-1\}$ denotes the Breuil-Kisin twist. }
Consequently, we have a natural quasi-isomorphism 
$$ A \Omega_{\fX} \otimes^\L_{\Ainf, \theta} \mO_C \cong \omega^*_{\ul{\fX}/\ul{\mO_C}}.$$
 
\et 

\br The key is to show that $H^1 (\sq \Omega_{\fX}) \cong \omega^1_{\fX/\mO_C} \{-1\}$. The proof in \cite{Kestutis} (Theorem 4.11) uses formal GAGA and a form of the Grothendieck existence theorems in the non-Noetherian setup to reduce to the case of good reduction (Theorem 8.3 of \cite{BMS}). The requirements in their argument are the following: write $\fU = \spf R \in \fX_{\ett, \textup{aff}}$ as the formal $p$-adic completion of the base change of an affine scheme $\fU_0$ over a discretely valued subring $\mO \subset \mO_C$, then we need 
\be
\item The $p$-adic completion morphism $\fU \ra \fU_0$ of log ringed $\etale$ sites is strict; 
\item The complement $\fU_0 \minus \fU_0^{\textup{sm}}$ of the smooth locus $ \fU_0^{\textup{sm}} \subset \fU_0$ over $\spec \mO$ has codimension at most $2$. 
\ee 
Both are satisfied for $\fX$ in our setup. 
\er 
 
\subsection{A basic example}  \label{ss:example} 
We record a simple example of the Koszul complex that computes $ L \eta_\mu R\Gamma_{\textup{ct}} (\Delta, \Ainf(R_\infty)) $, which both provides a summary for the discussion in Subsection \ref{sss:local_setup_of_R_infty} and is important for our application later. 

\beg  
\label{example:q_dR_Ainf}  In this example let $\varpi = p^q$ for some $q \in \Q_{> 0}$ and
$$R = R^\square = \mO_C \gr{T_0, T_1}/ (T_0 T_1 - \varpi).$$
Therefore $A(R) = A(R^\square) = \Ainf \gr{X_0, X_1}/ (X_0 X_1 - [\varpi^\flat]).$  In this setup, $\Delta \cong \Z_p (1)$ has a generator $\gamma $ which acts on $X_0$ by $[\epsilon]^{-1}$ and on $X_1$ by $[\epsilon]$. After applying $L\eta_\mu$, the  group cohomology $L \eta_\mu R \Gamma_{\textup{ct}} (\Delta, \Ainf(R))$ is computed by the complex
$$\eta_\mu \Big(A(R) \xrightarrow{\gamma - 1} A(R) \cdot e_1 \Big)  = 
A(R)   \xrightarrow{\gamma - 1} ([\epsilon] - 1) \cdot A(R) \cdot e_1 .$$ This is in turn quasi-isomorphic to the free complex 
$$ A(R)=  \Ainf \gr{X_0, X_1}/ (X_0 X_1 - [\varpi^\flat]) \xrightarrow{\frac{\gamma - 1}{[\epsilon] - 1}}  \Ainf \gr{X_0, X_1}/ (X_0 X_1 - [\varpi^\flat]) \cdot e_1, $$
where $e_1$ is our symbol for a dummy basis in degree $1$ (this is denoted by $\textup{dlog} X_1$ in \cite{BMS}). The differential in the last complex is explicitly given by 
\be[1).]
\item $X_0^m \mapsto -([\epsilon]^{-m} + \cdots + [\epsilon]^{-1}) X_0
^m \cdot e_1,$ 
\item $ X_1^m \mapsto (1 + [\epsilon]+ \cdots + [\epsilon]^{m-1}) X_1
^m \cdot e_1$. 
\ee

\noindent $\bullet$ \textbf{base change to $\mO_C$.} Note that the map $\sq \theta: \Ainf(R_\infty) \ra R_\infty$ sends $ X_i \mapsto T_i^{1/p}$, hence the base change $\big(L \eta_\mu R \Gamma_{\textup{ct}} (\Delta, \Ainf(R)) \big) \otimes^\L_{\sq \theta} \mO_C$ is given by 
$$ \mO_C \gr{T_0^{1/p}, T_1^{1/p}}/ (T_0^{1/p} T_1^{1/p} - \varpi^{1/p}) \longrightarrow \mO_C  \gr{T_0^{1/p}, T_1^{1/p}}/ (T_0^{1/p} T_1^{1/p} - \varpi^{1/p}) \cdot e_1 $$
with differentials 
\be[1).]
\item $ T_0^{m/p} \mapsto - (\zeta_p^{-m}+ \cdots + \zeta_p^{-1}) T_0 \cdot e_1 ,$
\item $ T_1^{m/p} \mapsto  (1 + \cdots + \zeta_p^{m-1}) T_1 \cdot e_1 .$
\ee
This complex is quasi-isomorphic to 
$$ \mO_C \gr{T_0, T_1} /(T_0 T_1 - \varpi)  \xrightarrow{\:\: 0 \: \: } \mO_C \gr{T_0, T_1} /(T_0 T_1 - \varpi) \cdot e_1.$$

\noindent $\bullet$ \textbf{base change to $W(k)$.} Similarly, since $\sq \vartheta([\varpi^\flat]) = 0$ the twisted crystalline base change $\big(L \eta_\mu R \Gamma_{\textup{ct}} (\Delta, \Ainf(R)) \big) \otimes^\L W(k)$ is computed by 
$$W(k) \gr{T_0^{1/p}, T_1^{1/p}}/ (T_0^{1/p} T_1^{1/p}) \longrightarrow W(k) \gr{T_0^{1/p}, T_1^{1/p}}/ (T_0^{1/p} T_1^{1/p}) \cdot e_1 $$ 
with differentials 
\be[1).]
\item $ T_0^{m/p} \mapsto -m T_0^{m/p} \cdot e_1$,
\item $ T_1^{m/p} \mapsto m T_1^{m/p} \cdot e_1 $.
\ee
Note that, since each $T_i^{p^{r-1}}$ maps to $\pm p^{r} T_i^{p^{r-1}} \cdot e_1$, we have 
$$T_i^{p^{r - 1}} \in H^0 \Big( L \eta_\mu R \Gamma_{\textup{ct}} \big( \Delta, \Ainf(R) \big) \otimes^\L W(k)/p^r \Big)$$ 
\eeg

We have presented the simplest possible case to simplify notations, though the example above easily generalizes to $R^\square$ with more variables and to the case where $s > 1$. For example, if $R = R^\square = \mO_C \gr{T_0, ..., T_r, T_{r+1}^{\pm 1}, ..., T_d^{\pm 1}}/(T_0 \cdots T_r - \varpi)$, 
$L \eta_\mu R \Gamma_{\textup{ct}} (\Delta, \Ainf(R))$ is computed by the Koszul complex
\begin{multline*} \mathbb{K}_{A(R)} ( \frac{\gamma_1 -1}{[\epsilon] - 1}, ..., \frac{\gamma_{d} -1}{[\epsilon] - 1}) =  
\bigg( A(R)   \xrightarrow{ \frac{\gamma_1 -1}{[\epsilon] - 1}, ..., \frac{\gamma_{d} -1}{[\epsilon] - 1} } A(R) \cdot  e_1\oplus \cdots \oplus e_d  \longrightarrow  \\
\cdots \longrightarrow   A(R) \cdot  \underset{1 \le i_1 < \cdots < i_k \le d}{\oplus} e_{i_1} \wedge \cdots \wedge e_{i_k}  \longrightarrow \cdots  \bigg)
\end{multline*}


\newpage 
\section{The log crystalline specialization} \label{section:log crystalline_specialization}

In this section we prove Theorem \ref{mthm:log_crys_compare_global} in the introduction, which we now restate. Let $\fX$ be a $p$-adic formal scheme of generalized semistable reduction over $\mO_C$, equipped with the divisorial log structure as in Subsection \ref{ss:semistable_setup}. By topological invariance of $\etale$ sites $\fX_{\ett} \cong \fX_{k, \ett}$, we view $\mW \omega^*_{\ul{\fX_k}/\ul k}$ as a sheaf on $\fX_{\ett, \textup{aff}}$. 

\bt \label{thm:log_cris_specialization} 
There exists a canonical, $\varphi$-compatible quasi-isomorphism 
$$\mW \omega^*_{\ul{\fX_k}/\ul k} \isom A\Omega_{\fX} \widehat \otimes^\L W(k)$$ in $\mD(\fX_{\ett}, W(k))$, which induces a quasi-isomorphism 
$$ R\Gamma_{\logcris} (\ul{\fX_k}/W(\ul k)) \cong R\Gamma_{\Ainf} (\fX) \widehat \otimes^\L_{\Ainf} W(k).$$
\et 

The general line of argument is similar to the one used in \cite{BLM}. We first restrict to a smaller site $\fX_{\ett, \textup{aff}}$ and use Theorem \ref{thm:BLM_L_eta_p_fixed_points} to equip the right hand side with a Dieudonn\'e algebra structure; we then appeal to the universal property of log de Rham--Witt complexes to obtain the desired morphism, which is easily checked to be an isomorphism using the Hodge--Tate comparison (cf. Theorem \ref{thm:HT_dR_specialization}). As mentioned in the introduction, aside from the complications caused by allowing the scheme to be more general, we face the additional difficulty of equipping the right hand side with suitable log structures. In order to remedy this, we will work locally on small enough charts and produce the log structures using semistable coordinates, and then show that the morphism on the underlying complexes is canonical (independent of coordinates).  We carry out this procedure in detail in the remaining of the article. 

Recall $A \Omega_{\fX}^{\textup{pre}}$ from Subsection \ref{ss:local_analysis_sq_A_Omega}, which is a commutative algebra object in $ \textup{Psh} (\fX_{\ett, \textup{aff}}, \mD (\Ainf))$. Taking its (completed) base change along $\sq \vartheta: \Ainf \ra W(k)$ we get $A \Omega_{\fX, W}^{\textup{pre}} := A \Omega_{\fX}^{\textup{pre}} \widehat \otimes^\L W(k) \in  \textup{Psh} (\fX_{\ett, \textup{aff}}, \widehat \mD (\Z_p))$, which is a commutative algebra object taking the value $A \Omega_R \widehat \otimes^\L W(k)$ on $\fU = \spf R \in \fX_{\ett, \textup{aff}}$. We also have another presheaf  $(A \Omega_{\fX} \widehat \otimes^\L W(k))^{\textup{pre}} = \iota \big(A \Omega_{\fX} \widehat \otimes^\L W(k) \big)$, where $\iota: \widehat \mD(\fX_{\ett, \textup{aff}}) \ra  \textup{Psh} (\fX_{\ett, \textup{aff}}, \widehat \mD (\Z_p))$ is the fully faithful embedding as described in  Subsection \ref{ss:local_analysis_sq_A_Omega}, which takes the value $R \Gamma(\fU, A \Omega_{\fX} \widehat \otimes^\L W(k))$ on $\fU = \spf R$. Both presheaves sheafify to $A \Omega_{\fX} \widehat \otimes^\L W(k)$ after taking the (derived) $p$-adic completion. In fact, we claim that the following stronger assertion holds: 
\bl \label{lemma:commuting_tensor_with_completed_global_sections}
The natural map $$ A \Omega_{\fX, W}^{\textup{pre}} = A \Omega_{\fX}^{\textup{pre}} \widehat \otimes^\L W(k)  \isom  \big(A \Omega_{\fX} \widehat \otimes^\L W(k) \big)^{\textup{pre}}$$ is an isomorphism. 
\el   

\bproof 
On each $\fU = \spf R \in \fX_{\ett, \textup{aff}}$, we have a canonical map 
\[
\begin{tikzcd}[
  row sep=1.2em,
  ar symbol/.style = {draw=none,"\textstyle#1" description,sloped},
  isomorphic/.style = {ar symbol={\cong}},
  ]
R \Gamma (\fU_{\ett}, A \Omega_{\fX}) \widehat \otimes^\L W(k) \arrow[r]  
& R \Gamma(\fU_{\ett}, A \Omega_{\fX} \widehat \otimes^\L W(k)) \ar[d,isomorphic] \\
&  \textup{Rlim}\:R \Gamma(\fU_{\ett}, A \Omega_{\fX} \otimes^\L W_n(k))
\end{tikzcd}
\]
To prove the lemma, it suffices to show that for each $n$, we have a quasi-isomorphism $R \Gamma (\fU_{\ett}, A \Omega_{\fX})  \otimes^\L W_n(k) \isom R \Gamma(\fU_{\ett}, A \Omega_{\fX}  \otimes^\L W_n(k))$. By Theorem 4.9 of \cite{BMS} we know that each $W_n(k)$ is a filtered colimit of perfect $\Ainf$-modules, since higher direct images of qcqs morphisms commute with filtered colimit (see \cite{SP} 07U6 or more generally 0739), it suffices to show that $A\Omega_R \otimes^\L B \isom R \Gamma(\fU_{\ett}, A \Omega_{\fX}  \otimes^\L B)$ for a perfect $\Ainf$-module $B$. Therefore we are reduced to the case when $B$ is finite projective, in other words, a direct summand of a finite free $\Ainf$-module, the lemma thus follows.  
\eproof 

By the same argument of the lemma, in order to prove Theorem \ref{thm:log_cris_specialization} it is enough to prove the first statement, for which it suffices to construct a quasi-isomorphism $\mW \omega^{*, \textup{pre}}_{\ul{\fX_k}/\ul k} \isom A \Omega_{\fX, W}^{\textup{pre}}$. In the subsequent subsections, we enhance the structure on the right hand side to provide such a map.

\subsection{The Dieudonn\'e algebra structure}

In this subsection we prove
\bp \label{prop:putting_in_correct_category} There exists a canonical isomorphism 
$$ A \Omega_{\fX, W}^{\textup{pre}} \isom L \eta_p  A \Omega_{\fX, W}^{\textup{pre}},$$ 
which puts $  A \Omega_{\fX, W}^{\textup{pre}}  \in  \textup{Psh} (\fX_{\ett, \textup{aff}}, \widehat \mD (\Z_p))^{L\eta_p}.$ Denote by $A_{\fX}^{\textup{pre}}$ the image of $ A \Omega_{\fX, W}^{\textup{pre}}$ in the following sequence of identifications %
$$ \textup{Psh} (\fX_{\ett, \textup{aff}}, \widehat \mD (\Z_p))^{L\eta_p} \cong \textup{Psh} (\fX_{\ett, \textup{aff}}, \widehat \mD (\Z_p)^{L\eta_p}) \cong \textup{Psh} (\fX_{\ett, \textup{aff}}, \textup{DC}_{\str}),$$
then $A_{\fX}^{\textup{pre}}$ takes values in the subcategory of strict Dieudonn\'e algebras. 
\ep

For this we need a sequence of lemmas, starting with 

\bl \label{lemma:Wk_is_colimit}
There exists a canonical isomorphism 
$$ \big(\textup{colim} \Ainf/\varphi^{-r}(\mu) \big)_p^{\widehat{\empty \:\: \empty}} \isom W(k) $$
from the (classical) $p$-adic completion of the filtered colimit of $\Ainf/\varphi^{-r} (\mu)$ along the natural quotient maps to $W(k)$. 
\el

\bproof 
The map $\sq \vartheta$ sends $\mu = [\epsilon] - 1$ to $0 \in W(k)$, therefore we have the canonical map as above. As both sides are $p$-adically complete and $W(k)$ is $p$-torsion free, we can check the isomorphism by reducing mod $p$. Note that taking filtered colimit is exact, so $\textup{colim} \Ainf/\varphi^{-r}(\mu) \cong  \Ainf/ \bigcup\varphi^{-r}(\mu)$, therefore it suffices to show that the canonical map (induced by $\sq \vartheta$) 
$$ \Ainf /\Big(\underset{r}{\cup} \varphi^{-r} (\mu), p \Big) \cong \mO_C^\flat / \underset{r}{\cup} \big([\epsilon]^{\frac{1}{p^r}} - 1\big)  \longrightarrow k$$
is an isomorphism, which in turn follows by considering the ($p^\flat$-adic) valuation of each $[\epsilon]^{1/p^r} - 1 \in \mO_C^\flat$.  
\eproof 

Also recall that the d\'ecalage operators behaves well under completions. 
\bl[\cite{BMS}  6.20] \label{lemma:BMS_6_20}
For any object $C \in D(W(k))$, the natural maps 
$$ \widehat{L \eta_p C} \isom L \eta_p \widehat C \isom \textup{R}\varprojlim L\eta_p (C \otimes^\L W_n(k)) $$ are quasi-isomorphisms, where all completions are derived $p$-adic completions. 
\el
 
\bc  \label{lemma:eta_p_fixed_point_after_basechange}
For any $\spf R \in \fX_{\ett, \textup{aff}}$, there is a natural Frobenius compatible quasi-isomorphism: 
$$ \big(L\eta_{\sq \xi} A \Omega_R \big) \widehat \otimes^\L W(k) \cong L\eta_{p} \big( A \Omega_R  \widehat \otimes^\L  W(k) \big).$$
Therefore, we get an isomorphism of presheaves (valued in $\widehat \mD(\Z_p)$):
$$\big( L \eta_{\sq \xi} A \Omega_{\fX}^{\textup{pre}} \big) \widehat \otimes^\L W(k) \isom L\eta_{p} \big( A \Omega_\fX^{\textup{pre}}  \widehat \otimes^\L  W(k) \big) = L\eta_p A \Omega_{\fX, W}^{\textup{pre}}. $$
\ec 

\bproof 
The proof is similar to that of 10.3.10 of \cite{BLM}. It suffices to show that 
\begin{equation} 
\Big(L \eta_{\sq \xi} A \Omega_R \Big) \otimes^\L \Ainf/\varphi^{-r} (\mu) \isom L \eta_{p} \Big(A \Omega_R  \otimes^\L \Ainf/\varphi^{-r} (\mu)\Big)  \tag{$\star \star$}\end{equation} is a quasi-isomorphism. For this it is enough to show that $H^i (A \Omega_R \otimes^\L \Ainf/\sq \xi)$ is flat over $\mO_C$ (hence $\varphi^{-r}(\mu)$-torsion free) by Lemma \ref{lemma:Bhargav_5_16}.  Now apply Theorem \ref{thm:HT_dR_specialization}, we have $H^i (A \Omega_R \otimes^\L \Ainf/\sq \xi)\cong \omega^i_{\ul R/\ul{\mO_C}}$, which is flat over $\mO_C$, since by assumption $\ul R$ is essentially log smooth over $\ul{\mO_C}$ as discussed in Remark \ref{remark:pullback_log_structure_1}. For completeness we finish the argument by the following sequence of quasi-isomorphisms: 
\begin{align*}
 \big(L\eta_{\sq \xi} A \Omega_R \big) \widehat \otimes^\L W(k)
       & \cong  \big[ \big(L\eta_{\sq \xi} A \Omega_R \big)  \otimes^\L  \varinjlim \Ainf/\varphi^{-r} (\mu) \big]^{\widehat{\:\:\:}}_p
       &\text{by Lemma } \ref{lemma:Wk_is_colimit}  \\
            & \cong   \big[ \varinjlim \Big(L\eta_{\sq \xi } A \Omega_R \otimes^\L  \Ainf/\varphi^{-r} (\mu)  \Big)  \big]^{\widehat{\:\:\:}}_p        
            & \quad (\star)  \\
                   &  \cong  \big[ \varinjlim L\eta_{p}  \Big(A \Omega_R \otimes^\L  \Ainf/\varphi^{-r} (\mu)  \Big)  \big]^{\widehat{\:\:\:}}_p          
                   & \text{by ($\star \star$) above}  \\
                        & \cong   \big[ L\eta_{p}  \Big(A \Omega_R \otimes^\L  \big( \varinjlim\Ainf/\varphi^{-r} (\mu) \big) \Big)  \big]^{\widehat{\:\:\:}}_p  
                        & \quad (\star)    \\                       
                             & \cong   L\eta_{p}  \Big(A \Omega_R \widehat \otimes^\L W(k) \Big)   
                             & \text{by Lemma } \ref{lemma:BMS_6_20}     
\end{align*}
where $(\star)$ uses the fact that both $L\eta$ and tensor commute with filtered colimit. 
\eproof  

This proves the first part of Proposition \ref{prop:putting_in_correct_category}:  
\bc 
There is  a canonical isomorphism 
$ A \Omega_{\fX, W}^{\textup{pre}} \isom \varphi_* L \eta_p  A \Omega_{\fX, W}^{\textup{pre}}$ of commutative algebra objects in $\textup{Psh} (\fX_{\ett, \textup{aff}}, \widehat \mD (\Z_p))$. 
\ec 

\bproof 
Recall from Subsection \ref{ss:local_analysis_sq_A_Omega} that on each $\spf R \in \fX_{\ett, \textup{aff}}$ we have the functorial quasi-isomorphism  induced from the Frobenius automorphism on $\AAinf$:
$$ L\eta_\mu R\Gamma(U, \AAinf) \isom \varphi_* L \eta_{\varphi(\mu)} R \Gamma(U, \AAinf) \cong  \varphi_* L \eta_{\sq \xi} L \eta_{\mu} R \Gamma(U, \AAinf),$$  which gives an isomorphism of commutative algebra objects $A \Omega_{\fX}^{\textup{psh}} \isom \varphi_* L \eta_{\sq \xi} A \Omega_{\fX}^{\textup{psh}}$. By the last assertion in Theorem \ref{prop:kestutis_3_20}, this transports to an isomorphism \footnote{Note that this crucially uses the definition of $\fX_{\ett, \textup{aff}}$, which only consists of small enough affine opens. } 
$$A \Omega_{\fX}^{\textup{pre}} \isom \varphi_* L \eta_{\sq \xi} A \Omega_{\fX}^{\textup{pre}}.$$ 
The corollary then follows by taking base change to $W(k)$ and applying Corollary \ref{lemma:eta_p_fixed_point_after_basechange}, as in the following diagram
\[
\begin{tikzcd}[column sep = 0em]
A \Omega_{\fX, W}^{\textup{pre}}  = A \Omega_\fX^{\textup{pre}}  \widehat \otimes^\L  W(k) \arrow[r, "\sim"] \arrow[rd, "\psi"] 
& \quad  \big(\varphi_* L \eta_{\sq \xi} A \Omega_{\fX}^{\textup{pre}} \big) \widehat \otimes^\L W(k)  \arrow[d, "\wr"] \qquad  \\
&  \varphi_* L \eta_p A \Omega_{\fX, W}^{\textup{pre}} =   \varphi_* L\eta_{p} \big( A \Omega_\fX^{\textup{pre}}  \widehat \otimes^\L  W(k) \big). \qquad
\end{tikzcd}
\]
Note that the vertical isomorphism preserves commutative algebra objects, since $L\eta_p$ is lax symmetric monoidal by \cite{BMS} 6.7.  
\eproof 

Ignoring the $W(k)$-linear structure, the isomorphism $ A \Omega_{\fX, W}^{\textup{pre}} \isom L \eta_p  A \Omega_{\fX, W}^{\textup{pre}}$  puts $  A \Omega_{\fX, W}^{\textup{pre}}  \in  \textup{Psh} (\fX_{\ett, \textup{aff}}, \widehat \mD (\Z_p)^{L\eta_p}) \cong  \textup{Psh} (\fX_{\ett, \textup{aff}}, \textup{DC}_{\str})$, where the isomorphism follows from Theorem  \ref{thm:BLM_L_eta_p_fixed_points}, and we denote the image of $ A \Omega_{\fX, W}^{\textup{pre}}$ in the latter category by  $A_{\fX}^{\textup{pre}}$ to emphasize that it takes value in strict Dieudonn\'e complexes. Note that $A\Omega_{\fX, W}^{\textup{pre}}$ (hence $A_{\fX}^{\textup{pre}}$) takes value in the subcategory of commutative algebra objects in $\mD (\Z_p)^{L\eta_p}$ (resp. in $\textup{DC}_{\str}$).   
Next we prove the remaining of Proposition \ref{prop:putting_in_correct_category}. 

\bl
For each $\spf R \in \fX_{\ett, \textup{aff}}$, the commutative algebra object in strict Dieudonn\'e complexes $A_R^* :=  R\Gamma_{\ett}(\spf R, A_{\fX}^{\textup{pre}})$ is a strict Dieudonn\'e algebra.  
\el
 
\bproof By Remark 3.1.5 in \cite{BLM}, we need to show that $A_R^*$ is concentrated in degree $\ge 0$ and that $F$ is the $p$-power Frobenius on $W_1 (A_R^0)=A_R^0/V(A_R^0)$. Recall from Section \ref{sec:log_dRW} Remark \ref{remark:BLM_L_eta_p_fixed_points} that we have $$A_R^* = \varprojlim_{R_r} H^*( A \Omega_R  \otimes^\L W(k)/p^r).$$  For the first assertion, it suffices to show that $$A \Omega_R  \otimes^\L W(k)/p = \sq \Omega_R \otimes^\L_{\mO_C} \mO_C/\fm $$ is coconnective (therefore by induction each $A \Omega_R  \otimes^\L W(k)/p^r$ is coconnective). This is follows from Theorem \ref{thm:dR_comparison} (note that since $H^0 (\sq \Omega_R) \cong R$ is $\mO_C$-flat, we have $H^0 (\sq \Omega_R) \otimes k = H^0 (\sq \Omega_R \otimes^\L k)$). 

Now it remains to show that the morphism induced by $F$ on $H^0(\sq \Omega_R \otimes^\L_{\mO_C} k)$ is the $p$-power Frobenius. The proof given here mimics the proof of 10.3.15 of \cite{BLM} with a minor difference in the last step.  We first choose a semistale coordinate $\square: R^\square \ra R$ and (hence) a formally $\etale$ morphism of perfectoid algebras $R_\infty^\square \ra R_\infty$ as in Subsection \ref{ss:local_analysis_sq_A_Omega}.  From Proposition \ref{prop:kestutis_3_20} we obtain a $\varphi$-equivariant map $A \Omega_R \ra \Ainf (R_\infty)$, since $A \Omega_R \cong L\eta_\mu R\Gamma_{\textup{ct}} (\Delta, \Ainf(R_\infty))$.  Note that since $R_\infty$ is perfectoid we have a $\varphi$-equivariant identification $\Ainf(R_\infty) \otimes^\L W(k) = W(R_{\infty, k})$, where $\varphi$ equals to the Witt-vector Frobenius on $W(R_{\infty, k})$, by Lemma 3.13 of \cite{BMS}. Therefore, after base change to $k = \mO_C/\fm$ we get a morphism $\sq \Omega_R \otimes^\L_{\mO_C} k \ra R_{\infty, k}$. Now taking $H^0$ we get the following $\varphi$-equivariant morphism
$$R_k = H^0( \sq \Omega_R \otimes^\L_{\mO_C} k) \longrightarrow R_{\infty, k}$$ where $\varphi$ agrees with $F$ on the left and identities with the usual Frobenius on the right. It suffices to show that the map $R_k \ra R_{\infty, k}$ is injective. For this we observe that the following diagram 
\[
\begin{tikzcd}[column sep=1.5em, row sep=1.5em]  R_{k}^\square \arrow[r] \arrow[d] & R^\square_{\infty, k} \arrow[d] \\ R_{k}  \arrow[r]  & R_{\infty, k}  \end{tikzcd}
\]
is cocartesian. Since $R_k^\square \ra R_k$ is $\etale$, and $R_k^\square \ra R_{\infty, k}^\square$ is injective, the proposition therefore follows. Note that unlike the situation in \cite{BLM}, here $R/p \ra R_\infty/p$ is not flat (at the singular point of the special fiber). 
\eproof

\subsection{A log Dieudonn\'e algebra structure (after fixing coordinates)} \label{ss:log_structure_after_coordinates}

Now we have two presheaves $\mW \omega^{*, \textup{pre}}_{\ul{\fX_k}/\ul k}$  and $A_{\fX}^{\textup{pre}}$ on $\fX_{\ett, \textup{aff}}$, both valued in the category $\textup{DA}_{\str}$ of strict Dieudonn\'e algebras. In order to obtain a functorial map between them, we want to enhance the right hand side with log structures. As mentioned in the beginning of this section, we will do this locally on small enough charts.  For this subsection we let $\spf R \in \fX_{\ett, \textup{aff}}$ and fix a choice 
$$\square: R^\square \longrightarrow R$$ of a semistable coordinate as in Subsection \ref{ss:semistable_setup}. Our goal is to equip $A_R^*$ with a suitable log structure and log derivation using this coordinate and the local analysis in Subsection \ref{ss:local_analysis_sq_A_Omega} and \ref{ss:example}, which upgrades $A_R^*$ to a log Dieudonn\'e algebra.   

\subsubsection{The log structure}   Recall the setup from Subsection \ref{ss:semistable_setup}, where  we have the log algebra $(R, M^\square)$ over $(\mO_C, N^\square)$ where $N^\square = \N^s$.  Let $A(R^\square) \ra A(R)$ be the (formally) $\etale $ morphism given by $\square: R^\square \ra R$ as in Subsection \ref{sss:local_setup_of_R_infty}, this determines the vertical map in the following diagram
\[
\begin{tikzcd}[column sep=1em, row sep=2.5em]
H^0(L\eta_\mu R\Gamma_{\textup{ct}} (\Delta, A(R^\square))   \otimes^\L W_r) \arrow[d] \arrow[r, equal, "\ref{example:q_dR_Ainf}"] \arrow[rd, "\kappa_r"] & H^0 \big( \Big(\underset{1\le h \le s}{\otimes} \mathbb{K}_{A(R^\square)} (\frac{\gamma_{h, i} - 1}{[\epsilon] - 1}) \Big)   \otimes W_r \big) \\ 
H^0(L\eta_\mu R\Gamma_{\textup{ct}} (\Delta, A(R))  \otimes^\L W_r) \arrow[r] &  H^0(A \Omega_R   \otimes^\L W_r) 
\end{tikzcd}
\]
where we have abbreviated notation by using $W_r$ for $W_r(k)$ and using $\mathbb{K}_{A(R^\square)} (\frac{\gamma_{h, i} - 1}{[\epsilon] - 1})$ to denote the Koszul complex 
$$ \mathbb{K}_{A(R^\square)} (\frac{\gamma_{h, 1} - 1}{[\epsilon] - 1}, ..., \frac{\gamma_{h, d_h} - 1}{[\epsilon] - 1}).$$
The natural composition is denoted by $\kappa_r$.   

\bd \label{def:giving_log_structure} Let $ \mathfrak{b}_{h, i} = (..., 0, 1, 0, ...)\in M^\square = \oplus_h \N^{r_h + 1}$ be the element with $1$ at the $(h, i)$-entry and $0$ otherwise. For each $r \ge 1$, we define the monoid maps  
$$\alpha_r = \alpha_r^\square: M^\square \ra H^0(A \Omega_R \otimes^\L W_r) $$ by sending each $\mathfrak{b_{h, i}}$ (where $1 \le h \le s, 0 \le i \le r_h$) to $\kappa_r (T_{h, i}^{p^{r-1}})$, where $\kappa_r$ is the composition as defined above.  Note that $T_{h, i}^{p^{r-1}} \in H^0(L\eta_\mu R\Gamma_{\textup{ct}} (\Delta, A(R^\square))  \otimes^\L W_r) $ by the discussion in Example \ref{example:q_dR_Ainf}. 
\ed

\br \label{remark:on_the_def_of_alpha} Since $\sq \xi \xmapsto{\sq \vartheta} p \in W(k)$, the map $\sq \vartheta$ induces  $\sq \vartheta: \Ainf/\sq \xi^r \ra W(k)/p^r$. This gives rise to the following natural commutative diagram
\[ 
\begin{tikzcd}  [row sep=2em]
H^0(L\eta_\mu R\Gamma_{\textup{ct}} (\Delta, A(R^\square)) \otimes^\L_{\sq \theta} \Ainf/\sq \xi^r)  \arrow[r] \arrow[d]
& H^0(L\eta_\mu R\Gamma_{\textup{ct}} (\Delta, A(R^\square)) \otimes^\L W_r)  \arrow[d]  \\
H^0(A \Omega_R \otimes^\L_{\sq \theta} \Ainf/\sq \xi^r) \arrow[r] 
 & H^0(A \Omega_R \otimes^\L W_r) 
\end{tikzcd}
\] 
The map $\alpha_r$ can be equivalently described as follows: $\alpha_r (\mathfrak{b_{h, i}})$ is the image of 
$$ T_{h, i}^{p^{r-1}} = \sq \theta (X_{h, i}^{p^r}) \in H^0(L\eta_\mu R\Gamma_{\textup{ct}} (\Delta, A(R^\square)) \otimes^\L_{\sq \theta} \Ainf/\sq \xi^r)$$ in $H^0(A \Omega_R  \otimes^\L W_r)$, where $\sq \theta: \Ainf (R_\infty^{\square, \flat}) \ra R_\infty^\square$ is defined similarly as before. 
\footnote{From Example \ref{example:q_dR_Ainf}, it is straightforward to check that $X_{h, i}^{p^r}$ is indeed an element in $ H^0(L\eta_\mu R\Gamma_{\textup{ct}} (\Delta, A(R^\square)) \widehat\otimes^\L_{\sq \theta} \Ainf/\sq \xi^r)$.} 
\er 

From the next lemma we use the same notation from Remark\ref{remark:BLM_L_eta_p_fixed_points}, it is convenient to also recall the diagram
\[\begin{tikzcd}
 H^i(A \Omega_R \otimes^\L W_r)) \arrow[r, "\mu_r"] \arrow[rd, swap, "R_r"] & H^i(L\eta_{\sq \xi} A \Omega_R  \otimes^{\L} W_{r-1}) \arrow[d, "\psi^{-1}"]  \\
 &   H^i( A \Omega_R \otimes^\L W_{r-1})
\end{tikzcd} \]
that defines the restriction map $R_r$, where $\mu_r$ sends $[y]$ to $[p^i y]$. Note that in the top right entry we have used the identification $ L\eta_{p} ( A \Omega_R \widehat \otimes^\L W ) \cong (L\eta_{\xi} A \Omega_R) \widehat \otimes^\L W  $ as proven in Corollary \ref{lemma:eta_p_fixed_point_after_basechange}. 

\bl  \label{lemma:supplying_the_log_structure}
The maps $\{\alpha_r\}$ is compatible with the restriction maps $R_r$, and satisfies $F(\alpha_r (m)) = \alpha_r (m)^p$. In other words the following diagrams commute 
\[
\begin{tikzcd}
M^\square \arrow[r, "\alpha_r"] \arrow[rd, swap, "\alpha_{r-1}"] & H^0(A \Omega_R \otimes^\L W_r) \arrow[d, "R_r"] \\
& H^0(A \Omega_R \otimes^\L W_{r-1})   
\end{tikzcd} \quad \qquad
\begin{tikzcd}
M^\square \arrow[r, "\alpha_r"] \arrow[d, "\times p"] & H^0(A \Omega_R \otimes^\L W_r) \arrow[d, "F"] \\
M^\square  \arrow[r, "\alpha_{r-1}"] & H^0(A \Omega_R \otimes^\L W_{r-1})   
\end{tikzcd}
\]  
Consequently, we obtain a monoid morphism 
$$\alpha = \alpha^\square: M^\square \ra A_R^0 = \varprojlim_{R_r} H^0(A \Omega_R \otimes^\L W_r)$$ which satisfies $F (\alpha (m)) = \alpha (m)^p$ for all $m \in M^\square$. 
\el

\bproof 
The commutativity of the second diagram follows directly from the definitions. For the commutativity of the first diagram, we reduce to the case of $R^\square$ by functoriality. To simply notation, write 
$$ A \Omega_{R^\square}^{\textup{gp}} := L\eta_\mu R\Gamma_{\textup{ct}} (\Delta, \Ainf(R_\infty^\square))   $$
By Proposition \ref{prop:kestutis_3_20} we know that $ A \Omega_{R^\square}^{\textup{gp}} \isom A \Omega_{R^\square}$ is a quasi-isomorphism, and by the same proof of \ref{lemma:eta_p_fixed_point_after_basechange}, we have a quasi-isomorphism 
$$\psi: A \Omega_{R^\square}^{\textup{gp}}  \widehat\otimes^\L W \longrightarrow L\eta_p \big( A \Omega_{R^\square}^{\textup{gp}} \widehat\otimes^\L W \big).$$  
Now consider the following commutative diagram  
\[
\begin{tikzcd}[column sep = 1.5em]
H^0(A \Omega_{R^\square}^{\textup{gp}} \widehat\otimes^\L W_r)    \arrow[r, "\mu_r"] \arrow[d]
&H^0(L\eta_p\Big(A \Omega_{R^\square}^{\textup{gp}}  \widehat\otimes^\L W\Big)/p^{r-1})  \arrow[r, "\psi^{-1}"] \arrow[d]
&H^0(A \Omega_{R^\square}^{\textup{gp}} \widehat\otimes^\L W_{r-1})  \arrow[d]  \\ 
H^0(A \Omega_{R^\square} \widehat\otimes^\L W_r) \arrow[r, "\mu_r"]
& H^0( L\eta_p \Big(A \Omega_{R^\square} \widehat\otimes^\L W\Big) /p^{r-1})  \arrow[r, "\psi^{-1}"]
& H^0(A \Omega_{R^\square} \widehat\otimes^\L W_{r-1})
\end{tikzcd}
\]
The top right arrow $\psi^{-1}$ comes from the quasi-isomorphism $\psi$ described above, and can be represented by 
$$\phi_F: \big(\otimes_h \mathbb{K}_{\Ainf(R_\infty^\square)} (\frac{\gamma_{h, i} - 1}{[\epsilon] - 1}) \big) \otimes W \longrightarrow \eta_p \Big(\big(\otimes_h \mathbb{K}_{\Ainf(R_\infty^\square)} (\frac{\gamma_{h, i} - 1}{[\epsilon] - 1}) \big) \otimes W \Big),$$
given by $p^k F$ in degree $k$ (cf. Definition \ref{def:usual_DA_2}). In particular, in degree $0$, the map sends $T_{h, i} \xrightarrow{\psi = F} T_{h, i}^p$. Therefore, unwinding definitions, we see that for each $r \ge 2$ the top arrow in the digram above sends 
$$ T_{h, i}^{p^{r-1}} \xmapsto{\: \mu \:}  T_{h, i}^{p^{r-1}} \xmapsto{\psi^{-1}} T_{h, i}^{p^{r-2}}.$$
This proves the commutativity of the triangle, hence the lemma. 
\eproof 

\br
Note that we have to use $ A \Omega_{R^\square}^{\textup{gp}}$ instead of  $L\eta_\mu R\Gamma_{\textup{ct}} (\Delta, A(R^\square))$, for otherwise  $\psi^{-1}$ is not well-defined. However, keeping track of the element $T_{h, i}^{p^{r-1}}$ only involves $L\eta_\mu R\Gamma_{\textup{ct}} (\Delta, A(R^\square))$. 
\er

\bl 
$(A_R^0, M^\square \xrightarrow{\alpha} A_R^0)$ is a log algebra over $W(\ul k^\square)$.  Here $\ul k^\square = (k, N^\square)$ is the log point with $0$ maps to $1$ and everything else maps to $0$.   
\el 

\bproof 
We need to show that for every $h$ such that $1\le h \le s$, we have  
$$\alpha (\mathfrak{b_{h, 0}} + \cdots + \mathfrak{b_{h, r_h}}) = 0.$$
Again it suffices to assume that $R = R^\square$, and suffices to show that for each $r \ge 1$, $\alpha_r (\mathfrak{b_{h, 0}} + \cdots + \mathfrak{b_{h, r_h}}) = \alpha_r (\mathfrak{b_{h, 0}}) \cdots  \alpha_r(\mathfrak{b_{h, r_h}}) = 0$, but $L\eta_\mu R\Gamma_{\textup{ct}} (\Delta, A(R^\square)) \widehat\otimes^\L W$ is computed by  the product (over $h$) of 
$$ \Big(\mathbb{K}_{A(R_h^\square)} (\frac{\gamma_{h, i}}{[\epsilon] - 1}) \Big) \otimes W $$ (cf. Example \ref{example:q_dR_Ainf}), and the degree $0$ term is (a product of)
$$W(k) \gr{T_{h, 0}^{1/p}, ..., T_{h, r_h}^{1/p}, ..., T_{h, d_h}^{\pm 1/p}}/ (T_{h, 0}^{1/p} \cdots T_{h, r_h}^{1/p}).$$
In particular, in $H^0(L\eta_\mu R\Gamma_{\textup{ct}} (\Delta, A(R^\square)) \otimes^\L W_r)$, 
$$ \alpha_r (\mathfrak{b_{h, 0}}) \cdots  \alpha_r(\mathfrak{b_{h, r_h}}) = T_{h, 0}^{p^{r-1}} \cdots T_{h, r_h}^{p^{r-1}} = 0.$$
\eproof 

\subsubsection{The log derivation} 
The log derivation is defined in a similar fashion. We start with the commutative diagram 
\[
\begin{tikzcd}[column sep = 1.5em, row sep = 2.5em]
H^1(L\eta_\mu R\Gamma_{\textup{ct}} (\Delta, A(R^\square))  \otimes^\L W_r) \arrow[d] \arrow[r, equal] \arrow[rd, "\lambda_r"] & H^1\big( \Big(\underset{1\le h \le s}{\otimes} \mathbb{K}_{A(R^\square)} (\frac{\gamma_{h, i} - 1}{[\epsilon] - 1}) \Big)   \otimes W_r \big) \\ 
H^1(L\eta_\mu R\Gamma_{\textup{ct}} (\Delta, A(R)) \otimes^\L W_r) \arrow[r] &  H^1(A \Omega_R  \otimes^\L W_r) 
\end{tikzcd}
\] 
As in Example \ref{example:q_dR_Ainf}, we specify a collection of dummy basis $\{e_{h, i} \}_{\substack{1 \le h \le s \\ 1 \le i \le d_h}}$ in degree $1$ for each 
\begin{multline*} 
\mathbb{K}_{A(R^\square}) (\frac{\gamma_{h, i} - 1}{[\epsilon] - 1}) = \\ A(R^\square) \xrightarrow{\frac{\gamma_{h, 1} - 1}{[\epsilon] - 1}, \cdots, \frac{\gamma_{h, d_h} - 1}{[\epsilon] - 1}} A(R^\square) \cdot e_{h,1} \oplus \cdots \oplus A(R^\square) \cdot e_{h, d_h}  \longrightarrow \cdots  $$
\end{multline*}
such that for each $i \ge 1$,  $X_{h, i}^m \mapsto (1+[\epsilon] + \cdots + [\epsilon]^{m-1}) X_{h, i}^m \cdot e_{h, i}.$ We denote its image in $H^1$ of the (base change of the) Koszul complex again by $e_{h, i}$.

\bd \label{def:giving_log_derivation}
For $r \ge 1$, we define a monoid morphism 
$$\delta_r = \delta_r^\square: M^\square \ra   H^1(A \Omega_R  \otimes^\L W_r)$$
 by sending each $\mathfrak{b}_{h, i}$ for $1 \le h \le s, 0 \le i \le r_h$ (cf. Definition \ref{def:giving_log_structure}) to  
$$\delta_r (\mathfrak{b}_{h, i}) = \begin{cases}  \gamma_r( e_{h, i}) &  1 \le i \le r_h \\
-  (\gamma_{e_{h, 1}} + \cdots + \gamma_{e_{h, r_h}}) & i = 0 
\end{cases}.$$ 
\ed

\bl \label{lemma:supplying_the_log_derivation}
For every $r > 1$, every $m \in M^\square$, we have 
$$R_r (\delta_r (m)) = F (\delta_r (m)) = \delta_{r-1}(m).$$ In other words, the maps $\{\delta_r\}$ are compatible with $R_r$. This gives rise to a map
$$\delta = \delta^\square: M^\square \ra A_R^1 := \varprojlim_{R_r} H^1(A \Omega_R \otimes^\L W_r)$$ of monoids,
which satisfies $F \delta = \delta$. 
\el

\bproof 
Similar to the proof of Lemma \ref{lemma:supplying_the_log_structure}, we need to show that the map 
\[
\begin{tikzcd}
H^1(A \Omega_{R^\square}^{\textup{gp}} \otimes^\L W_r)    \arrow[r, "\mu_r"]
&H^1(L\eta_p\Big(A \Omega_{R^\square}^{\textup{gp}}  \widehat\otimes^\L W\Big)/p^{r-1})  \arrow[r, "\psi^{-1}"]  
&H^1(A \Omega_{R^\square}^{\textup{gp}} \otimes^\L W_{r-1}) 
\end{tikzcd}
\]
sends $e_{h, i}$ to $e_{h, i}$ for each $h$ and $1 \le i \le r_h$. As before, $A \Omega_{R^\square}^{\textup{gp}}$ can be computed by the complex  $\otimes_h \mathbb{K}_{\Ainf(R_\infty^\square)} (\frac{\gamma_{h, i} - 1}{[\epsilon] - 1})$. Unwinding definitions, the map $\mu_r$ is given by $[x] \mapsto [px]$ for a cocycle $[x] \in H^1(A \Omega_{R^\square}^{\textup{gp}} \otimes^\L W_r)$. On the other hand, the map $\psi$ is given by 
$$\phi_F: \big(\otimes_h \mathbb{K}_{\Ainf(R_\infty^\square)} (\frac{\gamma_{h, i} - 1}{[\epsilon] - 1}) \big) \otimes W \longrightarrow \eta_p \Big(\big(\otimes_h \mathbb{K}_{\Ainf(R_\infty^\square)} (\frac{\gamma_{h, i} - 1}{[\epsilon] - 1}) \big) \otimes W \Big),$$
which is $p \cdot F$ in degree $1$, and sends $T_{h, i}^{m} \cdot e_{h, i} \mapsto p T_{h, i}^{m} \cdot e_{h, i}$. In particular it sends $e_{h, i} \mapsto p \cdot e_{h,i}$, therefore $\psi^{-1} (\mu_r (e_{h, i})) = e_{h, i}$ as desired. 
\eproof

\bl 
The monoid maps $\alpha$ and $\delta$ constructed above satisfy
\bi
\item $ \alpha(m) \delta (m) = d (\alpha (m)) $ for any $m \in M^\square$. 
\item $ \delta (\textup{diag} (n)) = 0 $ for any $n \in N^\square$. 
\item $d \delta(m) = 0$ for $m \in M^\square$. 
\ei
\el

\bproof 
The lemma follows from the corresponding statements for each $\alpha_r$ and $\delta_r$, and the last two claims directly follow from definitions. We need to show that for each $r \ge 1$, each $h$ and $1 \le i \le r_h$, $ \alpha_r (\mathfrak{b}_{h, i}) \delta_r (\mathfrak{b}_{h, i}) = \beta (\alpha_r (\mathfrak{b}_{h, i}))$ where $\beta$ stands for the Bockstein differential.  

Again by functoriality, it suffices to show that the Bockstein differential 
$$H^0 (A \Omega_{R^\square}^{\textup{gp}} \otimes^\L W_r) \xrightarrow{\beta_r} H^1 (A \Omega_{R^\square}^{\textup{gp}} \otimes^\L W_r) $$ 
sends $T_{h, i}^{p^{r-1}} \xmapsto{\:\: \beta_r \: \: } T_{h, i}^{p^{r-1}} \cdot e_{h, i}$. By Example \ref{example:q_dR_Ainf} (more precisely the computation of the base change to $W(k)$), we know that the differential $\beta$ on 
$$\Big(\otimes_h \mathbb{K}_{A(R^\square)} (\frac{\gamma_{h, i} - 1}{[\epsilon] - 1}) \Big) \otimes W(k)$$
sends $T_{h, i}^{p^{r-1}} = (T_{h, i}^{1/p})^{p^{r}}$ to $p^r T_{h, i}^{p^{r-1}}\cdot e_{h, i}$ for every $h$ and every $i \ge 1$. Therefore 
$$\beta_r (T_{h, i}^{p^{r-1}}) = \frac{\beta (T_{h, i}^{p^{r-1}})}{p^r} =  T_{h, i}^{p^{r-1}} \cdot e_{h, i}.$$
\eproof  

Let us summarize what we have shown so far in the following proposition. 
\bp \label{prop:supplying_the_log_structure_and_derivation} Upon choosing a semistable coordinate $\square: R \ra R^\square$, the strict Dieudonn\'e algebra $A_R^*$ equipped with $M^\square \xrightarrow{\alpha} A_R^0$ and $M^\square \xrightarrow{\delta} A_R^1$ is a strict $p$-compatible log Dieudonn\'e algebra over $W(\ul k^\square)$. We denote this log Dieudonn\'e algebra by $A_R^{*,\square}$ (the log data depend on coordinates). 
\ep


\subsection{The local comparison} Now we are ready to prove a local version of log crystalline comparison (upon fixing coordinates). Recall that $\sq \Omega_R \cong A \Omega_{R} \widehat \otimes^\L_{\sq \theta} \mO_C$. We first observe that there is a natural commutative diagram (the right hand side being a commutative square of rings)
\[
\begin{tikzcd}
M^\square \arrow[r] \arrow[dr, "\alpha_1"] & R^\square  \arrow[r] \arrow[d] 
&  R \arrow[d]  \\
& H^0 (\sq \Omega_{R^\square})   \arrow[r]
& H^0 (\sq \Omega_R) 
\end{tikzcd}
\]
where the monoid morphism $\alpha_1$ comes from Remark \ref{remark:on_the_def_of_alpha} (and gives rise to the map $\alpha_1: M^\square \ra H^0(A \Omega_{R^\square}) \otimes^\L k$ via base change), while the map $R^\square \isom H^0 (\sq \Omega_{R^\square})$ (resp. $R \isom H^0 (\sq \Omega_R)$) is the isomorphism given by part (2) of Theorem \ref{thm:HT_dR_specialization}, which sends each $T_{h, i}$ to $T_{h, i}$.  Since $H^0 (\sq \Omega_R) \otimes k = H^0 (\sq \Omega_R \otimes^\L k)$, the base change of the diagram above induces an isomorphism of log algebras over $\ul k^\square = (k, N)$:
$$\tau_1^{\square, 0}:  (R_k, M^\square) \isom (H^0(A \Omega_{R} \otimes^\L k), M^\square).$$
By the construction of $A_R^*$, we have $W_1(A_R^0) = H^0(A \Omega_{R} \otimes^\L k)$, therefore this isomorphism induces a map $$\tau^\square: \mW \omega^*_{(R_k, M^\square)/\ul{k}^\square} \ra A_R^{*, \square}$$
of strict $p$-compatible log Dieudonn\'e algebras over $W(\ul k^\square)$.

\bl  \label{lemma:tau_is_iso}
$\tau^\square$ is an isomorphism of log Dieudonn\'e algebras. 
\el 

\bproof  By Corollary \ref{cor:Cartier_implies_qi}, it suffices to show that the dotted arrow in the diagram 
\[
\begin{tikzcd} 
\mW_1 \omega^*_{(R_k, M^\square) /\ul{k}^\square} \arrow[r] & W_1(A_R^{*,\square}) \arrow[d, "\wr"]   \\
\omega^*_{(R_k, M^\square) /\ul{k}^\square}  \arrow[u, "\wr"]  \arrow[r, dashed, "\tau_1^\square"] & \Big(H^*( \sq \Omega_R \otimes^\L  k), \beta\Big)
\end{tikzcd}
\]
is an isomorphism of cdga's, where $W_1(A_R^{*, \square})$ is equipped with the Bockstein differential $\beta$. By Subsection \ref{ss:local_analysis_sq_A_Omega} we  have isomorphisms of graded algebras
$$ H^*( \sq \Omega_R \otimes^\L  k)  \cong H^*(\sq \Omega_R) \otimes k \cong \omega_{\ul R/\ul{\mO_C}}^* \otimes k \cong \omega^*_{\ul{R_k}/\ul k}.$$
Therefore, as a cdga, the target of $\tau_1^{\square}$ can be identified with $\Big(\omega^*_{\ul{R_k}/\ul k}, \beta\Big)$, where the differentials come from the Bockstein differentials on $H^*( \sq \Omega_R \otimes^\L  k)$, associated to $0 \ra \Z/p \xrightarrow{p} \Z/p^2 \ra \Z/p \ra 0$; and is equipped with a log structure coming from $\alpha_1^\square: M^\square \ra H^0( \sq \Omega_R \otimes^\L  k)$ as in Definition \ref{def:giving_log_structure} (and similarly for $\delta_1^\square$). 

Note that $ \omega^0_{(R_k, M^\square) /\ul{k}^\square} = R_k$ and receives the log structure $M^\square \ra R_k$.  Unwinding definitions (of the construction of $\tau^\square$), the morphism $\tau_1^\square$ on the log algebra $(R_k, M^\square)$ is given by $\textup{id}: (R_k, M^\square) \ra (R_k, M^\square)$.  Therefore, to prove the lemma, it suffices to identify the Bockstein differential $\beta$ with the de Rham differential on $\omega^*_{(R_k, M^\square) /\ul{k}^\square}   \cong \omega^*_{\ul{R_k}/\ul k}$. By the commutativity of 
\[ 
\begin{tikzcd}
 \sq \Omega_R =  A\Omega_R /\sq \xi \arrow[r, "\sq \xi"] \arrow[d]
& A \Omega_R/\sq \xi^2  \arrow[r] \arrow[d]
&  A\Omega_R /\sq \xi \arrow[d] \\ 
A\Omega_R  \otimes^\L W/p  \arrow[r, "p"] 
&A\Omega_R \otimes^\L W/p^2 \arrow[r] 
& A\Omega_R   \otimes^\L W/p
\end{tikzcd}
\] 
We are reduced to show that the Bockstein differentials on $H^*(\sq \Omega_R)$ agree with the de Rham differentials via $H^*(\sq \Omega_R) \cong \omega^*_{\ul{R}/\ul{\mO_C}}$, but this is part of Theorem \ref{thm:HT_dR_specialization}. 
\eproof

\br 
On the underlying rings, the map $\tau_1^{\square, 0}$ created in the paragraph above the lemma is given by the canonical isomorphism $R_k \isom H^0 (\sq \Omega_R \otimes^\L k)$, and has nothing to do with the choice of coordinates.  In particular, as a map of cdga's,  the isomorphism $\tau_1^\square$ in the proof above is independent of any choice of coordinates.  
\er 

\subsection{Independence of coordinates}

We freely use notations from Subsection \ref{ss:semistable_setup}. Let be $\fU = \spf R \in \fX_{\ett, \textup{aff}}$ as in previous subsections. Let $M_{\fU}= M_{\fX}|_{\fU} $ be the log structure on $\fU_{\ett}$ restricted from $\fX$, which is equal to the log structure associated to the monoid $M = \Gamma(\spf R, M_{\fX})$.  We write $\ul{R_k}$ for the log algebra $(R_k, M)$.  From Remark \ref{remark:pullback_log_structure_1}, for each choice of coordinates $\square$, we have an isomorphism of \textit{strict Dieudonn\'e algebras}
$$ \lambda^\square:  \mW \omega^*_{(R_k, M^\square) /\ul{k}^\square}  \isom \mW \omega^*_{\ul{R_k}/\ul k}$$
induced by $M^\square \ra M$ on the log structures on $R_k$.

\bp  \label{prop:two_maps_agree_without_log}
If we ignore the logarithmic structures, then the isomorphism 
$$\tau: \mW \omega^*_{\ul{R_k}/\ul k}  \xrightarrow{(\lambda^\square)^{-1}}  \mW \omega^*_{(R_k, M^\square) /\ul{k}^\square}   \xrightarrow{\tau^\square} A_R^*$$ on strict Dieudonn\'e algebras is independent of choice of coordinates. 
\ep

\bproof 
In the proof we use $\lambda^\square$ to identify the underlying Dieudonn\'e algebras of $\mW \omega^*_{(R_k, M^\square) /\ul{k}^\square}$ and of 
$\mW \omega^*_{\ul{R_k}/\ul k}$ for any coordinate $\square$. 
Let $R^\square \ra R$ and ${R^{\square}}' \ra R$ be two choices of coordinates. We need to show that, the morphism of Dieudonn\'e algebras
$\tau^\square: W  \omega^*_{(R_k, M^\square) /\ul{k}^\square}  \isom A_R^{*,\square}$ and $\tau^{\square'}:  W  \omega^*_{(R_k, M^{\square'}) /\ul{k}^{\square'}}  \isom A_R^{*,\square'}$ coincide. Note that for each $n$, both pre-log structures $M^{\square}$ and ${M^\square}'$ on $\spec W_n(R_k)$ factor through the pre-log structure $\alpha: M  \ra R_k \xrightarrow{[\cdot]} W_n(R_k) .$  From Remark \ref{remark:replace_with_finite_level} following the construction of the saturated log de Rham--Witt complex (cf. Lemma \ref{lemma:homomorphism_F_compatible_algebra} and the proof of Theorem \ref{thm:log_dRW_complex_exist}), it suffices to show that, after we identify 
\begin{align*}
& \lambda_1: \varprojlim \omega^1_{W_n(R_k, M^\square) /W_n(\ul{k}^\square)}\: \: \isom \: \varprojlim \omega^1_{W_n(\ul{R_k})/W_n(\ul k)}  \\
 \textup{and } \:\: & \lambda_2: \varprojlim \omega^1_{W_n(R_k, {M^\square}') /W_n({\ul{k}^\square }')}  \isom  \varprojlim \omega^1_{W_n(\ul{R_k})/W_n(\ul k)} \: ,\:\:\:  
\end{align*} 
the two maps $\tau^\square, \tau^{\square'}$ from $ \varprojlim \omega^1_{W_n(\ul{R_k})/W_n(\ul k)} $ to $A_R^1$ induced\footnote{see Lemma \ref{lemma:homomorphism_F_compatible_algebra} for their construction.} respectively from the log structures $M^\square$ and $M^{\square'}$ agree with each other. More precisely, the lemma will follow from the following
 
\textbf{Claim.} Consider the log differential $\textup{dlog}_1: M^\square \longrightarrow \varprojlim \omega^1_{W_n(R_k, M^\square) /W_n(\ul{k}^\square)}$. For every $m \in M^\square$,  let $\textup{dlog} (m):= \lambda_1(\textup{dlog}_1(m)) \in \varprojlim \omega^1_{W_n(\ul{R_k})/W_n(\ul k)} $. Then the following element in $A_R^1$ defined by 
$$\delta = \delta(m) := \tau^\square (\textup{dlog} (m)) - {\tau^\square}'(\textup{dlog} (m))$$  is equal to $0$. 

Indeed, $\omega^1_{W_n(\ul{R_k})/W_n(\ul k)}$ is generated by elements $dx$ for $x \in W_n(R_k)$ and $\textup{dlog} (m)$ for $m \in M^\square$, as all maps in discussion preserve differentials, the claim implies the lemma by Remark \ref{remark:replace_with_finite_level} and  Lemma \ref{lemma:homomorphism_F_compatible_algebra} (and its proof).

It remains to justify the claim. We first observe that for any $m \in M^\square$, both $\tau^\square (\textup{dlog} (m))$ and ${\tau^\square}'(\textup{dlog} (m))$ are fixed by the Frobenius $F$ on $A_R^1$.\footnote{For example, for $\tau^\square$, note that $M$ is generated by $M^\square$ and $R^\times$, so the image of $\textup{dlog}$ is generated by $\textup{dlog} (m_0)$ for $m_0 \in M^\square$ and $dx/x$ for $x \in R^\times$ as a monoid, both fixed by $F$.} Therefore, if we use $\tau$ to identify $\mW \omega^1_{\ul{R_k}/\ul k}  \isom A_R^1$, we have that 
$$\delta \in \ker (\mW \omega^1_{\ul{R_k}/\ul k} \xrightarrow{F-1} \mW \omega^1_{\ul{R_k}/\ul k} ).$$
By Corollary 2.14 of \cite{Lorenzon}, the Frobenius fixed points of $\mW \omega^1_{\ul{R_k}/\ul k}$ are precisely the Hodge--Witt forms $W \omega^1_{\log}$, consisting of sections which are $\etale$ locally sums of sections $\textup{dlog} (m_i)$ for $m_i \in (\mM^a)^{\textup{gp}}$, where $\mM^a$ is the log structure on $\spec R_k$ associated to the constant pre-log structure $M^\square$.  Now from the relations 
$$\alpha(m) \cdot \tau^\square (\textup{dlog} (m)) = \alpha(m) \cdot {\tau^\square}' (\textup{dlog} (m)) = d \alpha(m),$$ we know that $\delta$ satisfies $\alpha (m) \cdot \delta = 0$. The claim, therefore the proposition, follows from the next lemma. 
\eproof 

We retain the same setup from Proposition \ref{prop:two_maps_agree_without_log} in the following lemma.  For convenience we also denote $Y = \spec R_k$ and $\ul Y = (Y, M_Y)$ which the log scheme associated to the log algebra $(R_k, M)$. 
\bl
Let $m \in M^\square$, and $\delta \in W \omega^1_{\log}$ be a section of the Hodge--Witt forms, such that $\alpha(m) \cdot \delta = 0$, then $\delta = 0$. 
\el
\bproof 
We have the following exact sequences of sheaves from \cite{Lorenzon} (2.13 -- 2.14)
\[
\begin{tikzcd}
& p W \omega^1_{\ul Y/\ul k, \log} \arrow[r, hook]   \arrow[d, hook]  &\im (V+\textup{d} V)  \arrow[r] \arrow[d, hook] & K \arrow[d, hook]  \\
0 \arrow[r] & W \omega_{\ul Y/\ul k , \log}^1 \arrow[r] \arrow[d] & \mW \omega^1_{\ul{Y}/\ul k} \arrow[r, "1-F"] \arrow[d]  &  \mW \omega^1_{\ul{Y}/\ul k} \arrow[r] \arrow[d]  & 0 \\
0 \arrow[r] &  \omega_{\ul Y/\ul k, \log}^1 \arrow[r] &  \omega^1_{\ul{Y}/\ul k} \arrow[r, "1-F"] &  \omega^1_{\ul{Y}/\ul k}/ \textup{d} \mO_{Y} \arrow[r] & 0
\end{tikzcd}
\]
where $\textup{d} \mO_Y$ denotes the subsheaf of $\Omega^1_{Y/k}$ generated by $dy$ for $y \in \mO_Y$, and $K$ denotes the kernel of $ \mW \omega^1_{\ul{Y}/\ul k} \ra \omega_{\ul Y/\ul k}^1/\textup{d} \mO_Y$, which contains $\textup{Fil}^1= \im (V + \textup{d} V)$. 
 Note that the induced $V$-filtration (given by $V + \textup{d}V$) on $W \omega_{\ul Y/\ul k, \log}^1$ agrees with the $p$-filtration since $F = \text{id}$. Write $W\omega^1_{\log}$ for the ($\etale$) global sections, for any $\sigma \in W\omega^1_{\log}$, we write $\cl \sigma$ for its image (under reduction mod $p$) in $\omega^1_{\log}$.  

Now for contradiction we assume that $\delta \ne 0$. If $\cl \delta = 0$ then we can write $\delta = p \delta'$ for some $\delta' \in W \omega^1_{\log}$, which still satisfies $\alpha(m) \cdot \delta' = 0$ since $W \omega^1_{\log} \subset \mW \omega^1_{\ul{R_k}/\ul k}$ is $p$-torsion free. Note that $W \omega^1_{\log}$ is $p$-adically separated, so we may assume without loss of generality that $\cl \delta \ne 0$.  Now we have that $\cl \delta \ne 0$ but $\cl \alpha (m) \cdot \cl \delta = 0$, where $\cl \alpha$ denotes the monoid map $M^\square \ra R_k^{\square} \ra R_k$. This is impossible from the explicit description of (the Frobenius fixed points of)
$$\omega^1_{\ul{R_k}/\ul k} \cong \omega^1_{(R_k, M^\square)/\ul k^\square} \cong \omega^1_{(R_k^\square, M^\square)/\ul k^\square} \otimes_{R_k^\square} R_k$$ which is isomorphic to the free-$R_k$ modules with basis 
$$ \Big\{\textup{dlog} (T_{h, i}), \textup{dlog} (T_{h, j}):  0 \le i \le r_h - 1, r_h + 1 \le j \le d_h \Big\}_{1 \le h \le s}$$
for $R^\square$ as described in Subsection \ref{ss:semistable_setup}. Note that in the indices above we omit $r_h$ as  $\textup{dlog} (T_{h, 0}) + \cdots + \textup{dlog} (T_{h, r_h}) = 0$. 
\eproof 

\subsection{The log crystalline comparison} Now we are finally ready to prove the main theorem of this section. 

\bproof[Proof of Theorem \ref{thm:log_cris_specialization}]
For each $\fU = \spf R \in \fX_{\ett, \textup{aff}}$, we write $\tau_R:  \mW \omega^*_{\ul{R_k}/\ul k}  \isom A_R^*$ for the isomorphism obtained from Proposition \ref{prop:two_maps_agree_without_log}. We claim that this supplies a natural map
$$\tau^{\textup{pre}}: \mW \omega^{*, \textup{pre}}_{\ul{\fX_k}/\ul k} \isom A \Omega_{\fX, W}^{\textup{pre}}$$
of presheaves as desired in the beginning of this section. 

\bl 
Let $f: R \ra S$ be a morphism in $\fX_{\ett, \textup{aff}}$, then the following diagram of maps of Dieudonn\'e algebras commutes  
\[
\begin{tikzcd}
\mW \omega^*_{\ul{R_k}/\ul k} \arrow[r, "\tau_R"] \arrow[d, "\mW_f"]  &   A_R^* \arrow[d, "A_f"] \\
\mW \omega^*_{\ul{S_k}/\ul k} \arrow[r, "\tau_S"]  & A_S^* 
\end{tikzcd}
\]
\el 

\bproof 
Note that all maps have already been specified (with the vertical arrows obtained from functoriality). 
We now choose a coordinate $R^\square \ra R$, which also serves as a coordinate for $S$ via $R \xrightarrow{f} S$. This equips $R$ and $S$ with pre-log structures $M^\square$ as in Subsection \ref{ss:semistable_setup}. By Proposition \ref{prop:two_maps_agree_without_log}, $\tau_R$ does not depend on $\square$ and in particular $\tau_R = \tau^\square \circ (\lambda^{\square})^{-1}$ (and likewise for $\tau_S$). Therefore, it suffices to show that the right square in the following diagram commutes (as the left square commutes by functoriality of saturated log de Rham--Witt complexes)
\[
\begin{tikzcd}
\mW \omega^*_{\ul{R_k}/\ul k}  \arrow[d, "\mW_f"]  \arrow[r, "(\lambda^\square)^{-1}"]  & \mW \omega^*_{(R_k, M^\square)/\ul k^\square}   \arrow[r, "\tau_R^\square"] \arrow[d, "\mW_f^\square"]  &   A_R^{*,\square} \arrow[d, "A_f"] \\
\mW \omega^*_{\ul{S_k}/\ul k}  \arrow[r, "(\lambda^\square)^{-1}"]   & \mW \omega^*_{(S_k, M^\square)/\ul k
^\square}  \arrow[r, "\tau_S^\square"]  & A_S^{*,\square}
\end{tikzcd}
\]
Note the maps in the right square are enhanced to maps of log Dieudonn\'e algebras, where the log structures on $A_R^{*,\square}$ and $A_S^{*,\square}$ are constructed in Subsection \ref{ss:log_structure_after_coordinates} from the coordinate $\square$, and the maps on the monoid $M^\square$ are all given by identify. Now, both maps 
$A_f \circ \tau_R^\square$ and $\tau_S^\square \circ W_f^\square$ are morphisms $ \mW \omega^*_{(R_k, M^\square)/\ul k^\square}  \longrightarrow A_S^{*,\square}$ between log Dieudonn\'e algebras, therefore, by the definition (and construction) of the saturated log de Rham--Witt complex, it suffices to show that they correspond to the same map under the bijection 
$$ \Hom_{\ADpc} ( \mW \omega^*_{(R_k, M^\square)/\ul k^\square},  A_S^{*,\square}) \isom \Hom_{\textup{Alg}^{\log}_{\ul{k}}} (\ul{R_k^\square}, A_S^{0, \square}/ \im V).$$
This follows from the commutative diagram 
\[
\begin{tikzcd}
R_k \arrow[d]  \arrow[r,"\sim"]& H^0(\sq \Omega_R \otimes^\L k) \arrow[d]  \\
S_k   \arrow[r, "\sim"]& H^0(\sq \Omega_S \otimes^\L k)
\end{tikzcd}
\]
\eproof 

 Theorem \ref{thm:log_cris_specialization} now follows, as after sheafification, we get an isomorphism of sheaves valued in the derived $\infty$-categories 
 $$\tau: \mW\omega_{\ul{\fX_k}/\ul k}^* \isom A \Omega_{\fX, W} \isom A \Omega_{\fX} \widehat \otimes^\L W(k)$$
in $\mD(\fX_{\ett, \textup{aff}}, W(k)) \cong \mD(\fX_{\ett}, W(k))$. The quasi-isomorphism on the derived global sections follows from the same proof of Lemma \ref{lemma:commuting_tensor_with_completed_global_sections}. 
 \eproof

\newpage 
\appendix

 
\section{Log geometry} \label{sec:log}
In this section we fix notations we use in log geometry, following \cite{Kato}, and is slightly different from notations in  literature.\footnote{For example, we always use $X$ to denote a scheme and $\ul X$ a (pre-) log scheme, which is opposite to the convention adopted by some algebraic geometers.} Most of the material (if not everything) presented in this section is considered standard, but we make some remarks that might only be obvious for experts. 

\subsection{Log schemes (after Fontaine--Illusie--Kato)} 
 A pre-log structure on a scheme $X$ is a morphism  $\alpha: M \ra \mO_X$ of sheaves of monoids on $X_{\ett}$. The pre-log structure $\alpha$ is a log structure if $\alpha: \alpha^{-1}(\mO_X^{\times}) \isom \mO_X^{\times}$  is an isomorphism.  
We will suppress notation and write $(X, M)$ or $\underline{X}$ for a (pre-) log scheme. For a pre-log structure $M$ on $X$, we often denote by $M^a$ (or $\mM$) its associated log structure, which is the push-out of $\alpha^{-1} (\mO_X^\times) \ra M$ and $\alpha^{-1} (\mO_X^\times) \ra \mO_X^\times$ in the category of sheaves of monoids on $X_{\ett}$. We make the same definitions for formal schemes.  A log structure is \textit{fine} if it is coherent and integral. For most of the article we work with fine log schemes, but in the application to $\Ainf$-cohomology we relax the condition of coherence to quasi-coherence. A chart for a log scheme $\underline{X} = (X, M)$ is a morphism $P_X \ra M$ from a constant pre-log structure $P_X$ to $M$ which induces an isomorphism $P_X^a \isom M$. Charts exist $\etale$ locally for fine log schemes. 
A chart for a morphism $f: (X, M) \ra (Y, N)$ between log schemes is a triple $(P_X \xrightarrow{\psi_X} M, Q_Y \xrightarrow{\psi_Y} N, Q \xrightarrow{\beta} P)$ where $\psi_X, \psi_Y$ are charts for $(X, M), (Y, N) $ respectively, and $\beta$ a morphism of monoids making the obvious diagram commute. 

\bl \label{lemma:pushout_agrees}
Let $X$ be a scheme, $L_X \xrightarrow{\alpha} \mO_X$ be a constant pre-log structure with underlying monoid $L$. Let $\mL = (L_X)^a$ denote the log structure associated to $L_X$, and $L^{\textup{sh}} := \Gamma (X, \mL)$ be the global sections of $\mL$, then its associated log structure $\mL^{\textup{sh}}:=(L^{\textup{sh}}_X)^a$ is identified with $\mL$. 
\el

\bproof  
This is formal. We view $L^{\textup{sh}}$ as a constant pre-log structure $\beta: L^{\textup{sh}}_X \ra \mO_X$ by composing the natural map $L^{\textup{sh}}_X \ra \mL$ with the induced log structure $\alpha^a: \mL \ra \mO_X$. We have the following diagram of pre-log structures on $X_{\ett}$
\[
\begin{tikzcd} 
& L^{\textup{sh}}_X \arrow[rd, "i"] \arrow[rr, "\iota_\beta"] & & \mL^{\textup{sh}} \\ 
L_X \arrow[ru, "h"]  \arrow[rr, "\iota_\alpha"]&  & \mL
\end{tikzcd}
\]
where $h$ is given by $\iota_\alpha$ valued on $X$. By the universal property, $\iota_\beta \circ h: L_X \ra \mL^{\textup{sh}}$ induces a morphism $f: \mL \ra \mL^{\textup{sh}}$ of log structures, while $i$ induces  $g: \mL^{\textup{sh}} \ra \mL$. The maps $f$ and $g$ are inverses to each other. To see that $f \circ g = \textup{id}$, note that under the isomorphism $\Hom_{\text{log}} (\mL^{\textup{sh}}, \mL) \cong \Hom_{\text{pre-log}} (L^{\textup{sh}}_X, \mL)$, $g$ corresponds to $i$, therefore, after composing by $f: \mL \ra \mL^{\textup{sh}}$, we know that $f \circ g$ corresponds to $f \circ i$ under the isomorphism  
$\Hom_{\text{log}} (\mL^{\textup{sh}}, \mL^{\textup{sh}}) \cong \Hom_{\text{pre-log}} (L^{\textup{sh}}_X, \mL^{\textup{sh}})$. But $f \circ i = \iota_\beta$ by construction, which corresponds to $\textup{id}: \mL^{\textup{sh}} \ra \mL^{\textup{sh}}$, hence $f \circ g = \textup{id}.$ The other verification is similar. 
\eproof 

\bl \label{lemma:etale_chart} Let $(\spec R, \mM^a)$ be an affine log scheme associated to a log algebra $\ul R = (R, \alpha: M \ra R)$, and $R \ra S$ be an $\etale$ morphism, then $(\spec S, \mM^a|_{\spec S})$ is precisely the log scheme associated to $L \ra S$. In other words, the composition of monoid maps $L \ra R \ra S$ gives rise to a chart for $\mM^a|_{\spec R}$.
\el 

\bproof 
Let $U = \spec R$ and $V = \spec S$. By definition $\mM^a$ is the sheafification of the presheaf pushout of $\alpha^{-1} (\mO_U^\times) \ra L$ and $\alpha^{-1} (\mO_U^\times) \ra \mO_U^\times$, so its restriction $\mM^a|_V$ can be computed by first restricting the presheaf pushout to $V$, and then taking sheafification as presheaves on $V_{\ett}$. Unwinding definitions, the latter is precisely the log structure associated to $L \ra \mO_V$ on $V_{\ett}$. 
\eproof 
 
\bd A morphism $i = (i, \psi): \uX \ra \uY$ between fine log schemes $\uX = (X, \alpha: M \ra \mO_X)$ and $\uY = (Y, \beta: N \ra \mO_Y)$ is a closed immersion if 
\be
\item $i: X \ra Y$ is a closed immersion; and
\item $\psi^a: i^* N \ra M$ is surjective (here $\psi^a$ is induced from $\psi: i^{-1} (N) \ra M$).  
\ee
A closed morphism $i: \uX \ra \uY$  is exact if moreover in (2) above $\psi^a$ is an isomorphism. More generally, a morphism $f: (X, M) \ra (Y, N)$ is exact if the following diagram is Cartesian: 
\[
\begin{tikzcd}[column sep=1.2em,row sep=1.35em]
f^{-1} (N) \arrow{r}{}  \arrow{d}{} & M \arrow{d}{} \\
f^{-1}(N)^{\textup{gp}} \arrow{r}{}  & M^{\text{gp}}
\end{tikzcd}
\]
A log thickening (of order $n$) of a fine log scheme $(T, L)$ is an exact closed morphism 
$$ \iota: (T, L) \ra (T', L') $$
such that $T$ is defined in $T'$ by an ideal $I$ where $I^{n+1} = 0$. 
\ed 

\beg \label{sss:examples} 
We list some examples. 

$\bullet$  \textit{Divisorial log structure.}  Let $D \subset X$ be a divisor, define 
$$M_D(U) = \{f \in \mO_X(U) : f|_{U\minus D} \in \mO_X^{\times}(U\minus D) \}.$$
Equivalently, let $j: U = X \minus D \hookrightarrow X$ be the open inclusion, then this log structure is  $M_D = j_* \mO_U^{\times} \cap \mO_X \hookrightarrow \mO_X$.  

$\bullet$    \textit{Standard log points.} Let $k$ be a field, the standard log point associated to $k$ is $\spec k$ with pre-log structure $\N \ra k$ by sending $0 \mapsto 1$ and everything else to $0$. 

$\bullet$  \textit{Witt log structure.} Let $(X, M)$ be a log scheme on which $p$ is nilpotent. For each integer $r \ge 1$, then $r$-th Witt log scheme $W_r (X, M)$ consists of the underlying scheme $W_r (X)$ and the pre-log structure $W_r(\alpha): M \ra W_r (\mO_X)$ given by $m \mapsto [\alpha (m)]$, here $[\alpha (m)]$ is the $r$-th Teichmuller lift of $\alpha (m)$ . 

\eeg

\subsection{Log differentials}  
 Let  $f = (f, \psi): \uX \ra \uY$ be a morphism between fine log schemes $\uX = (X, \alpha: M \ra \mO_X)$ and $\uY = (Y, \beta: N \ra \mO_Y)$.

\subsubsection{Log derivations} \label{sss:log_derivation}  
 A log derivation of $\uX/\uY$ into a sheaf of modules $D$ is a pair $(d, \delta)$ where $d: \mO_X \ra D$ is a derivation of $X/Y$ into $D$, and $\delta: M \ra D$ is a morphism of sheaf of monoids, such that  
\bi
\item $d (\alpha (m)) = \alpha (m) \delta (m)$ for any local section $m$ of $M$;  
\item $\delta (\psi (n)) = 0 $ for any local section $n$ of $N$. 
\ei

\br \label{remark:pre_log_to_log_derivation}
Let $f: (X, \alpha: M\ra \mO_X) \ra (Y, \beta: N \ra \mO_Y)$ be a morphism of pre-log schemes. We define a pre-log derivation in the same way as above. Then a pre-log derivation $(d, \delta)$ of $(X, M)/(Y, N)$ in to $D$ extends uniquely to a log derivation $(d, \sq \delta)$ of $(X, M^a)/(Y, N^a)$. 
\er

\subsubsection{Log differentials}  \label{sss:log_differentials}
The relative log differential $\omega_{\uX/\uY}^1$ equipped with $(\text{d}, \textup{dlog})$ is the universal log derivation of $\uX/\uY$ representing the functor $D \mapsto \text{Der}^{\log}_{D}$, where $\text{Der}^{\log}_D$ is the set of log derivations of $\uX/\uY$ into the sheaf of $\mO_X$ module $D$.  More concretely, $\omega_{\uX/\uY}^1$ is the quotient of 
$$\Omega_{X/Y}^1 \oplus (\mO_X \otimes_{\Z} M^{\text{gp}})$$
by (the $\mO_X$ module generated by) local sections 
\bi
\item $(d \alpha(m), 0) = (0, \alpha (m) \otimes m)$
\item $(0, 1 \otimes \psi (n))$
\ei
In the construction $\text{dlog} (m) := (0, 1 \otimes m)$. 
 
\br If $f: \ul X \ra \ul Y$ is a morphism of pre-log schemes, then we can make the same definitions. Moreover, we have 
$$\omega^1_{(X, M)/(Y, N)} \isom \omega^1_{(X, M^a)/(Y, N)} \isom \omega^1_{(X, M^a)/(Y, N^a)}.$$
\er 

\subsubsection{Log smooth morphisms}  
A morphism $f = (f, \psi): \uX \ra \uY$ of fine log schemes is log smooth (resp. log $\etale$) if $X \xrightarrow{f} Y$ is locally of finite presentation and for any commutative diagram 
\[
\begin{tikzcd}
(T, L) \arrow{r}{s} \arrow{d}{\iota}& (X, M) \arrow{d}{f} \\
(T', L') \arrow{r}{t} \arrow[dashed]{ur}{g} & (Y, N)
\end{tikzcd}
\]
where $\iota$ is a second order log thickening, $\etale$ locally there exists (resp. there exists a unique) $g: (T', L') \ra (X, M)$ making all diagrams commute. 

As in the classical case, if $f: \uX \ra \uY$ is log smooth, then the log differential $\omega^1_{\uX/\uY}$ is locally free of finite type. 
 
\bp[\cite{Kato} 3.5]
Let $f: \uX \ra \uY$ be a morphism between fine log schemes as above, and $Q \ra N$ a chart for $\uY$. Then $f$ is log smooth if and only if $\etale$ locally on $X$ there exists a chart $(P_X \ra M, Q_Y \ra N, \beta: Q \ra P)$ of $f$ extending the chart for $Y$ such that  
\be 
\item $\ker (\beta^{\text{gp}})$ and $\coker(\beta^{\text{gp}})_{\text{tor}}$ are finite groups of order invertible on $X$;
\item the induced morphism $X \ra Y\times_{\spec (\Z[Q])} \spec \Z[P]$ is $\etale$. 
\ee
Similarly, $f$ is log $\etale$ if and only if a similar condition holds with the torsion part of the kernel $\coker(\beta)_{\text{tor}}$ replaced by $\coker(\beta)$ in $(1)$ above. 
\ep

\subsection{Log-Cartier type}
\subsubsection{Integral morphisms} \label{sss:integral_morphism}  
A log-smooth morphism might fail to be flat. For example, consider $\ul{X} = (\spec \Z [x, y], \N^2)$ with log structure $(1, 0) \mapsto x;  (0, 1)\mapsto y$. The morphism $\ul X \ra \ul X$  given by $x \mapsto x, y \mapsto xy$ is   log-smooth (even log-$\etale$) but not flat. This leads to the following definition: 
\bd A  morphism $f: \ul{X} \ra \ul{Y}$ is \textit{integral} if for any $\ul{Y'} \ra \ul{Y}$ where $\ul{Y'}$ is a fine log scheme, the base change $\ul{X} \times_{\ul{Y}} \ul{Y'}$ is a fine log scheme.  
\ed 
This is equivalent to requiring that $\etale$ locally on $X$ and $Y$, $f$ has charts given by $\beta: Q \ra P$ such that the induced morphism $\Z[Q] \ra \Z[P]$ is flat. 
If $f: \ul{X} \ra \ul{Y}$ is log-smooth and integral, then the underlying morphism $f: X \ra Y$ is flat. 

\subsubsection{log-Cartier type}  \label{sss:log_Cartier_type}

Let $\ul{X} = (X, M)$ be a log scheme  in characteristic $p$, the absolute ($p$-power) Frobenius $F_{\ul{X}}$ is given by the usual absolute Frobenius on $X$ and $M \xrightarrow{\times p } M$. Note that we have implicitly identified $F_{\ul{X}}^{-1} (M) \cong M$ on $X_{\ett}$. 

\bd A morphism $f: \ul{X} \ra \ul{Y}$ over $\F_p$ is of log-Cartier type if $f$ is integral and the relative Frobenius $F_{\ul{X}/\ul{Y}}$ in the diagram below is exact. 
\[
\begin{tikzcd}
\ul{X} \arrow[bend left]{rr}{F_{\ul{X}}} \arrow[dashed]{r}[swap]{F_{\ul{X}/\ul{Y}}} 
& \ul{X}^{(p)} \arrow{r}{h} \arrow{d}{} \arrow[dr, phantom, "\square", very near start] & \ul{X} \arrow{d}{f} \\ 
& \ul{Y} \arrow{r}{F_{\ul{Y}}} & \ul{Y}
\end{tikzcd}
\]
\ed 

The most important feature for a log-smooth morphism of log-Cartier type is that the Cartier isomorphism holds. This will be a key step to relate our log de Rham--Witt complex with the de Rham complex of a log Frobenius lift. 

\bp[\cite{HK} 2.12] \label{prop:HK_Cartier_criterion}
Let $f$ be a log-smooth morphism of log-Cartier type, then there exists a (Cartier) isomorphism 
$$C^{-1}: \omega^k_{\ul{X}^{(p)}/\ul{Y}} \cong \mH^{k} (\omega^*_{\ul{X}/\ul{Y}}),$$ which, on local sections $a \in \mO_X, m_1, ..., m_k \in M$, is given by 
$$ x \; \textup{dlog} (h^* m_1) \wedge \!\cdots \!\wedge \textup{dlog} (h^* m_k)  \longmapsto F_{\ul{X}/\ul{Y}} (x) \; \textup{dlog}(m_1) \wedge \!\cdots \!\wedge \textup{dlog}(m_k). $$
\ep

\newpage

\bibliographystyle{alpha}
\bibliography{dRW}

\end{document}